\documentclass[twoside]{article}
\usepackage{psfig,latexsym,graphicx, epsfig}
\usepackage{amssymb,amsmath,amsthm,amsopn,pifont} 
\usepackage{subfigure}
\numberwithin{table}{section}
 \newcommand\D{{\mathbb D}}

\newcommand\cS{{\mathcal S}}

\newcommand\cC{{\mathcal C}}
\newcommand\cE{{\mathcal E}}

\newcommand\I{{\mathcal I}}
\newcommand\C{{\mathbb C}}
\newcommand\R{{\mathbb R}}

\newcommand\Z{{\mathbb Z}}
\newcommand\Q{{\mathbb Q}}

\newcommand \bzeta{{\boldsymbol\zeta}}

\newcommand \btau{{\boldsymbol\tau}}
\newcommand \blam{{\boldsymbol\lambda}}

\newcommand \bzero{{\boldsymbol 0}}
\newcommand \bphi{{\boldsymbol\varphi}}
\newcommand \bF{{\mathbf F}}
\newcommand \bA{{\mathbf A}}
\newcommand \bc{{\mathbf c}}

\newcommand\bstar{{\textstyle\star}}

\newcommand\bxi{{\boldsymbol\xi}}

\renewcommand\tt{{\theta}}
\newcommand\idt{{\boldsymbol\xi}}

\newcommand\delbold{{\boldsymbol\delta}}
\newcommand\bdel{{\boldsymbol\delta}}

\newcommand\ssm{{\smallsetminus}}

\newcommand\x{{\mathbf x}}

\newcommand\ab{{\mathbf a}}
\newcommand\ub{{\mathbf u}}
\newcommand\y{{\mathbf y}}
\newcommand\m{{\mathbf m}}
\newcommand\z{{\mathbf z}}
\newcommand\vb{{\mathbf v}}

\newcommand\ord{{\rm ord}}
\newcommand\w{{\mathbf w}}
\newcommand\QED {$\quad\square$}
\newcommand\vc{\stackrel{\to}}

\refstepcounter{equation}

\newtheorem{lem}[table]{Lemma}
\newtheorem{theo}[table]{Theorem}

\newtheorem{coro}[table]{Corollary}

\theoremstyle{definition}
\newtheorem{definition}[table]{\indent Definition}
\newtheorem{rem}[table]{\indent Remark}

\newtheorem{ex}[table]{\indent Example}
\refstepcounter{table}

\begin{document}

\title{\bf Cubic Polynomial Maps\\ with Periodic Critical Orbit,\\
Part II: Escape Regions}

\author{\bf Araceli Bonifant\footnote{Partially supported by the Simons 
Foundation.}, Jan Kiwi\footnote{Supported by Research Network on Low 
Dimensional Dynamics PBCT/CONICYT, Chile.} and John Milnor}

\maketitle
\def\IMSmarkvadjust{0 pt}
\def\IMSmarkhadjust{0 pt}
\def\IMSmarkhpadding{0 pt}
\def\IMSpubltext{Published in modified form:}
\def\SBIMSMark#1#2#3{
 \font\SBF=cmss10 at 10 true pt
 \font\SBI=cmssi10 at 10 true pt
 \setbox0=\hbox{\SBF \hbox to \IMSmarkhpadding{\relax}
                Stony Brook IMS Preprint \##1}
 \setbox2=\hbox to \wd0{\hfil \SBI #2}
 \setbox4=\hbox to \wd0{\hfil \SBI #3}
 \setbox6=\hbox to \wd0{\hss
             \vbox{\hsize=\wd0 \parskip=0pt \baselineskip=10 true pt
                   \copy0 \break%
                   \copy2 \break%
                   \copy4 \break}}
 \dimen0=\ht6   \advance\dimen0 by \vsize \advance\dimen0 by 8 true pt
                \advance\dimen0 by -\pagetotal
	        \advance\dimen0 by \IMSmarkvadjust
 \dimen2=\hsize \advance\dimen2 by .25 true in
	        \advance\dimen2 by \IMSmarkhadjust

%
%
  \openin2=publishd.tex
  \ifeof2\setbox0=\hbox to 0pt{}
  \else 
     \setbox0=\hbox to 3.1 true in{
                \vbox to \ht6{\hsize=3 true in \parskip=0pt  \noindent  
                {\SBI \IMSpubltext}\hfil\break
                \input publishd.tex 
                \vfill}}
  \fi
  \closein2
  \ht0=0pt \dp0=0pt
 \ht6=0pt \dp6=0pt
 \setbox8=\vbox to \dimen0{\vfill \hbox to \dimen2{\copy0 \hss \copy6}}
 \ht8=0pt \dp8=0pt \wd8=0pt
 \copy8
 \message{*** Stony Brook IMS Preprint #1, #2. #3 ***}
}
\SBIMSMark{2009/3}{October 2009}{}
\vspace{-1cm}

\begin{abstract}
The parameter space $\cS_p$ for monic centered cubic polynomial maps with
a marked critical point of period $p$ is a smooth affine
algebraic curve whose
genus increases rapidly with $p$. 
  Each $\cS_p$  consists of a compact
connectedness locus together with finitely many escape regions, each of which is
biholomorphic to a punctured disk and is characterized by an essentially
 unique Puiseux series. This note will describe the topology of $\cS_p$,
and of its smooth compactification, in terms of these escape regions.
It concludes with a discussion of the real sub-locus of $\cS_p$.
\end{abstract}

\vspace{.2cm}
\noindent
{\bf Keywords:} cubic polynomials, canonical parametrization, escape regions,
\break
marked grid, Puiseux series.

\vspace{.2cm}
\noindent
  {\bf Mathematics Subject Classification (2000):} 37F10,  30C10, 30D05.

\setcounter{equation}{0}
\setcounter{table}{0}
\section{Introduction}\label{s-in}

This paper is a sequel to \cite{M4}, and will be continued in \cite{BM}.
 We consider cubic maps of the form
\begin{equation*}\label{e-1} 
F(z)~=~F_{a,v}(z)~=~z^3-3a^2z+(2a^3+v)\,, 
\end{equation*}
and study the smooth algebraic curve
$\cS_p$ consisting of all pairs $(a,\,v)\in\C^2$ such that the
\textbf{\textit{~marked critical point~}}
 $+a$ for this map
has period exactly $p\ge 1$. Here $v=F(a)$ is the \textbf{\textit{~marked
critical value}}.~  We will often identify $F$
with the corresponding point  $(a,\,v)\in\cS_p$, and write $a=a_F,~v=v_F$. 
For each critical point of such a map, there is a uniquely defined 
\textbf{\textit{~co-critical point}} which has the same image under $F$.  
The marked critical point $+a$ has co-critical point $-2a$,
 while the \textbf{\textit{~free critical point}}~
 $-a$ has co-critical point $+2a$.
\smallskip

Here is a brief outline. 
Section~\ref{s-can} introduces a convenient local
parametrization of $\cS_p$ which is uniquely defined up to translation.
 Section~\ref{s-esc} describes several preliminary
invariants of escape regions in $\cS_p$, namely the Branner-Hubbard marked
grid, as well as a pseudo-metric on the filled Julia set which is a sharper
invariant, and the kneading sequence which is a weaker invariant.
 It also presents counterexamples to an incorrect 
 statement in \cite{M4}. Section~\ref{s-pui}
describes the complete classification of escape regions
by means of associated Puiseux series.
Section~\ref{s-mon} provides a more detailed study of these Puiseux series,
centering around a theorem of Kiwi which implies  that the
 asymptotic behavior of the differences $F^{\circ j}(a)-a$ as $|a|\to\infty$,
provides a complete invariant. It also presents an effective algorithm
which shows that the asymptotic behavior of $F^{\circ j}(a)-a$,~
is uniquely determined, up to a multiplicative constant, by the marked grid.
Section~\ref{s-can-esc} relates the Puiseux
series to the canonical coordinates of Section~\ref{s-can}. 
Section~\ref{s-Eul}   computes the Euler characteristic 
of the non-singular compactification $~\overline\cS_p\,$,~ and
 Section~\ref{ss-top} provides further information about
 the topology of $\cS_p$ for small $p$.~
 Section~\ref{s-real} describes the subset of real maps in $\cS_p$.\smallskip

\setcounter{table}{0}
\section{Canonical Parametrization of $\cS_p$}\label{s-can}

For most periods $p$, the parameter curve $\cS_p$ is a many times
 punctured (possibly not connected ?)  
surface of high genus. (See Theorem~\ref{t-eul} and
\S\ref{ss-top}.) For example: 

\qquad$\cS_1~~$ has genus zero with one puncture (so that $\cS_1\cong \C$),

\qquad$\cS_2~~$ has genus zero with two punctures,

\qquad$\cS_3~~$ has genus one with 8 punctures,

\qquad$\cS_4~~$ has genus 15 with 20 punctures.

\noindent
Both the genus and the number of punctures grow exponentially with $p$.
At first we had a great
deal of difficulty making pictures in $\cS_p$, since it seemed hard to find
good local parametrizations. Fortunately however,  there is a very
simple procedure which works in all cases. (Compare Aruliah
 and Corless \cite{AC}.)
\smallskip

Let $\cS\subset\C^2$ be an arbitrary smooth curve, and let
$\Phi:U\to\C$ be a holomorphic function $\Phi(z_1,\,z_2)$
 which is defined and without critical points throughout some neighborhood
$U$ of $\cS$, with $\Phi|_\cS$ identically zero.
Then near any point of $\cS$ there is a {\it local parametrization}
$$  t~\mapsto~\vec z\,(t)~\in~\cS $$
which is well defined, up to translation in the $t$-plane, by the Hamiltonian
differential equation
\begin{equation*}
 \frac{dz_1}{dt}~=~\frac{\partial\Phi}{\partial z_2}\,,\qquad \frac{dz_2}{dt}
~=~-\frac{\partial\Phi}{\partial z_1}\,. 
\end{equation*}
Equivalently,  the total
differential of the locally defined function $t$ on $\cS$ is given by
\begin{equation*}\begin{matrix}\label{e-hde} 
 dt~=~~~\displaystyle\frac{dz_1}{\partial\Phi/\partial z_2} &{\rm whenever} &
 \partial\Phi/\partial z_2\ne 0\,, & {\rm and}\\
\noalign{\vspace{.2cm}}
 dt~= -\displaystyle\frac{dz_2}{\partial\Phi/\partial z_1}&{\rm whenever}&
 \partial\Phi/\partial z_1\ne 0\,.
\end{matrix}\end{equation*}
The identity $$d\Phi~=~(\partial\Phi/\partial z_1)\,dz_1+
(\partial\Phi/\partial z_2)\,dz_2~=~0$$
on $\cS$ implies that the last two equations are equivalent whenever both
partial derivatives are non-zero. \textbf{\textit{ Thus the curve $\cS$
has a canonical local parameter $t$, uniquely defined up  to 
translation.}}
\smallskip

\begin{figure}[ht]
\centerline{\psfig{figure=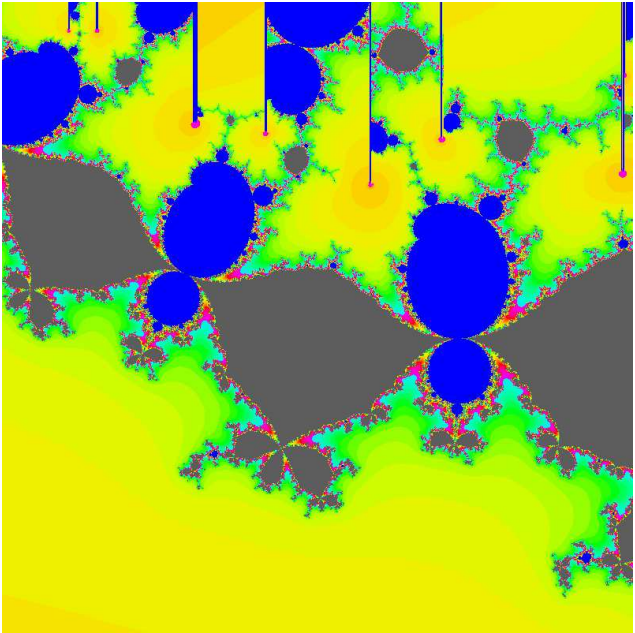,height=3.3in}}
\caption{\textsl{Part of the curve $\cS_4$, represented in the $t$
 parameter plane.}\label{f-S4}}
\end{figure}

In principle, there is a great deal of choice involved here, since we can
multiply $\Phi(z_1,\,z_2)$ by any function $\Psi(z_1,\,z_2)$ which is
holomorphic and non-zero throughout a neighborhood of $\cS$, and thus
obtain many other local parametrizations. However, in the case of the
period $p$ curve $\cS_p$, there is one natural choice which seems convenient.
Namely, using coordinates $(a,\,v)$ as above, we will work with the function
\begin{equation*}\label{e-2} 
\Phi_p(a,\,v)~=~F^{\circ p}(a)-a 
\end{equation*}
which by definition vanishes identically on $\cS_p$, and which
 has no critical points near $\cS_p$. (See \cite[Theorem~5.2]{M4}.)
\smallskip

\begin {rem}\label{r-can-glob}
\rm Although the canonical
 parameter $t$ is locally well
defined up to translation in the $t$-plane, it does not follow that it can
be defined globally. For example,
computations show that there is a loop $L$ in $\cS_3$ with
$$ \int_L dt~\ne~0\,.$$
It follows that $t$ cannot be defined as a single valued function on $\cS_3$.
Conjecturally, the same is true for any $\cS_p$ with $p\ge 3$.
\end{rem}

\setcounter{table}{0}
\section{Escape Regions and Associated Invariants. }\label{s-esc}

By definition, an
 \textbf{\textit{~escape region~}} $\cE_h\subset{\cal S}_p$
is a connected component of the open subset of $\cS_p$ consisting of maps
$F\in\cS_p$  for which the orbit of the free
 critical  point~ $-a$~ escapes to infinity. This section will
 describe some basic invariants for escape regions in $\cS_p$.
The topology of $\cE_h$ can be described as follows. 

\begin{lem}\label{l-multe}
 Each escape region $\cE_h$ is canonically diffeomorphic to the
 \hbox{$\mu_h$-fold} covering of the
complement $\C\ssm\overline\D$, where $\mu_h\ge 1$ is an integer called
the \textbf{\textit{multiplicity\/}} of $\cE_h$.
\end{lem}



\begin{figure}[ht]
\centerline{\psfig{figure=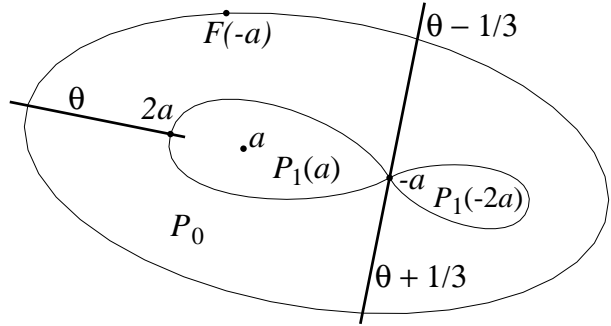,height=1.6in}}
\caption{\textsl{Sketch of the dynamic plane. Here $\tt\in\R/\Z$ is the
co-critical angle.}\label{f1}}
\end{figure}

{\bf Outline Proof.} 
(For details,
 see \cite[Lemma 5.9]{M4}.) The dynamic plane for a map $~F\in\cE_h$
is sketched in Figure \ref{f1}.~
 The equipotential through the escaping critical point $-a=-a_F$
 is a figure eight curve which also passes through the co-critical
point $2a=2a_F$, and which completely encloses the filled Julia set $K(F)$.
The B\"ottcher coordinate ${\mathfrak B}_F(z)\in\C\ssm\overline\D$
is well defined
for every $z$ outside of this figure eight curve, and is also well defined
 at the point $z=2a$. Setting $~{\mathfrak B}(F)\,=\,{\mathfrak B}_F(2a)$,~
we obtain the required covering map
$$ F~\mapsto~ {\mathfrak B}(F)\,\in\,\C\ssm\overline\D\;, $$
and define $~\mu_h\ge 1~$ to be the degree of this map. \qed\medskip

\begin{definition}\label{d-anch}
Any smooth branch of the $~\mu$-th root function
$$ F\mapsto \root{~\scriptstyle\mu} \of {{\mathfrak B}(F)}$$
 yields a bijective map
$~\cE_h\stackrel{\cong}{\longrightarrow} \C\ssm\overline\D$.~
Such a choice, unique up to multiplication by $~\mu$-th roots of unity,
will be called an \textbf{\textit{\,anchoring\/}} 
of the escape region $~\cE_h$.
\end{definition}

\begin{rem}\label{r-Sbar}
 It is useful to compactify $\cS_p$ by adjoining finitely many
 ideal points $\infty_h$, one for each escape region $\cE_h$,
 thus obtaining a compact complex 1-manifold
$\overline\cS_p$. Thus each escape region, together with its ideal
point, is conformally isomorphic to the open unit disk, with parameter
$~1/\root{\scriptstyle\mu}\of{{\mathfrak B}(F)}$. 

 Although $t$ is a local uniformizing parameter
 near any finite point  of $\cS_p$, the $t$-plane is
 ramified over most ideal points.\footnote   
{The behavior of these parametrizations near ideal points of $\cS_p$
will be studied in \S\ref{s-can-esc}.}
A typical picture of the $t$-plane for the curve $\cS_4$ is
 shown in Figure \ref{f-S4}. Here each ideal point in the figure is represented
by a small red dot at the end of a slit in the $t$-plane.
 Each ideal point is surrounded by an escape region which has been colored
 yellow. The various escape regions are separated by the
connectedness locus, which is colored blue\footnote{In grey-scale
  versions 
of these figures, ``blue'' appears dark grey while ``brown'' appears
light grey.}
for copies of the Mandelbrot set
and brown for maps with only one attracting orbit. 
\end{rem}

\begin{rem}\label{r-count} 
The number of escape regions, counted with
multiplicity, is precisely equal to the degree of the curve $\cS_p$.
This grows exponentially
with $p$.
In fact,  $~{\rm deg}(\cS_p)~=~3^{p-1}+O(3^{p/2})~\sim~3^{p-1}~$
 as $~p\to\infty$. (See~\cite[Remark~5.5]{M4}.)
\end{rem}

\begin{rem}\label{r-I}
The curve $~\cS_p~$ has a \textbf{\textit{canonical involution~}}
 $\I$ which sends the map $F:z\mapsto F(z)$ to
 the map $\I(F): z\mapsto -F(-z)$,
rotating the Julia set by $180^\circ$, and also rotating parameter
pictures in the canonical $t$-plane by $180^\circ$.
 In terms of the $(a,\,v)$ coordinates
for $F$, it sends $\,(a,v)\,$ to $(-a,\,-v)$.
This involution rotates some escape regions by $180^\circ$ around the ideal
point, and matches other escape regions in disjoint pairs.
(For tests to distinguish an escape region $\cE_h$ from the
\textbf{\textit{~dual~}} region
$\I(\cE_h)$, see Remark \ref{r-ord-one}.) We will also need the
complex conjugation operation $~F_{a,\,v}\mapsto F_{\overline a,\,\overline v}$,~
or briefly $~F\mapsto \overline F\,$. (Compare \S8.)
\end{rem}


\subsection*{The Branner-Hubbard Puzzle, and the
  Associated Pseudometric, Marked Grid and Kneading Sequence.}
 Branner and Hubbard introduced the structure which is
now called the Branner-Hubbard \textbf{\textit{\,puzzle\,,~}} in order
to study cubic polynomials which are outside of the connectedness
locus. (See \cite{BH}, \cite{Br1}.) They also introduced a
diagram called the \textbf{\textit{~marked grid~}} which captures many
of the essential properties of the puzzle. In this subsection we first
introduce a pseudometric that captures  somewhat more
 of the basic properties of the
 puzzle. This pseudometric determines the associated
marked grid. However, it does not distinguish between a region
$\cE_h$ and its complex conjugate region $\overline\cE_h$ or its dual
region $\I(\cE_h)$.
\smallskip 

For the moment, we do not need the hypothesis
that $F\in\cS_p$. However, we will assume that $F$ is monic and centered,
that its marked critical orbit
$$ a~=~a_0~\mapsto~a_1~\mapsto~a_2~\mapsto~\cdots $$
is bounded, and that the orbit of $-a_0$ is unbounded.

\begin{definition} The
\textbf{~\textit{puzzle piece $P_0$ of level zero~}}
is defined to be the open topological disk consisting
 of all points $z$ in the dynamic plane such that
$$ G_F(z)~<~ G_F\big(F(-a)\big)~=3\,G_F(-a) \,,$$
where $G_F$ is the Green's function (= potential function) which vanishes
only on the filled Julia set $K_F$. For $\ell>0$,
any connected component of the set
$$ F^{-\ell}(P_0)~=~\{z\in\C~~;~~ G_F(z)~<~ G_F(-a)/3^{\ell-1}\} $$
is called a \textbf{~\textit{puzzle piece of level $\ell$}}.
 The notation
$P_\ell(z)$ will be used for the puzzle piece of level $\ell$ which contains
some given point $z\in K_F$. Note that
$$ P_0~=~P_0(z)~\supset~ P_1(z)~\supset~ P_2(z)~\supset~\cdots \,.$$
See Figure~\ref{f1} for a schematic picture of the puzzle pieces
of  levels zero and one.
\end{definition}


\begin{definition}\label{d-pseud}
Given two points $x$ and $y$ in $K_F$, define the \textbf
{\textit{~greatest common level~}} $L(x,y)\in {\mathbb N}\cup\infty$
 to be the largest integer $\ell$ such that
$$ P_\ell(x)~=~P_\ell(y)\,,$$
setting $~L(x,y)=\infty~$ if $~P_\ell(x)=P_\ell(y)~$ for all levels $\ell$.
The \textbf{\textit{~puzzle pseudometric~}} on the filled Julia set $K_F$
is defined by the formula
$$ d_F(x,\,y)~=~2^{-L(x,y)}~\in~[0,\,1]\,.$$
Thus
$$ d_F(x,\,y)~=~
\begin{cases}
0 \qquad \qquad \quad {\rm whenever}\quad P_\ell(x)~=~P_\ell(y)\quad
{\rm for~every~level}~~\ell\,, \\
2^{-\ell}\,>\,0 \qquad {\rm if}~~\ell\ge 0~~{\rm
 is~the~largest~integer~with}~~P_\ell(x)\,=\,P_\ell(y)\,.
\end{cases}  $$
Since the special case $x=a_0$ will play a particularly important role,
we sometimes use the abbreviations
\begin{equation}\label{e-L0}
 L_0(y)~=L(a_0,\,y)\qquad{\rm and}\qquad d_0(y)~=~d_F(a_0,\,y)\,.
\end{equation}
\end{definition}

The basic properties of this pseudometric can be described as follows.

\begin{lem}\label{l-pseud}
For all $~x,\,y,\,z\in K_F$,~ we have:
\smallskip

\begin{itemize}
\item[{\bf(a)}] {\bf Ultrametric inequality.}
$~~~~d_F(x,z)\,\le\,\max\big(d_F(x,y),\;d_F(y,z)\big)$, with
 equality whenever $d_F(x,y)\ne d_F(y,z)$.

\item[{\bf(b)}]~~ $d_F(x,y)\,d_F(x,z)\,d_F(y,z)~<~1\,.$

\item[{\bf(c)}]~~ $d_F\big(F(x),\,F(y)\big)~\le ~2\,d_F(x,\,y)\,,$~ and furthermore

\item[{\bf(d)}]~ if~ $d_F\big(F(x),\,F(y)\big)\,< \,2\,d_F(x,\,y)$
~~ with~~ $d_F(x,y)<1$,~ then
$$d_F(x,y)~=~d_0(x)~=~d_0(y)~>~0\,.$$



\item[{\bf(e)}]~ $d_F(x,y)~=~0$ ~if and only if $x$ and $y$ belong to the same
connected component of $K_F$.
\end{itemize}
\end{lem}
\smallskip

\noindent As an example, applying {\bf(c)} and {\bf(d)} to the case $x=a_0$,
 it follows immediately that
\begin{equation}\label{e-da0}
 d\big(a_1,\,F(y)\big)\,=\,2\,d_F(a_0,\,y)\quad{\rm whenever}\quad
d_F(a_0,\,y)<1\,.
\end{equation}
\smallskip

{\bf Proof of Lemma \ref{l-pseud}.}
 Assertion {\bf (a)} follows immediately from the
definition, and {\bf(b)} is true because there are only two puzzle
pieces of level one.
 The proof of {\bf(c)} and {\bf(d)} will be based on the following
observation.

 Call a puzzle piece \textbf{\textit{~critical~}} if it contains
the critical point $a_0$, and \textbf{\textit{~non-critical~}} otherwise.
Then $F$ maps
each puzzle piece $P_\ell(x)$ of level $\ell>0$ onto the
puzzle piece $P_{\ell-1}\big(F(x)\big)$ by a map which is a two-fold branched
covering if $P_\ell(x)$ is critical, but is a diffeomorphism if $P_\ell(x)$ is
non-critical.
If $x$ and $y$ are
contained in a common piece of level $\ell>0$,
then it follows that $F(x)$ and $F(y)$
are contained in a common piece of level $\ell-1$; which proves {\bf(c)}.
Now suppose that $d_F(x,y)=2^{-\ell}<1$,
and that $d\big(F(x),\,F(y)\big)\le 2^{-\ell}$. Then the two distinct
puzzle pieces $P_{\ell+1}(x)$ and $P_{\ell+1}(y)$, both contained in
$P_\ell(x)$, must map onto a common puzzle piece $P_\ell\big(F(x)\big)$.
Clearly this can happen only if $P_\ell(x)$ also contains the critical point
$a_0$, but neither $P_{\ell+1}(x)$ nor $P_{\ell+1}(y)$ contains $a_0$.
Thus we must have $d_F(x,y)=d_F(a_0,x)=d_F(a_0,y)$, which proves {\bf(d)}.
For the proof of {\bf(e)}, see  \cite[\S5.1]{BH}.
\qed\bigskip

A fundamental result of the Branner-Hubbard theory can be stated
  as follows, using this pseudometric terminology. 

\begin{theo}\label{t-bh}
 Let $K_0\subset
K_F$ be the connected component of the filled Julia set
which contains the critical point $a_0$.
Suppose that $d_F(a_0,\,a_n)=0$ for some
 integer $n\ge 1$. $($In other words,
suppose that $P_\ell(a_0)=P_\ell(a_n)$ for all levels $\ell\,.)$~
If $n$ is minimal,
 then $F^{\circ n}$ restricted to some neighborhood of $K_0$ is polynomial-like,
and is hybrid equivalent to a uniquely defined quadratic polynomial $Q$
with connected Julia set. In this case, countably many
 connected components of $K_F$ are homeomorphic copies of $K_Q$, and all
other components are points.
\end{theo}

{\bf Proof.\/} See
 Theorems 5.2 and 5.3 of \cite{BH}. (In their terminology,
the hypothesis that $d_F(a_0,\,a_n)=0$ is expressed by saying that the
\textbf{\textit{~marked grid~}} of the critical orbit has period $n$. Compare
Definition \ref{d-mg} below.)
 \qed\medskip

By definition, $Q$ will be called the \textbf{\textit~{associated quadratic
map}}.
If $F$ belongs to the period $p$ curve $\cS_p$, then evidently the
hypothesis of Theorem \ref{t-bh} is always satisfied, for some smallest
integer $n\ge 1$ which necessarily divides $p$. It then follows that
the critical point $a_0$ has period $p/n$ under the map $F^{\circ
  n}$. Thus we obtain the following.

\begin{coro}\label{c-bh}
If $F\in\cS_p$ satisfies the hypothesis of Theorem $\ref{t-bh}$,
then the associated quadratic map $Q$ is critically periodic of period
$p/n$.
\end{coro}
\smallskip

\begin{rem}[{\bf Erratum to \cite{M4}}]\label{r-err}
 Unfortunately,
an incorrect version of Corollary \ref{c-bh} was stated
 in  \cite[Theorem 5.15]{M4},
 based on the erroneous belief
than every escape region with kneading sequence of period $n$
must also have a marked grid of period $n$. (See Definitions \ref{d-mg}
and \ref{d-kn}  below.) For counterexamples, see
 Example \ref{ex-p6}.
\end{rem}
\medskip

 The \textit{marked grid} associated with an orbit 
$z_0\mapsto z_1\mapsto\cdots$ in $K_F$ is a graphic method of
visualizing the sequence of numbers
$$ L_0(z_0),~L_0(z_1),~L_0(z_2),~\ldots~\in~{\mathbb N}\cup\infty~,$$
which describe the extent to which this orbit approaches the marked
critical point $a_0$. The marked grid for the critical orbit
$a_0\mapsto a_1\mapsto a_2\mapsto\cdots$
 plays a particularly important role in the 
Branner-Hubbard theory, and
 is the only grid that we will consider.
\smallskip

\begin{definition}\label{d-mg}
{\rm  The
\textbf{\textit{critical marked grid}\/} $M=[M(\ell,\,k)]$,
associated with the critical orbit $a_0\mapsto a_1\mapsto\cdots$
in $K_F$, can be described as
 an infinite matrix of zeros and ones, indexed by
 pairs $(\ell,\,k)$ of non-negative integers.
By definition,  $M(\ell,\,k)=1$ 
if and only if $\ell\le L_0(a_k)$; that is, if and only if
 the puzzle piece $P_\ell(a_k)$ is critical, that is,
if and only if
$$ P_\ell(a_0)~=~P_\ell(a_k)\qquad\Longleftrightarrow\qquad
d_F(a_0,\,a_k)~\le~2^{-\ell}\,.$$
Grid points $(\ell,\,k)$ with $M(\ell,\,k)=1$ are said to be  \textbf{\textit{
 marked}}.
This matrix is represented graphically by an infinite tree where the
 marked points $(\ell,\,k)$ are represented
by heavy dots joined vertically, and joined horizontally along the
entire top line $\ell=0$. This marked grid remains constant
 throughout the escape region.}
\end{definition}

\begin{figure} [ht]
\centerline{\psfig{figure=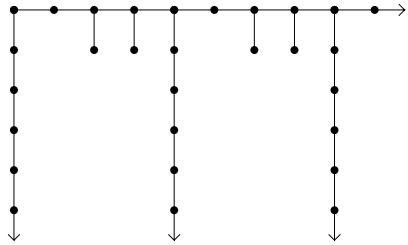,height=1.3in}}
\caption{\label{f-4grid}\textsl{Marked grid for the critical
  orbit of a map belonging to any one of four different
 escape regions in $\cS_4$. For each $k\ge 0$, there are $L_0(a_k)$ vertical
 edges
in the $k$-th column.}}
\end{figure}

\begin{figure}[ht]
\centerline{\psfig{figure=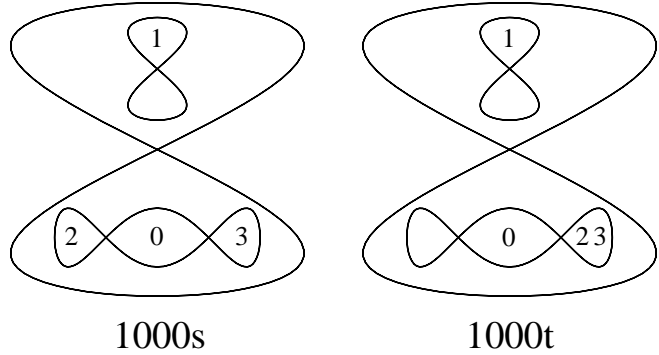,height=1.8in}}
\caption{{\textsl{Schematic picture illustrating the puzzle pieces of
level one and two for the two possible types of map
in $\cS_4$ 
with kneading sequence $\overline{1000}$. In both cases,
 the equipotentials through $-a$ and  through $F^{-1}(-a)$ are shown.
Each point $~a_j=F^{\circ j}(a)~$ of the
critical orbit is labeled briefly as $~j$.~ Thus the points $a_2$ and $a_3$
are separated in the $\overline{\rm 1000}${\rm s} case, but are together in
 the $\overline{\rm 1000}${\rm t} case.}}\label{f1000}}
\end{figure}

As an example, 
Figure~\ref{f-4grid} shows one of the possible critical grids of period 4
for maps belonging to an escape region in $~\cS_4$.~ In fact, there are four
disjoint escape regions which give rise to this same marked grid.
However, up to isomorphism there are only two distinct types, since the
$180^\circ$ rotation $~\I~$ of Remark~\ref{r-I} carries each of these regions
to a disjoint but isomorphic region.
See Figure~\ref{f1000}, which sketches levels one and two of
the puzzles corresponding to the two distinct types. (Compare 
Example~\ref{ex-1000}.)

Since the critical orbit has period 4, the associated pseudometric on the
critical orbit is completely described by the symmetric $4\times 4$ matrix of
distances $[d_F(a_i,\,a_j)]$,~ with $0\le i,\,j<4$. The matrices
corresponding to these two puzzles are given by
$$ \left[\begin{matrix}
0 & 1 & 1/2 & 1/2\\ 1 & 0 & 1 & 1\\
1/2 & 1 & 0 & 1/2\\ 1/2 & 1 & 1/2 & 0
\end{matrix}\right]\qquad{\rm and}\qquad\left[\begin{matrix}
0 & 1 & 1/2 & 1/2\\ 1 & 0 & 1 & 1\\
1/2 & 1 & 0 & 1/4\\ 1/2 & 1 & 1/4 & 0
\end{matrix}\right] $$
respectively. Thus $d_F(a_2,\,a_3)$ is equal to $1/2$ in one case and
$1/4$ in the other. However
 the marked grid, which is a graphic representation of
 the top row of the
matrix, is the same for these two cases.\medskip

 Every such critical
marked grid must satisfy three basic rules, as stated in
\cite{BH}, and also a fourth rule as stated in
\cite{K1}.\footnote{There are similar rules comparing the 
critical marked grid with the marked grid for an arbitrary orbit in
$K_F$. Different versions of a fourth rule have been given by Harris \cite{H}
and by DeMarco and Schiff \cite{DMS}.  (Kiwi and DeMarco-Schiff also consider
 the more general situation where the orbits of both critical points
may be unbounded.)}

\begin{theo}{\bf (The Four Grid Rules)\/}
\begin{itemize}
\item[{\bf (R1)}~] $ 1\,=\,M(0,\,k)~\ge~M(1,\,k))~\ge~M(2,\,k)~\ge~\cdots~\ge~0\,.$
\item[{\bf (R2)}~] If $M(\ell,\,k)=1$, then $M(\ell-i,\,k+i)=M(\ell-i,\,i)$ for
$0\le i\le \ell$.  
\item[{\bf (R3)}~] Suppose that $d_F(a_0,\,a_m)=2^{-\ell}$, that 
$d_F(a_0,\,a_i)>2^{i-\ell}$ for $0<i<k$ and that $~d_F(a_0,\,a_k)<2^{k-\ell}$.
 Then $~d_F(a_0,\,a_{m+k})=2^{k-\ell}$.
\item[{\bf (R4)}~] Suppose that $d_F(a_0,\,a_k)= 2^{-\ell}$, that
 $d_F(a_0,\,a_{k+i})>2^{i-\ell}$ for $0<i<\ell$, and that
 $~d_F(a_0,\,a_\ell)=1$.  Then $~d_F(a_0,\,a_{\ell+k})<1$.
\end{itemize}
\end{theo}
\medskip

{\bf Proof.\/}  All four  rules are consequences of 
Lemma~\ref{l-pseud}. \medskip

\noindent{\bf R1:} The First Rule follows immediately
 from Definition~\ref{d-pseud}.\medskip

\noindent{\bf R2:} For the Second Rule, the statement $M(\ell,\,k)=1$ means that
 $$~ d_F(a_0,\,a_k)~\le~ 2^{-\ell}~.$$
 It then follows
inductively, using Lemma~\ref{l-pseud}-(c), that
$d_F(a_i,\,a_{k+i})\le 2^{i-\ell}$. The ultrametric inequality
then implies that
 $$d_F(a_0,\,a_i)~\le~2^{i-\ell}~\qquad\Longleftrightarrow
\qquad d_F(a_0,\,a_{k+i})~\le~2^{i-\ell}\,,$$
which is equivalent to the required statement.\medskip

\noindent{\bf R3:}
To prove the Third Rule, we will first show by induction on $i$ that
$$  d_F(a_i, \,a_{m+i})~=~2^i\,d_F(a_0, \,a_m)~=~2^{i-\ell} $$
for $0\le i\le k$. The assertion is certainly true for $i=0$. If it is true for
 $i$, then it follows for $i+1$ by Lemma~\ref{l-pseud} items (c) and (d), 
unless
$$ d_F(a_i,\,a_{m+i})~=~d_F(a_0,\,a_i)~=~d_F(a_0,\,a_{m+i})\,.$$
This last equation is impossible since $d_F(a_i,\, a_{m+i})=2^{i-\ell}$ by
 the induction hypothesis but  $d_F(a_0,\,a_i)>2^{i-\ell}$.
In particular, it follows that $d_F(a_k,\,a_{m+k})~=~2^{k-\ell}$. Since
$d_F(a_0,\,a_k)<2^{k-\ell}$, the required equation
 $d_F(a_0,\,a_{m+k})~=~2^{k-\ell}$~
 follows by the ultrametric inequality.
\medskip

\noindent{\bf R4:} Since $d_F(a_0,\,a_k)=2^{-\ell}$, 
it follows inductively that
$d_F(a_i, a_{i+k})=2^{i-\ell}$ for $i\le \ell$. For otherwise
by Lemma~\ref{l-pseud} items (c) and (d), we would have to have
$$ d_F(a_0,\,a_i)~=~d_F(a_0,\,a_{i+k})~=~d_F(a_i,\,a_{i+k}) $$
for some $i<\ell$, which contradicts the hypothesis.
In particular,  this proves that $d_F(a_\ell,\,a_{\ell+k})=1$.
Since $d_F(a_0,\,a_\ell)=1$, it follows from Lemma~\ref{l-pseud}-(b) that
 $d_F(a_0,\,a_{\ell+k})<1$, as asserted.\QED
\medskip
\medskip

\begin{figure}[!ht]
  \makebox[0pt][l]{\includegraphics[width=\textwidth]{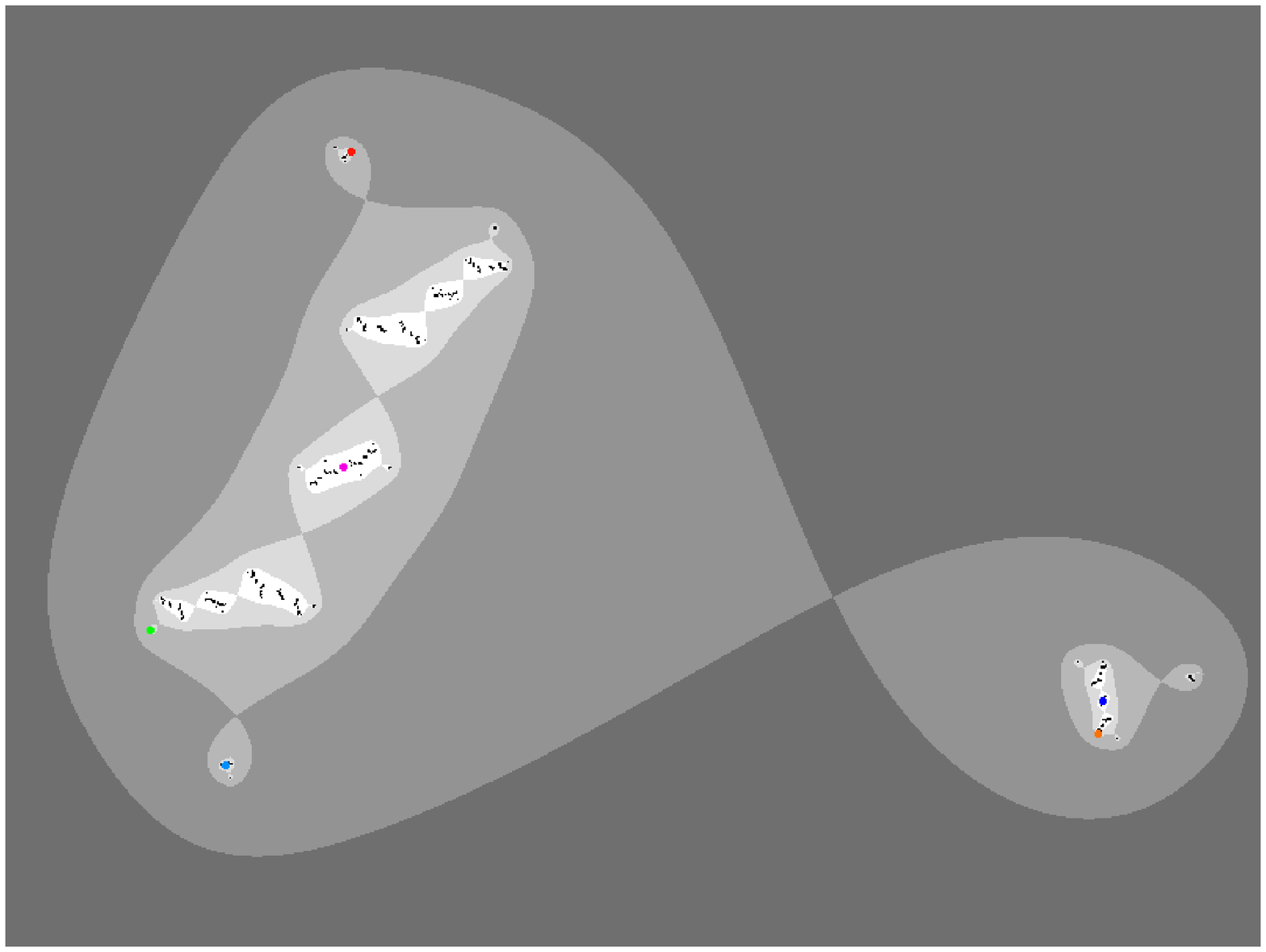}}%
  \includegraphics[width=\textwidth]{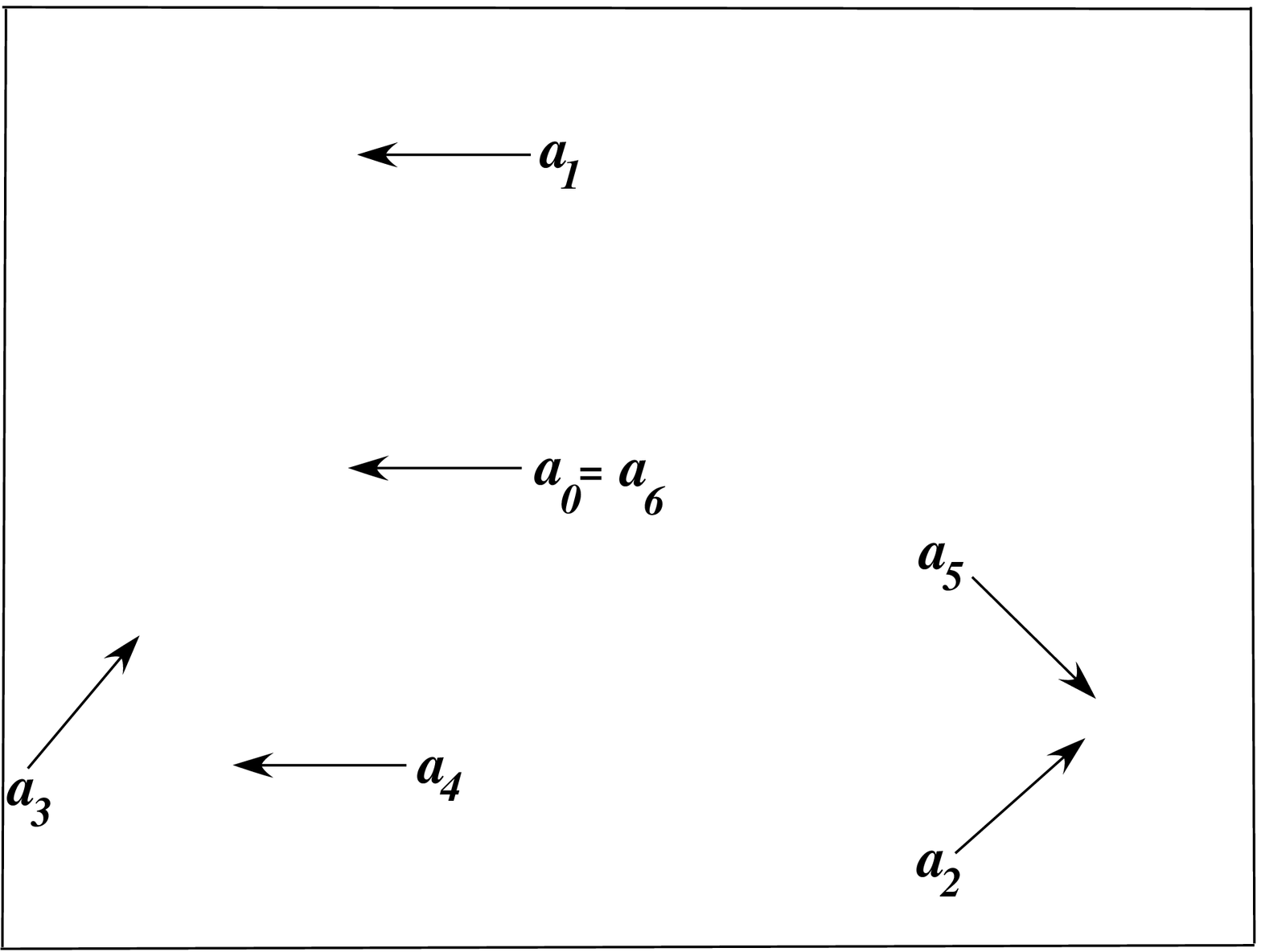}
\caption{\textsl{ Julia set of a polynomial map in $~\cS_6$ with marked grid
of period~$6$ but  kneading sequence $\overline{010010}$ of period~$3$,
showing the puzzle pieces of level one, two and three.}
\label{f-010} }
\end{figure}

For some purposes, it is convenient to work with an even weaker invariant.
 (Compare \cite[\S 5B]{M4}.)

\begin{definition}\label{d-kn}
The \textbf{\textit{~kneading sequence~}} of an orbit $z_0\mapsto z_1\mapsto
\cdots$ in $K_F$ is the sequence
$$\sigma(z_0),\,\sigma(z_1),\,\sigma(z_2),\,\cdots$$
of zeros and ones, where 
$$\sigma(z_j) ~=~\begin{cases} 0 \qquad {\rm if} ~~d_F(z_j,\,a_0)<1~,\\
1 \qquad {\rm if} ~~d_F(z_j,\,a_0)=1 ~.                     
\end{cases}$$
In other words, $~\sigma(z_j)=0~$ if $z_j$ and $a_0$ belong to the same
puzzle piece of level one, but  $~\sigma(z_j)=1~$ if they belong to different
 puzzle pieces of level one.\smallskip

 As an immediate application of the inequalities (b) and (c) of 
Lemma~\ref{l-pseud},  we have the following.
\begin{equation*} 
{\rm If}~ \sigma_{j+n}(z_0)\ne\sigma_{k+n}(z_0)~~ {\rm for~some}~~ n\ge 0\,,~
{\rm then}~~ d_F(z_j,\,z_k)\ge 2^{-n}>0\,. 
\end{equation*}

\begin{figure}[ht!]
\centerline{\makebox[0pt][l]{\scalebox{.902}{\includegraphics[width=\textwidth]{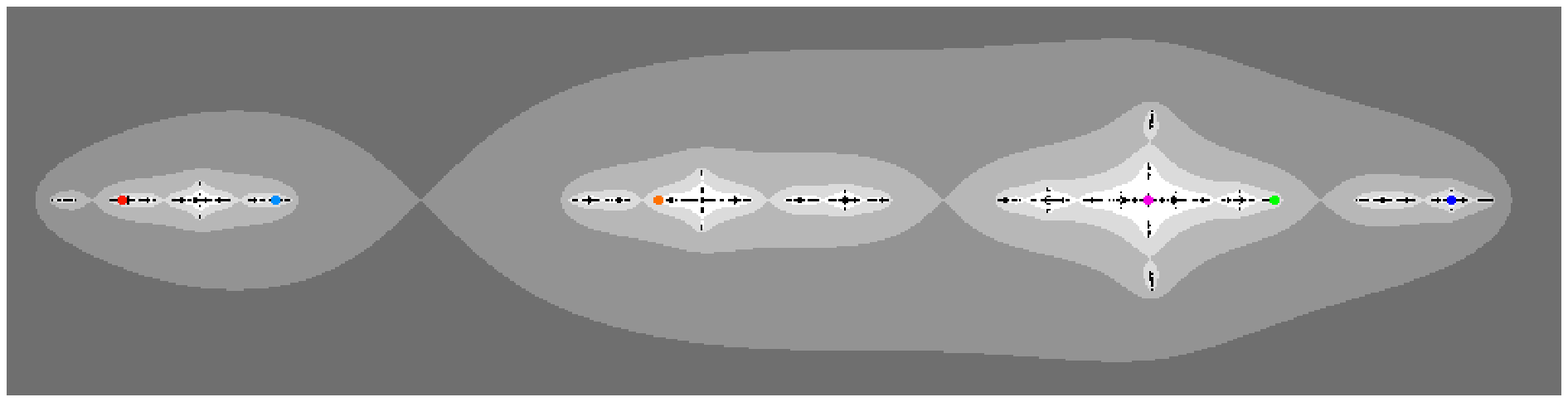}}}%
{\scalebox{.902}{\includegraphics[width=\textwidth]{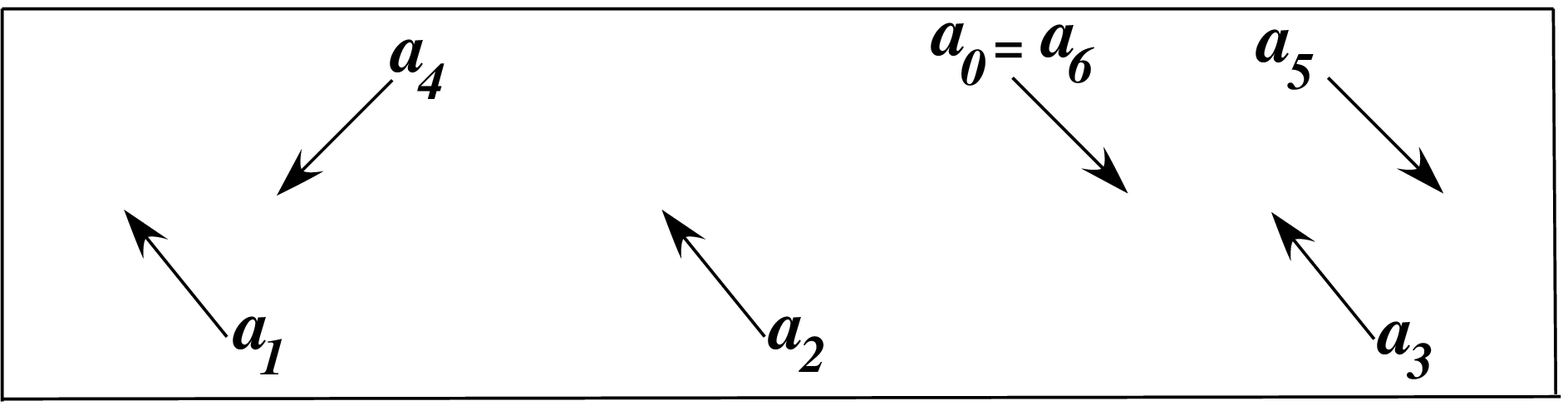}}}}\vskip .4in
\centerline{\scalebox{1.14}{\psfig{figure=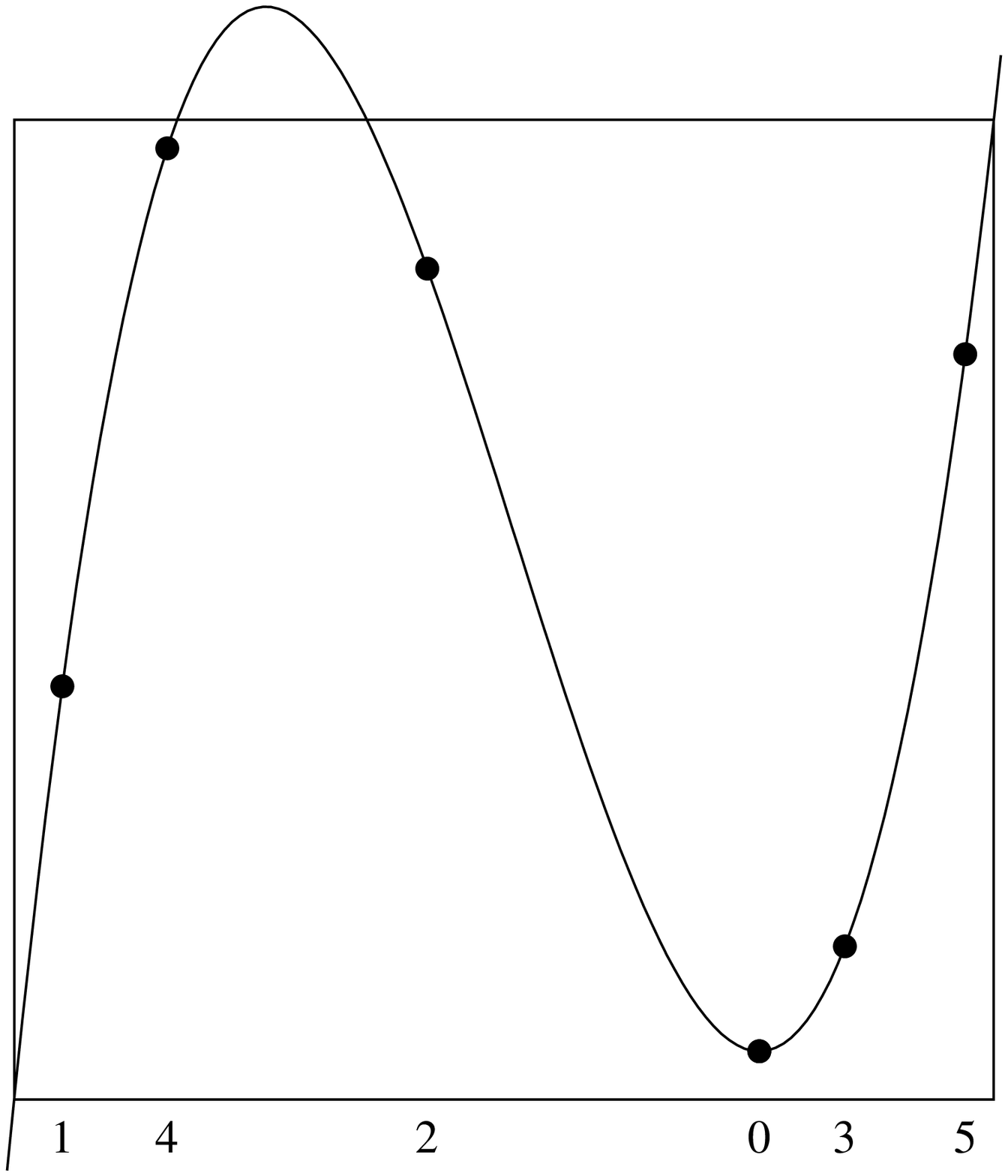,height=3.83in}}}
\caption{\textsl{A real example in $\cS_6$ with kneading sequence
  $\overline{100100}$ of period $3$ and marked grid of period $6$. The
  Julia set and puzzle structure
are shown above, and the graph of $F:\R\to\R$ below.
$($Similarly, Figure $5$ could be replaced by a pure imaginary
example, with $~a,\,v\in i\,\R\,.)$}
\label{f-100}}
\end{figure}

Again we are
 principally interested in the critical orbit $a_0\mapsto a_1\mapsto\cdots$.
 The knowledge of the critical kneading sequence is completely
 equivalent to the 
knowledge of the level one row of the critical marked grid; in fact
 $$\sigma(a_j)~ =~ 1-M(1,j)\,.$$ The notation
$$\vc{\sigma}~=~(~\sigma(a_1),\,\sigma(a_2),\,\sigma(a_3),\,\ldots~) $$
will be used for the kneading sequence of the critical orbit
(omitting the initial $\sigma(a_0)$ 
since it is always zero by definition). In the period $p$ case, we will
also write this as
$$ \vc\sigma~=~\overline{\sigma_1\sigma_2\cdots\sigma_{p-1}0}\,,$$
where the overline refers to infinite  repetition;
or informally just refer to the periodic sequence
 $\sigma_1\sigma_2\cdots\sigma_{p-1}0$.
 With this notation, note that the final bit must always be zero.
\end{definition}\smallskip

\begin{ex}\label{ex-1000}
For periods
 $~p~\le~3$,~ the kneading invariant together with the associated
quadratic map suffices
 to characterize the escape region, up to canonical
involution. However, for escape regions in $~\cS_4~$ with kneading sequence
$~\overline{1000}$,~ a ``secondary kneading invariant'' is needed to
specify whether the second and third forward images of $~a~$ are
``separate'' or ``together.'' (See Figure~\ref{f1000}.) For higher 
periods, many more such distinctions are necessary. (Compare 
Example~\ref{ex-1..0}.)
\end{ex}

\begin{figure}
\centerline{\psfig{figure=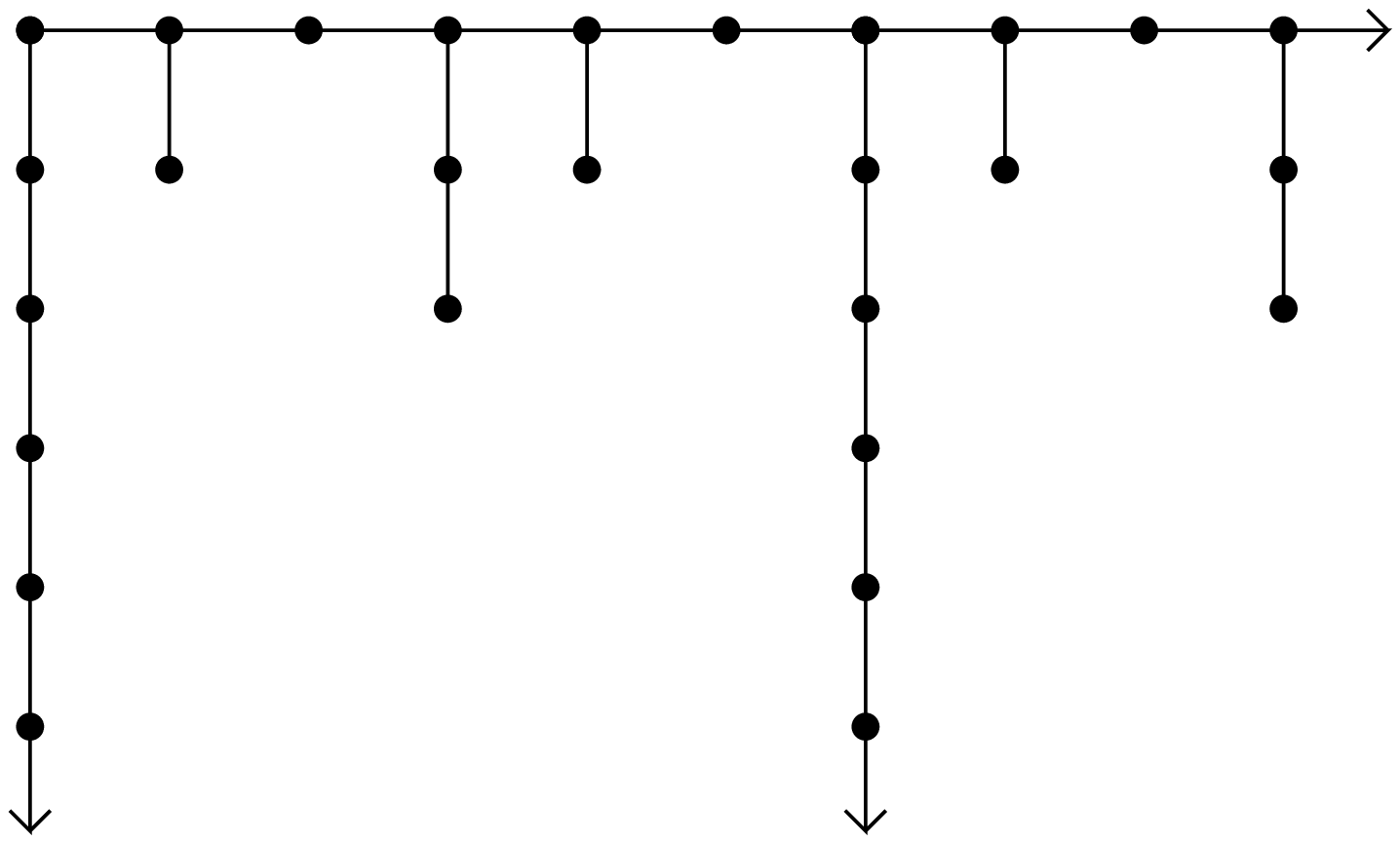,height=1.in}\qquad\qquad
\psfig{figure=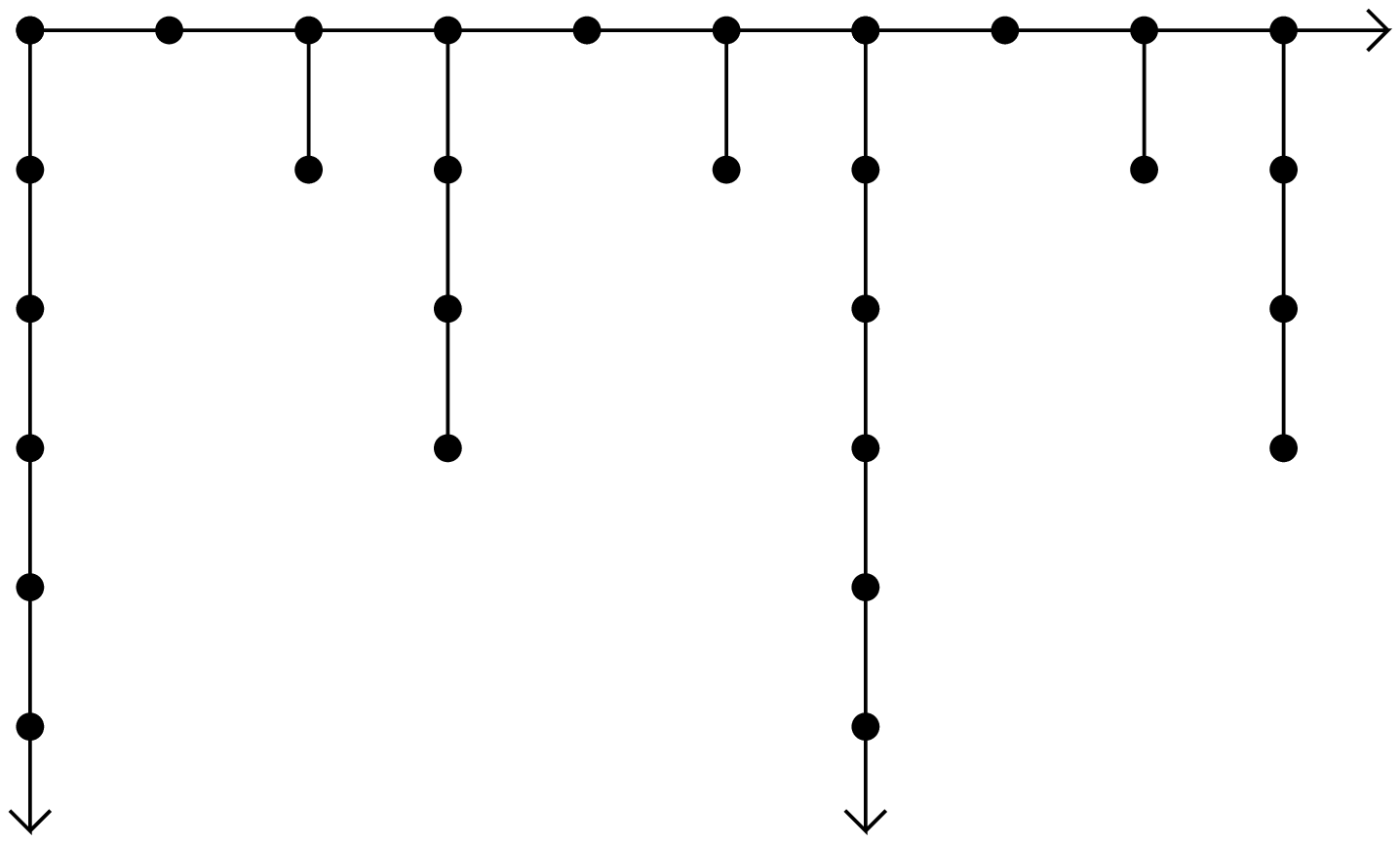,height=1.in}}
\caption{\textsl{ Marked grids corresponding to Figures~\ref{f-010}
and~\ref{f-100}.}\label{f-6grid}}
\end{figure}

\begin{ex}\label{ex-p6}
 Figures \ref{f-010} and \ref{f-100} show examples of maps
in $\cS_6$ which have kneading sequence of period 3 but marked grid
of period 6.
(Parameter values for
Figure 5: $~a\,=\,-0.8004\, +\,
0.2110\,i\,,\quad v\,=\,-0.77457275\,+\, 1.24396437\,i$,\break 
\noindent Figure 6: $~a\,=\,1.028778\,,\quad v\,=\,-1.877412$.)~~
 The corresponding marked grids are shown in Figure \ref{f-6grid}.
\end{ex}

\setcounter{table}{0}
\section{Puiseux Series.}\label{s-pui} 

(Compare \cite{K1}.)
Recall from Lemma~\ref{l-multe} that for each escape region 
$\cE_h\subset\cS_p$ the projection map $(a,v)\mapsto a$ from $\cE_h$ to
the complex numbers has a pole of order $\mu\ge 1$ at the
 ideal point $\infty_h$. It will be more convenient to work with the variable
$$ \xi~=~1/(3a) $$
which is a bounded holomorphic function throughout a neighborhood of
 $\infty_h$ in $\cS_p$. (The
 factor 3 has been inserted here in order to simplify
later formulas.) Since $\xi$ has a zero of order $\mu$, we can
choose some $\mu$-th root $~\xi^{1/\mu}~$ as a local uniformizing
parameter near the ideal point.

\begin{rem}\label{r-anch}
 This choice of local parameter $~\xi^{1/\mu}~$ is completely equivalent to the
choice of anchoring $~{\mathfrak B}(F)^{1/\mu}$~ 
 in Definition \ref{d-anch}.
In fact, the function $~F_{a,v}\mapsto {\mathfrak B}(F_{a,v})\,=\,
{\mathfrak B}_{F_{a,v}}(2\,a)~$ is asymptotic to
$~2\,a=2/(3\,\xi)$~  as $~|a|\to\infty$,~ so that the product
$~\xi\,\,{\mathfrak  B}(F)$~ converges to $~2/3\,$.
 Hence we can always choose the
$~\mu$-th roots so that their product converges to the positive root
$~(2/3)^{1/\mu}>0$.
\end{rem}
\smallskip

Let
 $$ a_0\mapsto a_1\mapsto a_2\mapsto\cdots\mapsto a_p=a_0 $$
be the periodic critical orbit, with $a=a_0$ and $F(a)=v=a_1$. Then each $a_j$
can be expressed as a meromorphic function of $~\xi^{1/\mu}~$, with a pole
at the ideal point. More precisely, according to
\cite[Theorem 5.16]{M4}, each $a_j$ can be expressed as a meromorphic
function of the form
\begin{equation*}\label{e-ajs}  
 a_j = \begin{cases}
\quad  a + O(1)  \qquad  {\rm if}\quad \sigma_j=0~,\cr
 -2a + O(1) \quad\;\;  {\rm if} \quad \sigma_j=1
\end{cases} 
\end{equation*}
where each $O(1)$ term represents a holomorphic function of
$~\xi^{1/\mu}~$ which is bounded for small $|\xi|$. (Compare 
Lemma~\ref{l-u-asym} below.) In order to replace the $a_j$ by holomorphic 
functions, we introduce the new variables
\begin{equation}\label{e-uj}
u_j~=~\frac{a-a_j}{3\,a}\,.
\end{equation}
Evidently each $u_j$ is a globally defined meromorphic function on $\cS_p$.
Within a neighborhood of the ideal point $\infty_h$,
each $u_j$ is a bounded holomorphic function of the
local uniformizing parameter $~\xi^{1/\mu}~$. In fact, this particular 
expression~(\ref{e-uj}) has been chosen so that $u_j$ takes the convenient value
$$u_j(0)~=~\sigma_j~\in~\{0,\,1\}$$
at the ideal point $\xi=0$. More precisely, each $u_j$ has a 
power series of the form
\begin{equation}\label{e-ujpui}
u_j~=~\sigma_j~+~a_\mu\xi~+~a_{\mu+1}\xi^{1+1/\mu}~+~a_{\mu+2}\xi^{1+2/\mu}~+~
\cdots
\end{equation}
which converges for small $~|\xi|$,~ 
with $\sigma_j$ as in Definition~\ref{d-kn}.
We will refer to this as the \textbf{\textit{~Puiseux expansion~}} of $u_j$,
since it is a power series in fractional powers of $\xi$.

Note that $u_0=u_p=0$ by definition. Recall the equation
$$ a_{j+1}~=~F(a_j)~=~a_j^3-3a^2a_j+2a^3+v~=~ (a_j-a)^2(a_j+2a)+a_1\,.$$
Substituting $a_j=a(1-3u_j)=a-u_j/\xi$,
 this reduces easily to the following equations, which will play a
 fundamental role.
\begin{equation}\label{e-u-condition}
\xi^2(u_{j+1}-u_1)~=~u_j^{\,2}(u_j-1)\,, 
\qquad{\rm or}\qquad u_{j+1}~=~u_1+u_j^{\,2}(u_j-1)/\xi^2\,.
\end{equation}
Thus, if we are given the Puiseux series for $u_1=(a-v)/3a$, then the
series for $u_2,\,u_3,\,\ldots,\,u_{p}$ can easily be computed
 inductively.\smallskip

\begin{ex}[{\bf Escape Regions in $\cS_2$}]\label{ex-p=2} 
{\rm Consider the case $p=2$. Since $u_2=0$, there is just one unknown
function $u_1$, which must satisfy Equation~(\ref{e-u-condition}), that is
 $~0=u_1+u_1^{\,2}(u_1-1)/\xi^2$,~
or in other words
$$u_1^{\,3}~-~u_1^{\,2}~+~\xi^2u_1~=~0\,.$$
This cubic equation in $u_1$ has three solutions, namely $u_1=0$
 corresponding to the unique escape region in
$\cS_1$, and the two solutions
$$u_1~=~\textstyle{\frac{1}{2}}\Big(1\pm\sqrt{1-4\xi^2}\Big)
~=~\textstyle{\frac{1}{2}}\Big(1\pm(1-2\xi^2-2\xi^4-4\xi^6-\cdots)
\Big) $$
or in other words
\begin{equation}\label{e-20}
u_1~=~\begin{cases}
1-\xi^2-\xi^4-2\xi^6-\cdots\, \quad & {\rm if\; the\; kneading\; sequence\;
 is~\;} \overline{10}\,, \cr
0 +\xi^2+\xi^4+2\xi^6+\cdots\,\quad   & {\rm if\; it\; is~\;} \overline{00}\,,
\end{cases}
\end{equation}
corresponding to the two escape regions in $\cS_2$.}
\end{ex}

Using the Equations (\ref{e-u-condition}), it is not difficult
to obtain asymptotic estimates.

\begin{lem}\label{l-u-asym}
For $~F\in\cE_h\subset\cS_p~$ and for $~0<j<p$,~
we have the following asymptotic estimates as $~\xi\to 0~$
or as $~|a|\to\infty$. 
\begin{eqnarray} 
u_j~&\sim&~\pm\xi\sqrt{u_1-u_{j+1}}\,\qquad{\rm if}\qquad \sigma_j=0\,,\label{e-u_0}\\ 
u_j-1~&\sim&~~\xi^2(u_{j+1}-u_1)\quad\qquad{\rm if}\qquad \sigma_j=1\,. \label{e-u_1}
\end{eqnarray}
Given only the kneading sequence $~\vc{\sigma}$,
we can state these estimates as follows.
\begin{eqnarray}
u_j~&=&~ \pm\xi\sqrt{\sigma_1-\sigma_{j+1}}~+~O(\xi^{3/2})\quad\qquad{\rm if}
\quad \sigma_j=0\,, \label{e-sig-0}\\
u_j~&=&~1~+~\xi^2(\sigma_{j+1}-\sigma_1)~+~O(\xi^3)\qquad\;{\rm if}
\quad \sigma_j=1\,.\label{e-sig-1}
\end{eqnarray}
\end{lem}

(In special cases, it is often possible to improve these
error estimates.) The proofs are easily supplied. 
\qed

\begin{rem}\label{r-ord-one}
As an example, still assuming that $~0< j< p$,~
 it follows that $~u_j~$ has the form
\begin{equation}\label{e-ord-one} 
u_j~=~\beta\,\xi~+~{\rm (higher~order~terms)}\qquad{\rm
 with}\qquad\beta\ne 0
\end{equation}
if and only if $~\sigma_j=0$~
and $~\sigma_{j+1}\ne \sigma_1$. Here $$\beta~=~\pm 1\qquad{\rm when}
\qquad\sigma_1-\sigma_{j+1}~=~+1\,,$$
 but $$\beta~=~\pm i\qquad{\rm when}
\qquad\sigma_1-\sigma_{j+1}~=~-1\,.$$
 In both cases, the equation
$~a_0-a_j=u_j/\xi$,~ implies that we have the following limiting formula
 as $~\xi\to 0~$ (or as $~|a|\to\infty\,$):
$$ \lim_{\xi\to 0}~ (a_j~-~a_0)~=~\beta~\ne~0\,.$$
It follows that $\cE_h\ne\I(\cE_h)$, since the difference $a_j-a_0$
changes sign when we replace the region $\cE_h$ by its dual $\I(\cE_h)$.
Still assuming (\ref{e-ord-one}), the approximate equality 
$$ a_0~-a_j~\approx~\beta\,, $$
 holds for all maps $~F\in\cE_h~$ with 
 $~|a_0|~$ large. As an example, Figure~\ref{f-010} illustrates a
 case with $a_1\approx a_0+i$.
In fact, assuming (\ref{e-ord-one}),
we can usually distinguish between maps in
the regions $\cE_h$ and $\I(\cE_h)$ simply by checking the sign of
the real part  $~\Re\big((a_0-a_j)/\beta\big)$.
\end{rem}\smallskip

\begin{rem}\label{r-u-cond} Recall that $u_p=0$  by definition.
The requirement that $u_p=0$, together with the set of 
Equations~(\ref{e-u-condition}),
imposes a very strong restriction as to which series
$u_1$ can occur. In fact this is the {\it only\/} restriction. A priori,
there could be formal power series solutions to the 
Equations~(\ref{e-u-condition}) with $u_p=0$ which
have zero radius of convergence, and hence do not correspond to any actual
escape region. However, every such solution satisfies a polynomial
 equation of degree
$3^{p-1}$ in the field of formal Puiseux series, and a counting argument
shows that every solution corresponds to an actual escape region
in some $\cS_n$, where $n$ must divide $p$. In particular, every solution has
 a positive radius of convergence, and the corresponding functions
$$ a~=~1/\,3\,\xi\,,\qquad v~=~\big(1-3u_1(\xi^{1/\mu})\big)/3\,\xi $$
parameterize a neighborhood of the ideal point in a corresponding escape region.
\end{rem}

\begin{rem}
In the case of an escape region of multiplicity $\mu>1$, some comment is needed.
There are $\mu$ possible choices for a preferred $\mu$-th root of the
holomorphic function $\xi$, and these will give rise to $\mu$ different
 power series for the function $u_1$. We can think of these various solutions
as elements of the field $\C((\idt^{1/\mu}))$ of formal power series.
 (The bold $\idt$ symbol
is used to emphasize that we are now thinking
 of $\idt$ as a formal indeterminate rather
than a complex variable.) However these different series for
the same escape region are all conjugate to each other under a Galois
automorphism $\idt^{1/\mu}\mapsto\alpha\idt^{1/\mu}$ from
$\C((\idt^{1/\mu}))$ to itself which fixes every point of the sub-field
 $\C((\idt))$. Here $\alpha$ is to be an
 arbitrary $\mu$-th root of unity.\footnote{This algebraic 
   construction has a geometric analogue.
The group of $\mu$-th roots of unity acts holomorphically on a
neighborhood of the ideal point. Let $\eta=\xi^{1/\mu}$ be a choice of
local parameter, and let $~\eta\mapsto\phi(\eta)\in\cE_h\subset\cS_p~$
be the holomorphic parametrization. Then the action is described by
$~\phi(\eta)\mapsto \phi(\alpha\,\eta)$,~ mapping a point $(a,\,v)$ to a
point of the form $~(a,\,v')$.}
 Thus we have outlined a proof of the following.

\begin{theo}
There is a one-to-one correspondence between escape regions in
$\cS_p$ and Galois conjugacy classes of solutions to the set of 
Equations~$(\ref{e-u-condition})$ with $u_p=0$, where $p$ is minimal.
\end{theo}
\end{rem}




\setcounter{table}{0}
\section{The Leading Monomial.}\label{s-mon}

Clearly it is enough to specify finitely many
terms of the Puiseux series $u_1,\,\ldots,$
$u_{p-1}$ in order to
 determine all of the $u_j$ uniquely, and hence determine the
 corresponding escape region uniquely. In fact, it turns out 
 that it is always enough to know
just the leading term of each $u_j$. (Equivalently, it is enough to
known the asymptotic behavior of each difference
 $a_j-a_0$ as $|a|=|a_0|$ tends to infinity. Compare
 Remark~\ref{r-ord-one}.)  
\smallskip

 It will be convenient to work with
formal power series. Let
$$\ub_j~=~\ub_j(\idt^{1/\mu})~\in~\C[[\idt^{1/\mu}]]$$
 be the formal power series
expansion for the holomorphic function $u_j=u_j(\xi^{1/\mu})$.
If a power series $~\x\in\C[[\idt^{1/\mu}]]~$ has the form
$$\x~=~\beta\,\idt^q~+~({\rm higher~order~terms})\qquad{\rm with}
\qquad \beta\ne 0\,,$$
then the monomial $~\m=\beta\,\idt^q~$ is called the 
 \textbf{\textit{\,leading term~}} of $~\x$,~ and the rational
 number $~q~$ is called the  \textbf{\textit{order}},
$$\ord(\x)~=~\ord(\m)~=~q\,.$$
For the special case $~\x={\bf 0}\,$, the order
is defined to be $~\ord({\bf 0})=+\infty$.\smallskip

 Closely related is the \textbf{\textit{~formal
power series metric}}, associated with the norm
\begin{equation}\label{e-norm}
\|\x\|~=~e^{-\ord(\x)}\,.
\end{equation}
This satisfies the \textbf{\textit{~ultrametric inequality}}
\begin{equation}\label{e-umi}
\|\x\pm\y\|~\le~ {\rm max}(\|\x\|,\,\|\y\|)\,.
\end{equation}
Two formal power series are \textbf{\textit{~asymptotically equal\/}},~
$\x~\sim~\x'$, if they have the same leading term.
 We will also use the $O(\idt^q)$ and $o(\idt^q)$
notations, defined as follows in the formal power series context: 
$$\begin{matrix}
\x=O(\idt^q)&\Longleftrightarrow&
\x~\equiv 0~~~~({\rm mod}~\idt^q)&\Longleftrightarrow&
\ord(\x)\ge q\,,&{\rm and}\\
\x=o(\idt^q)&\Longleftrightarrow&
\x~\equiv 0~~({\rm mod}~\idt^{q+1/\mu})&\Longleftrightarrow&
\ord(\x)> q\,.&
\end{matrix}$$
\smallskip

 It will be convenient to use notations such as
$~\stackrel{\to}{\w}~$ for a $p$-tuple\break $~(\w_1,\ldots,\,\w_{p-1}\,,0)$~
with $\w_p=0$,~ and  
 to define the ``error'' $E_j=E_j(\stackrel{\to}{\w})$ by the formula
\begin{equation}\label{e-err}
 E_j(\stackrel{\to}{\w})~=~\idt^2(\w_{j+1}-\w_1)\;-\;\w_j^{\,2}(\w_j-1)\,
\end{equation}
for $~0<j<p$.
Thus the required Equation (\ref{e-u-condition}) can be written briefly as
$~~E_j(\stackrel{\to}{\ub})~=~0\,.$ 

\begin{theo}[{\bf Kiwi}]\label{t-lmon} 
\it For any escape region $~\cE_h\subset\cS_p$,~ the vector
$$~\vc\ub~=~(\ub_1,\,\ub_2,\,\ldots,\,\ub_{p-1},\,\bzero)~$$
of Puiseux series is uniquely determined by the associated vector
$$~\vc\m~=~(\m_1,\,\m_2,\,\ldots,\,\m_{p-1},\,\bzero)~$$
of leading monomials.
\end{theo}

The proof will appear in $\cite{K2}\,$.\qed
\bigskip

Evidently these leading monomials $~\m_j~$ determine
and are determined by the asymptotic behavior
 of the meromorphic function
$~a-a_j=3au_j~$ as $~|a|\to\infty$ within the given escape region.
As an example, it follows from Remark \ref{r-ord-one}
 that $~\m_j=\pm\idt~$ if and only if
$$\lim_{|a|\to\infty} (a-a_j)~=~\pm 1\,.  $$

This theorem has a number of interesting consequences:

\begin{coro}\label{c-kiwi}
\begin{enumerate}
\item[{\bf(a)}] The multiplicity $~\mu~$ of the escape region 
$~\cE_h~$ is 
 equal to the least common denominator of the exponents
$$~q_j~=~{\rm ord}(\m_j)~=~{\rm ord}(\ub_j)\,.$$

\item[{\bf(b)}] Let $~E\supset \Q~$ be the smallest subfield of the algebraic
closure $\overline \Q$ which contains all of the coefficients $\beta_j$ of
the monomials $~\m_j=\beta_j\bxi^{q_j}$. Then the $~\ub_j~$ all
 belong to the ring $~E[[\bxi^{1/\mu}]]~$ of formal power series.\smallskip

\item[{\bf(c)}]
 In the special case where each $q_j$ is an even integer, it follows
that $~\ub_j \in E[[\bxi^2]].$~ This happens if and only if the escape region
 is invariant under the involution $\I$.
\end{enumerate}
\end{coro}

{\bf Proof.}  First note
 that the power series $\ub_j$ always belong to the subring
$$ \overline\Q[[\idt^{1/\mu}]]~\subset~\C[[\idt^{1/\mu}]]$$
consisting of series whose coefficients are algebraic numbers.
This statement follows from Puiseux's Theorem, which
states that the union over $\mu>0$ of the quotient fields
$\overline\Q((\xi^{1/\mu}))$ is algebraically closed,
together with the fact that
$\ub_1$ satisfies a polynomial equation with coefficients in
$\Q[\idt^2]$. (Compare Equation~(\ref{e-u-condition}) together with
 Remark \ref{r-u-cond}.)

Next, using the hypothesis that  the $~\ub_j~$ are uniquely determined
by the leading
terms $~\m_1,\,\ldots,\,\m_{p-1}$, we will show
 that all of the $~\ub_j~$ belong to the subring
${\mathbb E}[[\idt^{1/\mu}]] \subset \overline\Q[[\bxi^{1/\mu}]]$, as described
above. 
 To prove this statement, let ${\mathbb E}'$ be the smallest Galois
extension of $\mathbb E$ which contains all of the coefficients of
terms in the $\ub_j$. If $~{\mathbb E}'\ne{\mathbb E}$,~ then any
Galois automorphism of $~{\mathbb E}'~$ over $~\mathbb E~$ could be
applied to all of the coefficients of the $~\ub_j$ yielding a distinct
solution to the required equations. A similar argument shows that the
multiplicity $\mu$ is equal to the smallest positive
 integer such that every $\m_j$ is contained in $\overline\Q[\idt^{1/\mu}]$.
This proves assertions $\bf(a)$ and $(\bf b)$, and the proof of $(\bf c)$
is completely analogous.\qed\medskip

Here we will prove only the following
 very special case of Theorem \ref{t-lmon}.
However, this will cover many interesting examples, and the argument can
easily be used to describe an iterative algorithm for constructing
any finite number of terms of the series $~\ub_j$.
\smallskip

\begin{lem}\label{l-lead-mon}
Consider a solution $~\vc{\ub}$~ to the equations $~E_j(\vc{\ub})=0$,~
and let $\m_j$ be the leading term of $\ub_j$.
If $$~\ord(\m_j)~=~\ord(\ub_j)~<~2~$$ for $1\le j<p\,,$~
 then the power series $\ub_j$
are uniquely determined by these leading monomials $\m_j$.
\end{lem}

For the proof,
we will need to compare two vectors $~\stackrel{\to}{\w}~$ and 
$~\vc{\w}+\vc{\delbold}~$ in $\C[[\bxi^{1/\mu}]]$, which
satisfy $\delbold_j=o(\w_j)$ for all $j$, so that $\w_j$ and $\w_j+\delbold_j$
have the same leading monomial $\m_j$. As usual, we set
$\w_p=\delbold_p=0$, and for $1\le j<p$ assume that either
$\m_j=\sigma_j=1$,~ or that $~\m_j=O(\idt)~$ with $\sigma_j=0$.
It will be convenient to introduce the abbreviation
\begin{equation}\label{e-m*}
\m_j^{\textstyle\star}~~=~~\begin{cases}~~~~ \m_j~=~1 &
 ~~{\rm if}~~\sigma_j=1\,,\\
-2\,\m_j~=~O(\bxi) & ~~{\rm if}~~\sigma_j=0\,,\end{cases}
\end{equation}
or briefly $~\m_j^{\textstyle\star}=(3\sigma_j-2)\m_j$.

\begin{lem}\label{l-err-diff}
If $~\vc{\w}$ and $~\vc{\w}+\vc{\delbold}~$ satisfy
$$\w_j~\sim~\w_j+\delbold_j~\sim~\m_j $$
for all $j$, then
\begin{equation}\label{eq-error}
 E_j(\vc{\w}+\vc{\delbold})~=~ E_j(\w)\,+\,\idt^2(\delbold_{j+1}-\delbold_1)
\,-\,\m_j^{\textstyle\star}\delbold_j\,+\,o(\m_j^{\textstyle\star}\delbold_j)\,.
\end{equation}
\end{lem}

{\bf Proof.} Straightforward computation shows that
$~ E_j(\vc{\w}+\vc{\delbold})~$ is equal to
$$ E_j(\w)\,+\,\idt^2(\delbold_{j+1}-\delbold_1)
\,-\,(3\w_j^2-2\w_j)\delbold_j\,-\,(3\w_j+1)\delbold_j^2\,-\,\delbold_j^3\,.$$
It is not hard to check that $~3\w_j^3-2\w_j^2\sim\m_j^{\textstyle\star}$.~ 
Since $~\delbold_j=o(\m_j^{\textstyle\star})$, the conclusion follows.\qed\medskip

{\bf Proof of Lemma \ref{l-lead-mon}.}
Suppose there were two distinct solutions with the same leading terms,
so that
 $$ E_j(\vc{\ub})\,=\,E_j(\vc{\ub}+\vc{\delbold})~=~0\qquad{\rm with}
\qquad \ub_j\sim\ub_j+\delbold_j\sim\m_j$$
 for all $j$. 
It follows immediately from Lemma \ref{l-err-diff} that
$$ \idt^2(\delbold_{j+1}-\delbold_1)~\sim~\m_j^{\textstyle\star}\delbold_j \,,$$
hence 
$$\ord(\idt^2) +\ord(\delbold_{j+1}-\delbold_1)
~=~\ord(\m_j^{\textstyle\star})+\ord(\delbold_j) \,,$$
But by hypothesis,
$$ \ord(\m_j^{\textstyle\star})~=~\ord(\m_j)~<~\ord(\idt^2)~=~2\,,$$
so it follows that
$$\ord(\delbold_{j+1}-\delbold_1) ~<~\ord(\delbold_j)\,.$$
Thus, if we assume inductively that $\ord(\delbold_j)\ge q$ for all $j$,
then it follows that $\ord(\delbold_j)\ge q+1/\mu$ (using the fact that
$\ord$ is always an integral multiple of $1/\mu$).
 Iterating this argument, it follows that $~\ord(\delbold_j)~$
is greater than any finite constant, hence 
$~\delbold_j=0$,~ as required.\qed\bigskip

Still assuming that $\ord(\m_j)<2$, the next theorem provides a necessary
and sufficient condition for the existence of a solution to the equations
$~E_j(\stackrel{\to}{\ub})~=~0~$ with $\ub_j\sim\m_j$.

\begin{theo}\label{t-nas}
If $~\ord(\m_j)<2~$ for $0<j<p$, and
if we can find series $~\w_j\sim\m_j~$ with
\begin{equation}\label{e-approx-sol}
E_j(\stackrel{\to}{\w})~\equiv~0~~({\rm mod}~\idt^q)
\quad{\rm for~some}\quad q> 2\,\max_j\big(\ord(\m_j)\big)\,, 
\end{equation}
then there are $($necessarily unique$)$ series
$~\ub_j\sim\m_j$~ which satisfy the required
equations $~~E_j(\stackrel{\to}{\ub})~=~0\,.$
\end{theo}

{\bf Proof.} 
For any $(p-1)$-tuple $~\vc{\delbold} =\big(\delbold_1,\ldots, \delbold_{p-1})$~
satisfying $~\delbold_j=o(\m_j)$,~  we have
\begin{equation*}\label{e-correction}
E_j\big(\vc{\w}+\vc{\delbold}\big)~~=~~E_j(\vc{\w})
\,+\,\bxi^2(\delbold_{j+1}-\delbold_1)
 \,-\,\m^{\textstyle\star}_j\delbold_j\, +\,o(\m^{\textstyle\star}_j\delbold_j)\,,
\end{equation*}
by Lemma \ref{l-err-diff},
 with $~\m_j^{\textstyle\star}~$ 
as in Equation~(\ref{e-m*}). In particular, if we set
 \begin{equation}\label{e-w'}
 \w'_j~=~\w_j+\delbold_j\qquad{\rm with}
\qquad \delbold_j~=~ E_j(\vc{\w})/\m_j^{\textstyle\star}\,,
\end{equation}
then the terms $~E_j(\vc{\w})~$ and $~-\,\m^{\textstyle\star}_j\delbold_j~$  will
cancel, so that we obtain
\begin{equation}\label{e-E'}
 E_j(\vc{\w}{}')~=~
\idt^2\Big(\frac{E_{j+1}}{\m^{\textstyle\star}_{j+1}}-\frac{E_1}{\m^{\textstyle\star}_1}\Big)\,+\,o(E_j)\,.
\end{equation}
 Using Equation~(\ref{e-approx-sol}) together with the condition that
$~\ord(\m^{\textstyle\star}_j)<2$,~
it follows that $$E_j(\vc{\w}{}')= o(\idt^q)\,, \qquad{\rm or~in~other~words}
\qquad E_j(\vc{\w}{}')~\equiv~0\quad({\rm mod}~ \idt^{q+1/\mu})\,.$$
Furthermore, it follows that $~\w'_j\sim\w_j$.~
Iterating this construction infinitely often and passing to the limit,
we obtain the required series $\ub_j$ satisfying\break $E_j(\vc{\ub})=0$.
Since $\m_j\sim\w_j\sim\w'_j\sim\cdots\sim\ub_j$, this
completes the proof.\qed
\bigskip
                                             
\begin{ex}[{\bf Kneading sequence $\bf{\overline{1\cdots 10\cdots 0}}$}]
\label{ex-1..0}
 To illustrate
Theorem~\ref{t-nas}, suppose that the period $p$ kneading sequence satisfies
$$ \sigma_j~=~\begin{cases} 1\qquad{\rm for}\qquad  1~\le ~j~< ~k\,,\\
0  \qquad {\rm for} \qquad k~\le ~j~\le  ~p\,,
\end{cases}$$
with $~k>1~$ so that $~\sigma_1=1$.~ Then using Remark \ref{r-ord-one}
together with Equation~(\ref{e-ujpui}),
we see that the leading monomials have the form
$$\m_j~=~\begin{cases}~~ 1\qquad{\rm for} \qquad  1~\le ~j~< ~k\,,\\
\pm\idt \qquad {\rm for} \qquad k~\le ~j~<  ~p\,.
\end{cases}$$
Thus there are $2^{p-k}$ independent choices of sign. We will show that
each of these $2^{p-k}$ choices corresponds to a uniquely defined escape
region. In fact, to apply Theorem \ref{t-nas}, we simply choose the
approximating polynomials $\w_j$ to be\footnote{More precisely, 
the series $\ub_j$ have the form
 $~1-\idt^{2(k-j)}+{\rm (higher~terms)}~$ for $j<k$.}
$$ \w_j~=~\m_j\qquad{\rm for}\qquad j~\ne~k-1\,,\qquad{\rm but}
\qquad \w_{k-1}~=~1-\idt^2\,.$$
The required congruences $~~E_j(\vc{\w})~\equiv~0~~({\rm mod}~\idt^3)$~~
are then easily verified, and the conclusion follows.

Thus we have $~2^{p-k}~$ distinct escape regions, all with the same
kneading sequence (and in fact all with the same marked grid). For the case
$\bf{\overline{1000}}$, see Figure \ref{f1000}. The canonical involution $\I$
reverses the sign of each $\pm\idt$. For the cases with $k<p$, it follows that
$\I(\cE_h)\ne \cE_h$, so that there are only $2^{p-k-1}$ distinct regions
up to $180^\circ$ rotation.

On the other hand, for cases of the form
 $~\bf{\overline{11\cdots 110}}~$ with $k=p$,
there is only one escape region, invariant under the rotation $\I$. Closely
related is the fact that the Puiseux series for the $\ub_j$ all
contain only even powers of $\idt$ when $k=p$. In fact, it is not
 hard to show that these series all have integer coefficients,
 and hence belong to the subring~ $\Z[[\idt^2]]$. In fact, we 
 have $\m_j^{\textstyle\star}=1$ in Equation~(\ref{e-m*}),   so that
there are no denominators in Equations~(\ref{e-w'}) and (\ref{e-E'}).
\end{ex}

For further examples to illustrate Theorem~\ref{t-nas}, see Table~\ref{t-2}.

\subsection*{Computation of $~\ord(\m_j)$}

The following subsection will describe an algorithm 
which  computes the exponents of the leading monomials in terms of
the marked grid. 
In order to explain it, 
we must first review part of the Branner-Hubbard
theory.\medskip

 Recall from \S\ref{s-esc} that the puzzle piece $\,P_\ell(z_0)$ of
 level $\,\ell\ge 0$ for a map $\,F\in\cE_h$ is
the connected component containing $\,z_0$ in the open set
$$ \{\;z\in\C~~~;~~~G_F(z)~<~G_F(-a)/3^{\ell-1}\;\}\,.$$\,
The \textbf{\textit{~associated annulus~}} $\,A_\ell(z_0)\subset P_\ell(z_0)~$
is the outer ring
$$ \{\;z\in P_\ell(z_0)~~~;~~~G_F(z)~>~G_F(-a)/3^{\ell}\;\}\,.$$
(For $\ell=0$, this annulus is independent of the choice of $z_0$,
 and will be denoted simply by $A_0$.)\smallskip

For $\ell>0$, Branner and Hubbard show that the ``critical'' annulus
$A_\ell(a_0)$ maps onto $A_{\ell-1}(a_1)$ by a 2-fold covering map. It follows
that the moduli of these two annuli satisfy
$$    {\rm mod}\big(A_{\ell-1}(a_1)\big)~=~2\;{\rm mod}\big(A_\ell(a_0)\big)\,.$$
On the other hand in the non-critical case, $a_0\not\in P_\ell(z)$, the annulus
$A_\ell(z)$ maps by a conformal isomorphism onto $A_{\ell-1}\big(F(z)\big)$,
so the moduli are equal. In this way, they prove the following statement.

\begin{lem}[{\bf Branner and Hubbard}]\label{l-BH}
 For an arbitrary element $z\in K_F$, we have 
$$ {\rm mod}\big(A_\ell(z)\big)~=~{\rm mod}(A_0)/2^k\,,$$
where $k$ is the number of indices $i<\ell$ for which the puzzle piece
$F^{\circ i}\big(P_\ell(z)\big)$ contains the critical point $a_0$.
\end{lem}

 Here $~{\rm mod}(A_0)~$ is a non-zero constant (equal to
$~G_F(-a)/\pi$).~ It will be convenient to use the notation
 $~{\rm MOD}_{\textstyle\ell}(F)~$
for the ratio 
\begin{equation}\label{e-MOD}
 {\rm MOD}_{\textstyle\ell}(F)~=~{\rm MOD}\big(A_\ell(a_0)\big)~= 
 ~\frac
{2\;{\rm mod}\big(A_\ell(a_0)\big)}{{\rm mod}(A_0)}\,.\end{equation}
It follows from Lemma \ref{l-BH} that this is always a rational number
of the form $~1/2^{k-1}>0\,$,~
which is uniquely determined
 and easily computed from the associated marked grid (Definition \ref{d-mg}).
\smallskip

As in Equation~(\ref{e-L0}), the notation $~L_0(z)\ge 0~$ will denote the
supremum of levels $~\ell~$ such that $~P_\ell(z)=P_\ell(a_0)$.\medskip

\begin{theo}\label{t-kiwi}
For $~F\in\cE_h\subset\cS_p$,~  and for $~0<j<p$,~ the rational number
$~ {\rm ord}(\m_j)~=~{\rm ord}(\ub_j)~\ge ~0 ~$
is given by the formula
\begin{equation}\label{e-kiwi} {\rm ord}(\ub_j)~=
~\sum_{\ell=1}^{L_0(a_j)}~{\rm MOD}_{\textstyle\ell}(F)\,.
\end{equation}
\end{theo}\medskip

Here are three interesting consequences of this statement,

\begin{coro}\label{c0-pow2}
The multiplicity $~\mu~$ is always a power of two. 
\end{coro} 

{\bf Proof.} It follows immediately from Theorem \ref{t-kiwi} that the
denominator of $~\ord(\m_j)$,~ expressed as a fraction in lowest terms,
 is a power of two. The conclusion then follows from Theorem \ref{t-lmon}.\qed

\begin{coro}\label{c1-kiwi}
If an escape region in $\cS_p$ has trivial
kneading sequence $~\overline{0\cdots0}$, then 
\begin{equation}\label{e-sig=0}
{\rm ord}(\ub_j)~=~1+1/2+1/4+\cdots~=~2\,,\qquad{\rm for}\quad 0<j<p\,.
\end{equation}
\end{coro}

See Example~\ref{ex-0..0} below for more  
 precise information about this case.\medskip

{\bf Proof of Corollary \ref{c1-kiwi}.}
 As a first step, we show that $~d(a_0,\,a_j)=0~$ for all $j$.
(This is equivalent to the statement that $L_0(a_j)=\infty$ for all $j$;
or that all
grid points are marked.) Otherwise there would exist some $j$ with the
 largest value, say $D>0$,
of $d(a_0,\,a_j)$. It would then follow from the ultrametric inequality that
$~d(a_j,\,a_k)\le D~$ for all $j$ and $k$. But this is impossible
since it follows from
Equation (\ref{e-da0}) that $d(a_1,\,a_{j+1})=2D$.

Since all grid points are marked, the equation $~{\rm MOD}_\ell(F)=1/2^{\ell-1}$~
follows from Lemma \ref{l-BH}.~ Summing over $\ell>0$, the required equality
(\ref{e-sig=0}) then follows from the Theorem.\qed\smallskip

\begin{coro}\label{c2-kiwi}
For an arbitrary escape region, we have\footnote{For an alternative 
 necessary and sufficient condition, compare  Remark \ref{r-ord-one}.}

\begin{equation*}\label{e-ord-1}
{\rm ord}(\ub_j)=1\qquad{\rm if~ and~ only~ if}\qquad L_0(a_j)=1\,.
\end{equation*}
 \end{coro}

{\bf Proof.} This follows easily, using the observation that
$~{\rm MOD}_{\textstyle 1}(F)=1$~ in all cases.\qed
\smallskip

\medskip

{\bf Proof of Theorem \ref{t-kiwi}.} Let  $~\Q((\btau))~$ be the field
of Laurent series in one formal variable, and let $~\mathbb S~$ be the metric
completion of the algebraic closure of $~\Q((\btau))$,~ 
using the metric associated with the norm
$$ \|\z\|~=~e^{-{\rm ord}(\z)}\,, $$
with $~{\rm ord}(\btau)=1$.~ (Note that the field
 $~\overline\Q$~ of algebraic numbers is a discrete subset of $~\mathbb S$,~
with $~\|q\|=1$~ for $~q\ne 0$~ in $~\overline\Q$.)
\smallskip

 The argument will be based on \cite{K1}, which
 studies the dynamics of polynomial maps from this field
 $~{\mathbb S}$~ to itself, using methods developed by Rivera-Letelier 
\cite{RL1} for the analogous $p$-adic case.
Identify the indeterminate $\btau$ with $1/\ab$, and consider the map
\begin{equation}\label{e-bF}
 \bF(\z)~=~ \bF_\vb(\z)~=~ (\z-\ab)^2(\z+2\ab)\;+\;\vb 
\end{equation}
from $~\mathbb S~$ to itself, where $~\vb~$ is some constant
in $~{\mathbb S}$.~
(Kiwi actually works with the conjugate map
 $~\bzeta\mapsto \bF(\ab\,\bzeta)/\ab$,~ with critical points $~\pm 1$,~
but this form (\ref{e-bF}) will be more convenient for our purposes.)
Note that
 $$  \log\|\ab\|~=~\log\|3\,\ab\|~=~+1\,,\qquad{\rm and~in~general}\quad
\log\|\z\|=-{\rm ord}(\z)\,.$$
Following \cite{RL1}, any subset of the form
$$ \bA~=~ \bA(\lambda,\,\lambda'\,,~\widehat\z)~=~\{\;\z\in{\mathbb S}~~;~~ \lambda<
\log\|\z-\widehat\z\|<\lambda'\;\} $$
where $~\lambda<\lambda'$,~ is
 called an \textbf{\textit{~annulus}}~ surrounding $~\widehat\z$,~ with
 \textbf{\textit{~ modulus~}}
$$\mathfrak{Mod}(\bA)~=~\lambda'-\lambda~>~0\,.$$
By definition,  the \textbf{\textit{~level zero annulus~}}
associated with any such map $~\bF$~ is the annulus
$~\bA_0\,=\,\bA(1,\,3,~\bzero)$,~ 
 with $~\mathfrak{Mod}(\bA_0)= 2$.  Kiwi shows that each iterated
pre-image $~\bF^{-\ell}(\bA_0)~$ can be uniquely expressed as the union
of finitely many disjoint annuli, which are called \textbf{\textit{~annuli
of level $\ell$}}. Given any point $~\widehat\z\in K_\bF~$ (that is any
point with bounded orbit), there is a uniquely defined
infinite sequence of nested annuli of the form
$$\bA_\ell(\widehat\z)~=
 ~\bA(\lambda_{\ell+1},\,\lambda_\ell\,;\;\widehat\z)\,,$$
each $~\bA_\ell(\widehat\z)~$
surrounding $~\bA_{\ell+1}(\widehat\z)$,~ and all surrounding
 $~\widehat\z$.~ Note that $~\lambda_0=3$ and $~\lambda_1=1$;~
but otherwise the numbers $~\lambda_0>\lambda_1>\lambda_2>\cdots~$
depend on the choice of $~\widehat\z$.~  An
 annulus of level $~\ell~$ is called \textbf
{\textit{~critical~}} if it surrounds the critical point $~\ab~$(or in other words is equal to $~\bA_\ell(\ab)$),~ and is
\textbf{\textit{~non-critical~}} otherwise.
\smallskip

We can also describe this construction in terms of the Green's function
$$ G(\z)~=~\lim_{n\to\infty}\; \frac{1}{3^n}\;\log^+\|F^{\circ n}(\z)\|\,.$$
This satisfies $~G\big(\bF(\z)\big)\,=\,3\,G(\z)$, with $~G(\z)\le 1~$
whenever $~\log\|\z\|\le 1$,~ but with
$$ G(\z)~=~\log\|\z\|~>1\qquad{\rm whenever}\quad \log\|\z\|>1\,.$$
It follows that the set $~\bF^{-\ell}(\bA_0)~$ can be identified with
$$ \left\{\;\z\in{\mathbb S}~~;
~~\frac{1}{3^\ell}~<~G(\z)~<~\frac{3}{3^\ell}\;\right\}\,.$$

Just as in the Branner-Hubbard theory,
 there are only two annuli of level one, namely 
$$\bA_1(\ab) ~=~\bA(0,\,1,\,\ab)\quad{\rm of~modulus~one}, $$
and
$$ \bA_1(-2\,\ab)~=~\bA(-1,\,1,\,-2\,\ab)\quad{\rm of~modulus~two}\,.$$
Furthermore, each annulus of level $~\ell>0~$
maps onto an annulus of level $~\ell-1$,
$$ \bF\big(\bA_\ell(\widehat\z)\big)~=~\bA_{\ell-1}\big(\bF(\widehat\z)\big)\,,$$
with
$$\mathfrak{Mod}\Big(\bF\big(\bA_\ell(\widehat\z)\big)\Big)~=~\begin{cases}
2\,\mathfrak{Mod}\big(\bA_\ell(\widehat\z)\big)\qquad{\rm if}
\quad \bA_\ell(\widehat\z) \quad {\rm is~critical}\,,\\
 ~~\mathfrak{Mod}\big(\bA_\ell(\widehat\z)\big)
 \qquad{\rm otherwise}\,.\end{cases}$$
Furthermore, just as in the Branner-Hubbard theory, if the critical orbit
$$\ab=\ab_0~\mapsto~\ab_1~\mapsto~\ab_2~\mapsto~\cdots $$
is contained in $~K_\bF$,~ 
 then it gives rise to a marked grid which describes which
annuli $~\bA_\ell(\ab_j)~$ are critical. This grid can then be used to compute
the moduli of all of the annuli $~\bA_\ell(\ab_j)$.

\begin{figure}
\centerline{\psfig{figure=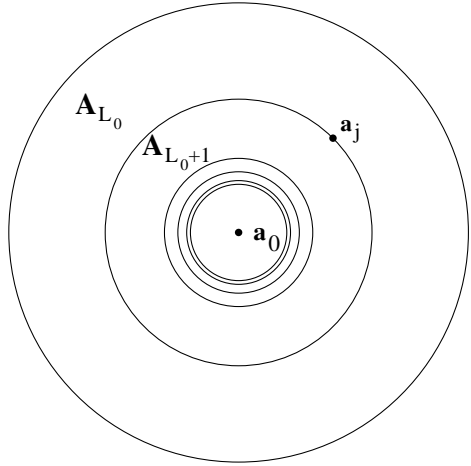,height=2.4in}}
\caption{\textsl{Schematic picture of concentric annuli 
 in $~\mathbb S\,$.}\label{f-annuli}}
\end{figure}

Now consider the sequence of critical annuli. Let
$$ \bA_\ell(\ab_0)~=~\bA(\lambda_{\ell+1},\,\lambda_\ell\,,\;\ab_0)\,,
\qquad{\rm so~that} 
\qquad \mathfrak{Mod}\big(\bA_\ell(\ab_0)\big)~=\lambda_\ell-\lambda_{\ell+1}\,.$$
Since $\lambda_1=1$, we can write
$$\mathfrak{Mod}\big(\bA_1(\ab)\big)\,+\,
\cdots\,+\,\mathfrak{Mod}\big(\bA_\ell(\ab)\big)~=~(\lambda_1-\lambda_2)
\,+\,\cdots\,+\,(\lambda_\ell-\lambda_{\ell+1})~=~1-\lambda_{\ell+1}\,.$$
As in Equation (\ref{e-L0}), 
set $~L_0=L_0(\ab_j)~$ equal to the
 supremum of values of $~\ell~$ such
 that the annulus $~\bA_\ell(\ab_j)~$ is critical. First suppose that
 $~L_0~$ is finite. This means that the annulus 
$~\bA_{L_0}(\ab_0)~$ surrounds $~\ab_j$,~ but that $~\bA_{L_0+1}(\ab_0)~$
does not surround $~\ab_j$.~
 It then follows easily that $~\log\|\ab_0-\ab_j\|~$ must be precisely
equal to $~\lambda_{L_0+1}$.~ (Compare Figure \ref{f-annuli}.)
 Since $~\ub_j=(\ab_0-\ab_j)/3\ab_0$,~ this means that
$$ \log\|\ub_j\|~=~\lambda_{L_0+1}-1\,,$$
hence
\begin{equation}\label{e-kiwi2}
{\rm ord}(\ub_j)~=~-\log\|\ub_j\|~ =~ 1-\lambda_{L_0+1}~=~
\sum_{\ell=1}^{L_0}\;\mathfrak{Mod}\big(\bA_\ell(\ab_0)\big)\,.
\end{equation}

Next we must prove this same formula (\ref{e-kiwi2}) in the case 
$L_0~=~L_0(\ab_j)~=~\infty\,$.
Note that $L_0(\ab_j)$ can be infinite only in the renormalizable case,
with $j$ a multiple of the grid period $n$, which
 is strictly less than the critical period $p$.\smallskip
 
It is not hard to check that the limit
 $~\lambda_\infty=\lim_{\,\ell\to\infty} \lambda_\ell$~ exists, and is a finite
rational number. {\it We must prove that
\begin{equation}\label{e-lam-inf}
\log\|\ab_0-\ab_j\|~=~\lambda_\infty
\end{equation}
whenever $j$ is a multiple of $n$,~ with $~0<n<p$.~} Given this equality,
 the proof  of Equation (\ref{e-kiwi2}) will go through just as
 before, so that
$${\rm ord}(\ub_j)~ =~ 1-\lambda_{\infty}~=~
\sum_{\ell=1}^{\infty}\;\mathfrak{Mod}\big(\bA_\ell(\ab_0)\big)\,.$$

For $~j\equiv 0~({\rm mod}~ n)$~ with $~0<j<p$,~ it
 is not hard to check that $~\bF^{\circ j}~$ maps the disk
$~\{\z\in{\mathbb S}~;~\log\|\z-\ab_0\|\le\lambda_\infty\}~$ onto itself,
 with $~\ab_0~$ as a critical point. It will be convenient to make a change
of variable, setting
 $~\bzeta=(\z-\ab_0)\,\btau^{\lambda_\infty}$,~ where $\btau=1/\ab$~
 so that
$$ \log\|\bzeta\|~=~\log\|\z-\ab_0\|-\lambda_\infty\log\|\ab\|
~=~\log\|\z-\ab_0\|-\lambda_\infty\,, $$
and hence
$$ \|\bzeta\|~\le~1\qquad\Longleftrightarrow 
 \qquad \log\|\z-\ab_0\|~\le~\lambda_\infty\,. $$
In this way, we see that  $~\bF^j~$
 is conjugate to a map $~\bzeta\mapsto\bphi(\bzeta)~$  which sends the
``unit disk'' $~{\mathbb U}=\{\bzeta~;~\|\bzeta\|\le 1\}~$ onto itself.
 Here $~\bphi~$ is a polynomial map, say
\begin{equation}\label{e-pui-poly}
 \bphi(\bzeta)=\bc_0+\bc_1\bzeta+\cdots+\bc_N\bzeta^N\,.
\end{equation}
Since $\bphi({\mathbb U})\subset{\mathbb U}$,~
it is not hard to check that
all of the coefficients must satisfy $~\|\bc_k\|\le 1$.~
Using this change of variable, the marked critical point $~\ab_0~$
and its image $~\bF^{\circ j}(\ab_0)=\ab_j~$ correspond to the 
critical point $~\bzeta=\bzero~$ and its image
$~\bc_0=\bphi(\bzero)$.~
Thus the required equality~(\ref{e-lam-inf})
translates to the statement that $~\|\bc_0\|=1$.

Since the origin is a critical point of $\bphi$, it follows that $\bc_1=0$.
Suppose that $~0<\|\bc_0\|<1$.~ Then it follows inductively that 
$$ \bF^{\circ k}(\bc_0)~=~ \bc_0\;+\;O(\bc_0^{\,2}) $$
for all $k$.
Hence the origin cannot be periodic.~ This contradiction proves that
$\|\bc_0\|=1$,~ hence $\log\|\ab-\ab_n\|=\lambda_\infty$. The proof of equation
(\ref{e-kiwi2}) then goes through just as before.\medskip

To finish the argument, we need only quote \cite[Proposition 6.17]{K1},
which can be stated as follows in our terminology. Suppose that $~F~$
belongs to an escape region $\cE_h\subset \cS_p$, and suppose that
$~\vb\in{\mathbb S}~$ is the Puiseux series attached to the ideal point of
$~\cE_h~$ which expresses the parameter $~v~$ as a function of $~1/a$.~
Then the marked grid associated with the critical orbit of the complex
map $~F~$ is identical 
with the marked grid for the critical orbit of the map $~\bF=\bF_\vb$~
of formal power series. It follows that
the algebraic moduli $~\mathfrak{Mod}\big(\bA_\ell(\ab_0)\big)~$ are identical
 to the geometric moduli  $~{\rm MOD}\big(A_\ell(a_0)\big)~$
of Equation (\ref{e-MOD}). The conclusion of Theorem \ref{t-kiwi}
then follows immediately.\qed
\medskip

\subsection*{Examples}

\def\2dots{\cdot\,\cdot}

Table \ref{ta-prim}  lists the first 
{\it two\/} terms of the series $\ub_j$ for each primitive
escape region in $\cS_p$, with $2\le p\le 4$. (An
escape region in $\cS_p$ is called
 \textbf{\textit{~primitive~}} if its marked grid
 has period exactly $p$.)

\begin{table}[!ht]
\begin{center}
\begin{tabular}{lccccl} 
$\vc{\sigma}$ & $\ub_1$ & $\ub_2$ & $\ub_3$ &  $\#~$ & $\mu$\cr \cr
$\overline{10}$ & $1-\idt^2+\2dots$ &         &      & 1 & 1  \cr
$\overline{110}$ & $1- \idt^4+\2dots$ & $1- \idt^2+\cdots$ & & 1 & 1\cr
$\overline{100}$ & $1 - \idt^2+\2dots$ & $\pm  \idt+\idt^2/2+\cdots$ &  & 2 & 1\cr
$\overline{010}$ & $\pm  i\idt+\idt^4/2+\cdots$ &$ 1 \mp  i\idt^3+\cdots$ & & 2 & 1\cr
$\overline{1110}$ & $1-\idt^6+\2dots$ & $1-\idt^4+\2dots$ & $1-\idt^2+\2dots$ & 1 & 1\cr
$\overline{1100}$ & $1-\idt^4+\2dots$ & $1-\idt^2+\2dots$ & $\pm \idt+\idt^2/2 +\2dots$ & 2 & 1\cr
$\overline{0110}$ & $\pm i(\idt+\idt^3/2)+\2dots$ & $1+\idt^2+\2dots$ &
 $1\mp i\idt^3+\2dots$ & 2 & 1\cr
$\overline{1000}$s & $1-\idt^2+\2dots$ & $\pm \idt+\idt^2+\2dots$ & $\mp \idt+\idt^2/2+\2dots$ & 2 & 1\cr
$\overline{1000}$t & $1-\idt^2+\2dots$ & $\pm( \idt-3\idt^3/4)+\2dots$ & $\pm \idt+\idt^2/2+\2dots$ & 2 & 1\cr
$\overline{0100}$ & $\omega^2\idt+\idt^4/2+\2dots$ & $1-\omega^2\idt^3+\2dots$ &
 $-\omega\idt^{3/2}+\omega^2\idt^3/2+\2dots$ & 2 & 2\cr
$\overline{0010}$ & $\omega\idt^{3/2}+\idt^2/2+\2dots$ & $-\omega^2\idt-\idt^2/2+\2dots$ &
 $1-\omega\idt^{7/2}+\2dots$ & 2 & 2\cr
\end{tabular}
\end{center}
\vspace{-.3cm}
\caption{\textsl{Primitive orbits of periods 2, 3 and 4. Here \label{t-2}
$\vc{\sigma}$ is  the kneading sequence,~ 
$\#$ denotes the number of solutions up to conjugacy, $~\mu$ is the
multiplicity, and
$\omega$ denotes an arbitrary 4-th root of $-1$, or in other words
$~\omega=(\pm 1\pm i)/\sqrt 2$.\label{ta-prim}}}
\end{table}
\medskip

For all of the cases in this table, the hypothesis that $~\ord(\m_j)<2~$ is\break
satisfied (compare Theorem \ref{t-nas});  although this hypothesis
 is certainly not satisfied in general.
For some of these entries, it would suffice to take
$~\w_j~$ equal to the leading term $~\m_j\,$,~ in order to satisfy
 the requirement of Equation~(\ref{e-approx-sol}).
However, in all of these cases it would suffice to take the two 
initial terms, as listed in the table.
In many cases, it is quite easy to compute the second term of $\ub_j$,
provided that we know the initial terms for $\ub_1,\,\ub_j\,,$ 
and $\ub_{j+1}$.~ In particular, if  $\m_j=1$, and $\m_{j+1}\ne\m_1$,
then we can compute the second term of $\ub_j$ from the estimate
$$ \ub_j-1 ~\sim~\idt^2(\ub_{j+1}-\ub_1)~\sim
~\idt^2(\m_{j+1}-\m_1) $$
of Equation (\ref{e-u_1}).
\medskip

\begin{rem}
Note that the canonical involution $\I$ of
 Remark \ref{r-I} corresponds\footnote{ 
The case $\mu>1$ is more complicated.
The involution $~\I:\idt\mapsto -\idt~$ of $~\C[[\idt]]~$ lifts to an
isomorphism
$$\I_\mu:\idt^{1/\mu}~\mapsto~e^{\pi i/\mu}\,\idt^{1/\mu}\qquad
{\rm of}\qquad \C[[\idt^{1/\mu}]] 
$$
which is not an
 involution, although $~\I_\mu\circ\I_\mu~$ is Galois conjugate
 to the  identity.}
 to the 
involution $\idt\leftrightarrow -\idt$ of $\C[[\idt]]$.
In a number of these cases, there
 are two distinct solutions in $\cS_p$ which maps to each other under this
 involution. Geometrically, this
means that two different escape regions in $\cS_p$ fold together
into a single escape
region in the moduli space
$\cS_p/\I$. The case with kneading sequence $\overline{1000}$
is exceptional among examples with $~p\le 4$,~ 
since in this case there are four distinct escape regions
in $\cS_4$, corresponding to two essentially different escape regions
in $\cS_4/\I$. (Compare Figure~\ref{f1000} and Table~\ref{t-2}, as well
as Example \ref{ex-1..0}.) 
\smallskip

In each case, according to \cite[Lemma 5.17]{M4}, the total number of solutions
associated with a given kneading sequence
 $~\overline{\sigma_1\cdots\sigma_{p-1}0}$,~
 counted with multiplicity, is equal to\begin{equation}\label{e-count} 
 2^{(1-\sigma_1)+\cdots+(1-\sigma_{p-1})}\,.\end{equation}
\end{rem}
\smallskip

\begin{ex}[{\bf Kneading sequence $\bf\overline{000\cdots 0}\,$}]\label{ex-0..0}
In the case of a trivial kneading sequence 
 $~\sigma_1=\cdots=\sigma_{p}=0$,~
we know from Corollary \ref{c1-kiwi}
that $~{\rm ord}(\ub_j)=2~$ for $1<j<p$.
  Thus the hypothesis that $\ord(\m_j)<2$ is not satisfied, so the arguments
above do not work. However, we can still classify the solutions $\ub_j$
 by working just a little bit harder. 

Even without using Theorem \ref{t-lmon}, it follows easily from the
Equations~(\ref{e-u-condition}) that $\ub_j\equiv 0~~({\rm
mod}~\idt^2)$. If we set $\m_j=-\idt^2c_j$, then it follows from these
equations that
$$ c_{j+1}~=~c_j^{\;2}\,+\,c_1$$
with $c_p=0$. In other words, the complex numbers $c_j$ must form
a period $p$ orbit
$$ 0~\mapsto~c_1~\mapsto~c_2~\mapsto~c_3\mapsto~\cdots~\mapsto~c_{p-1}
~\mapsto~0 $$
under the quadratic map $~Q(z)=z^2+c_1$, so that $c_1$ is the center
of a period $p$ component in the Mandelbrot set. Let ${\mathbb E}=
\Q[c_1]$ be the field generated by $c_1$.
\medskip

\begin{table}[!ht]
\begin{center}
\begin{tabular}{cccllc}
$p$ & & $c_1$ & & nickname & $\psi_p=\lim\, a/t$\\ \\
1 & & 0 & & $z^2$ & 1\\
2 & & -1 & & basilica& -1\\
3 & & -.12256+.74486\,{ i} & & rabbit & -1.675-1.125\,i\\
3 & & -1.75488 & &airplane&-5.649\\
4 & &  -0.15652 + 1.03225\,{ i}  & & kokopelli&-9.827-1.392\,i\\
4& & 0.28227 + 0.53006\,{ i}  & & $(1/4)$-rabbit&-2.273-2.878\,i  \\
4 & &  -1.94080   & & worm &-25.534 \\
4 & &  -1.31070  & & double-basilica& 1.734
\end{tabular}
\end{center}
\vspace{-.3cm}
\caption{\textsl{ \label{t-2a} List of quadratic Julia sets with critical period
$p\le 4$ and with critical value in the upper half-plane. 
$($Compare Figure $\ref{f-quad-j}.)$  The last column
gives the limit of $~a/t$ as $a\to\infty$
 for the associated cubic escape region with trivial kneading sequence.
$($See Theorem $\ref{t-dt}$ and Equation $(\ref{e-t-sim2}).)$}}
\end{table}
\bigskip

\begin{figure}[!ht]
\centering
\subfigure[\textsl{Kokopelli}]{\psfig{figure=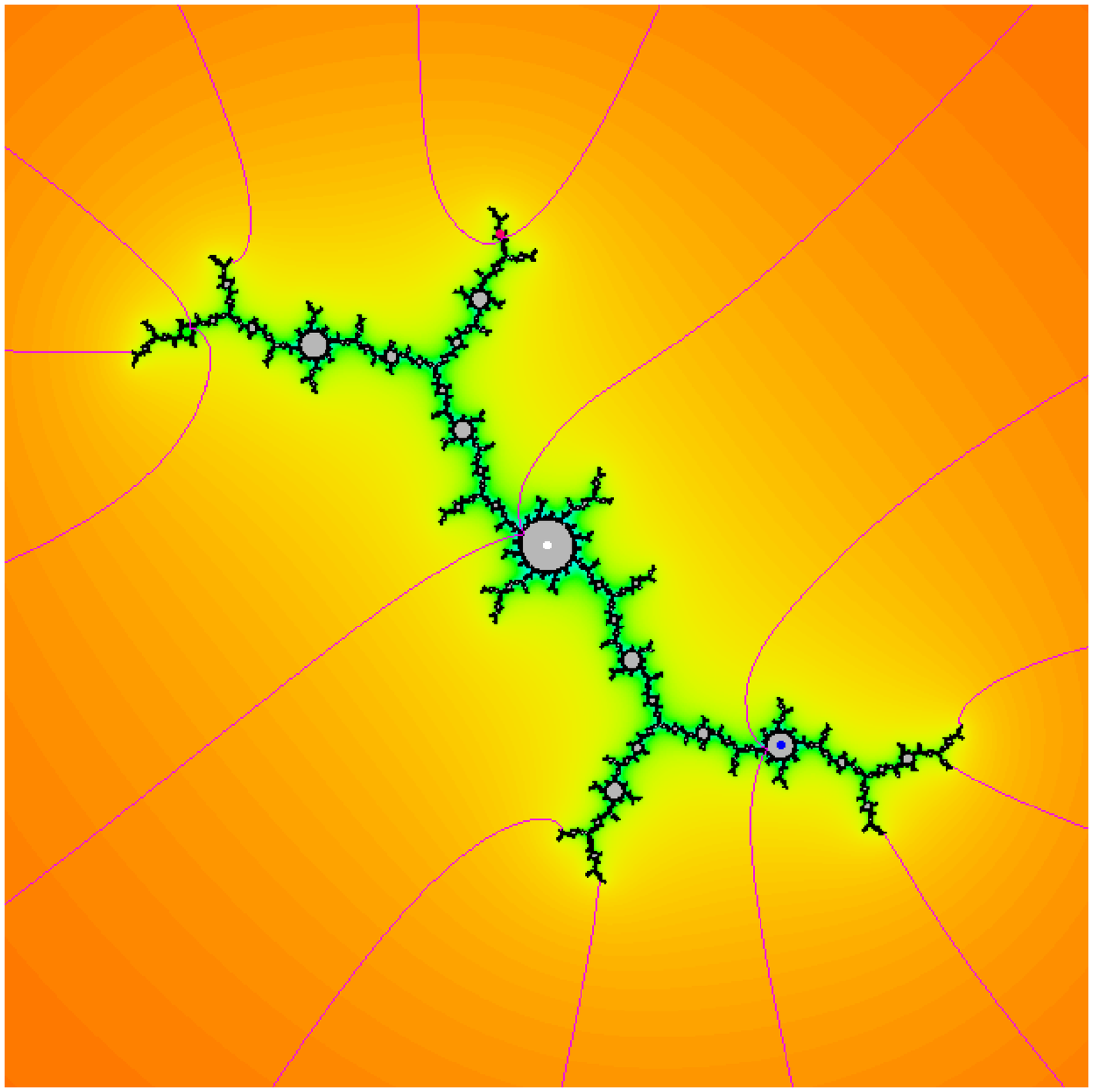,height=2.2in}}\qquad
\subfigure[\textsl{(1/4)-rabbit}]{\psfig{figure=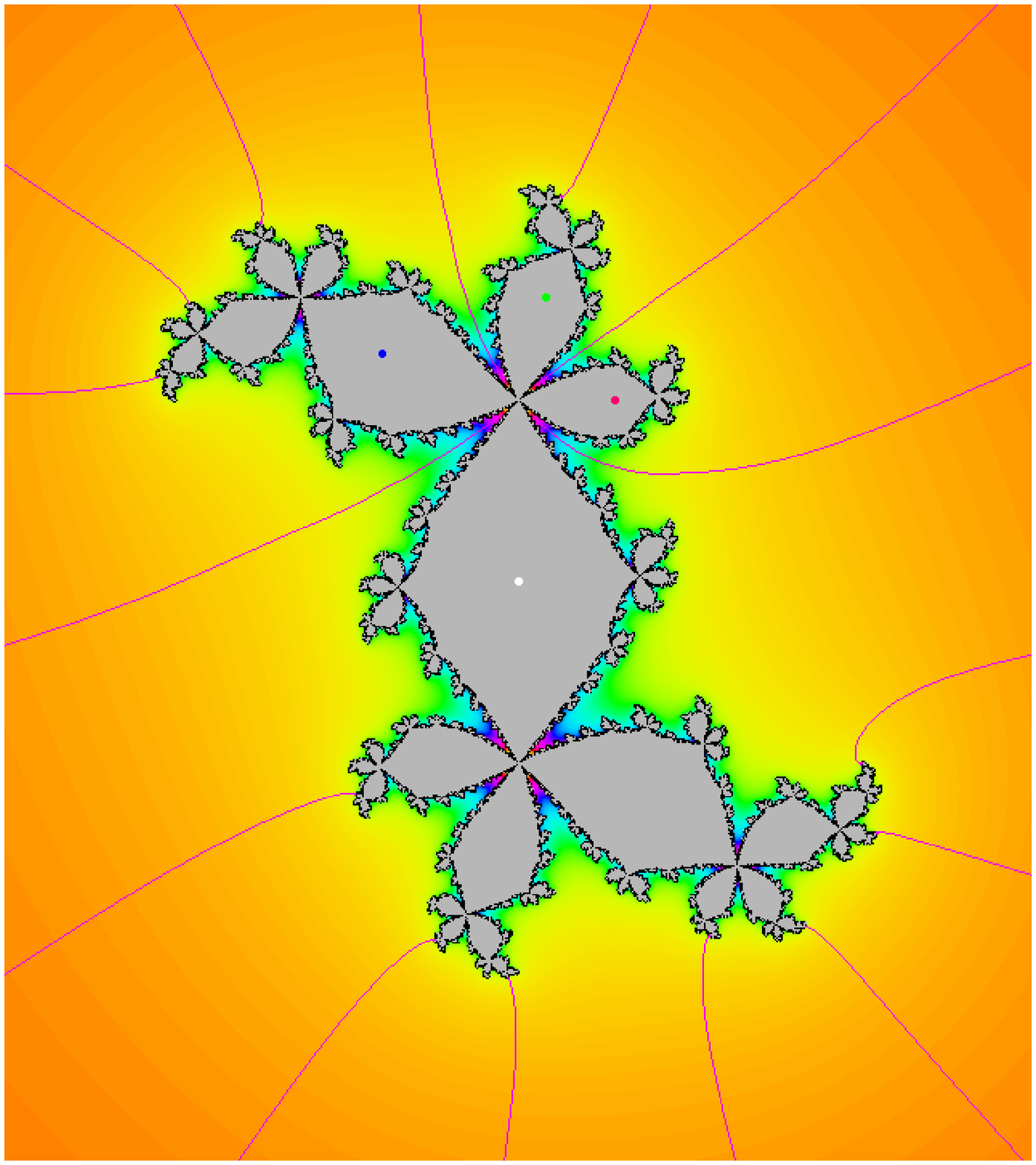,height=2.2in}}\\
\subfigure[\textsl{worm}]{\psfig{figure=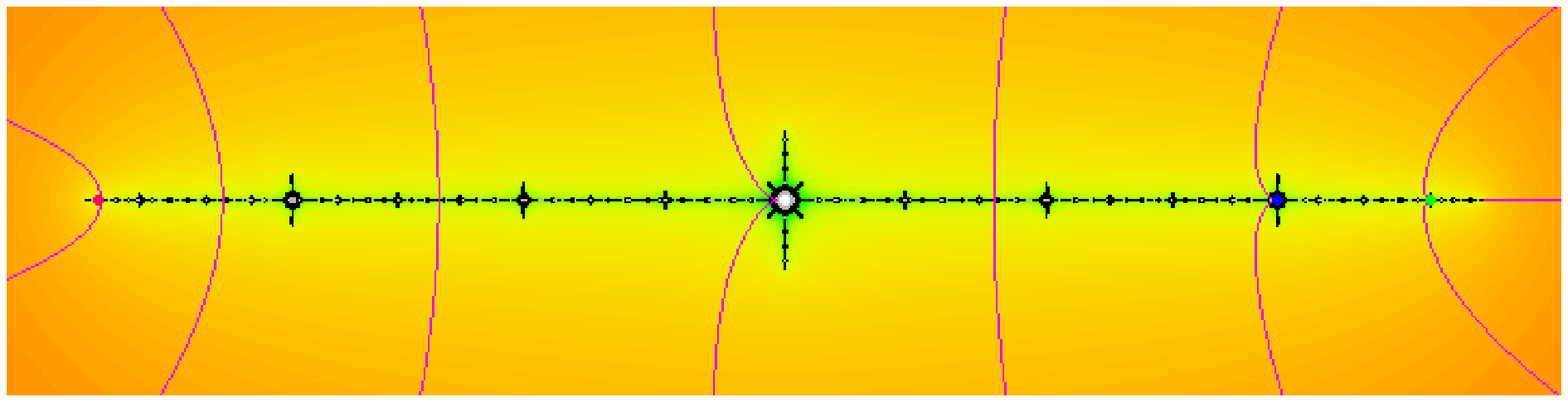,height=.7in}}\qquad
\subfigure[\textsl{double basilica}]{\psfig{figure=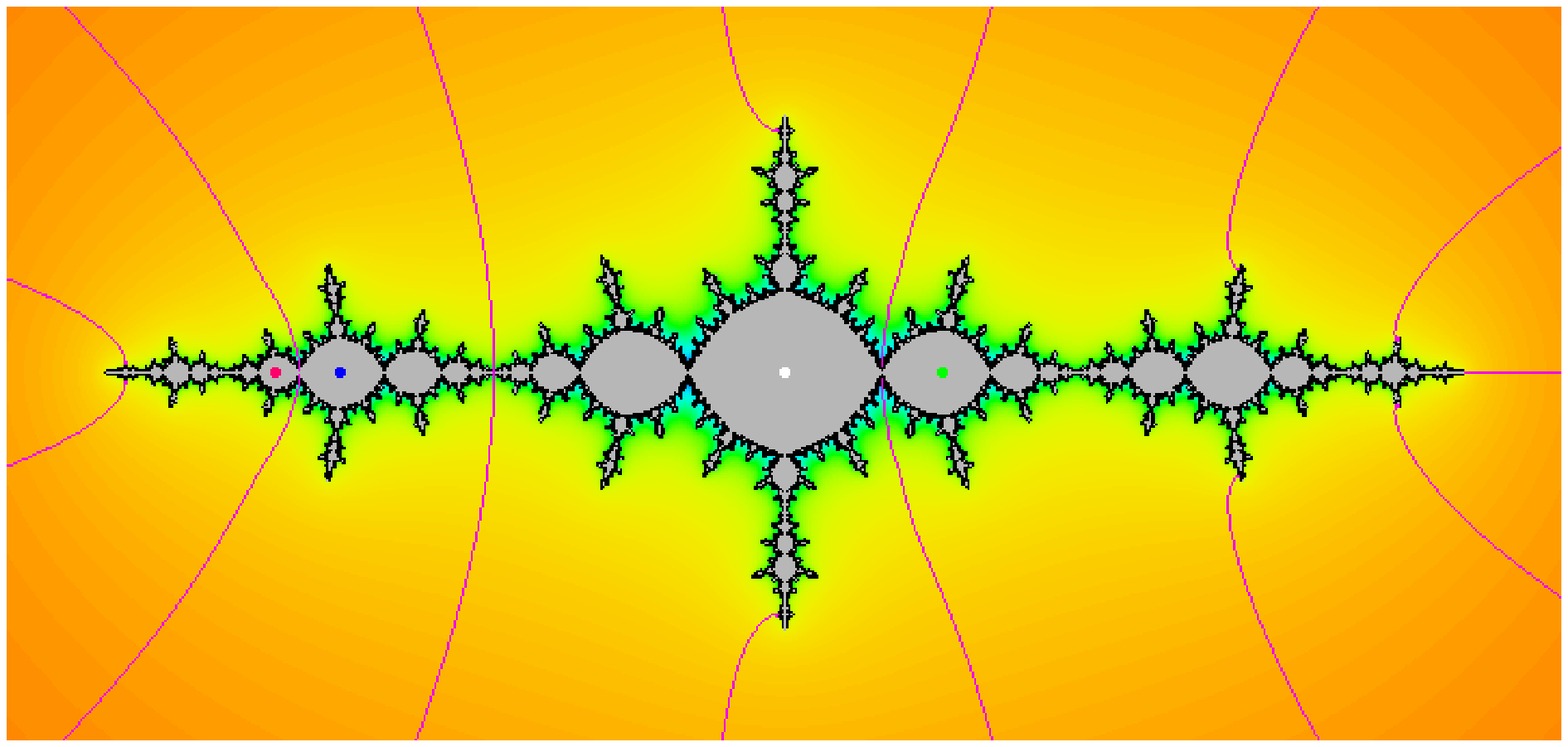,height=.7in}}
\caption{\textsl{ Four quadratic Julia sets with period $4$ critical orbit.}
\label{f-quad-j}}
\end{figure}
\medskip

\begin{lem}\label{l-k=0}
With $c_j$ and $\mathbb E$ as above, there are unique power series
$$ \ub_j~=~-\idt^2c_j\,+\,({\rm higher~order~terms})
~\in~{\mathbb E}[[\idt^2]] $$
which satisfy the required Equations $(\ref{e-u-condition})$. The 
corresponding escape regions are always $\I$-invariant.
\end{lem}
\smallskip

{\bf Proof Outline.} Suppose inductively that we can find series 
$\w_1,\ldots\w_{p-1}$
satisfying $~E_j(\vc{\w})\equiv 0~~({\rm mod}~\idt^{2n})~$. To start the
induction, for $n=3$ we can take $~\w_j=-\idt^2c_j$.

Setting $~\w'_j=\w_j+\idt^{2n-2}\delta_j~$, we find that
$$E_j(\vc{\w{}'})~\equiv~E_j(\vc{\w})\,+\,\idt^{2n}(\delta_{j+1}-\delta_1
+2c_j\delta_j)\quad({\rm mod}~\idt^{2n+2})\,.$$
(Compare Lemma \ref{l-err-diff}.) 
Now consider the linear map $~\vc{\delta}\,\mapsto\,\vc{\delta{}'}~$ from
${\mathbb E}^{p-1}$ to itself, given by
\begin{equation}\label{e-del'}
\delta{}'_j~=~\delta_{j+1}-\delta_1+2c_j\delta_j\,,
\end{equation}
where $\delta_p=0$. If this transformation has non-zero determinant,
then we can choose the $\delta_j$ so that $~E_j(\vc{\w{}'})~\equiv~0~~
({\rm mod}~\idt^{2n+2})$, thus completing the induction.

Now note that $c_1$ satisfies a monic polynomial equation with integer
coefficients, and hence is an algebraic integer. Thus the ring $\Z[c_1]\subset
{\mathbb E}$ is well behaved, with $~~1\not\equiv 0~~({\rm mod} ~2\,\Z[c_1])$.~
If we work modulo the ideal $2\,\Z[c_1]$, then the transformation of 
Equation~(\ref{e-del'}) reduces to
$$ \delta{}'_j~\equiv~\delta_{j+1}-\delta_1\,,$$
which is easily seen to have determinant one. Thus Equation~(\ref{e-del'}) also
has non-zero determinant, which completes the proof.\qed
\end{ex}

\begin{ex}[\bf Non-Primitive Regions]\label{ex-np}
More generally, consider any escape region $\cE_h\subset\cS_p$~ such that
the marked grid has period $n$ strictly less than the critical
 period $p$.  Branner and Hubbard showed that every such region can
be considered as a ``satellite'' of an escape region in $~\cS_n$,~ associated
with some quadratic map of critical period $p/n$. (Compare Theorem \ref{t-bh}.)
Thus it is natural to start
 with a period $n$ solution $\vc{\w}$ to the equations
$~E_j(\vc{\w})=0$,~ and look for a period $~p~$ solution of the form
 $~\vc \ub\;=\;\vc{\w}+\vc \delbold$.

\begin{theo}\label{t-np} 
Suppose that $~\vc\w~$ is a period $n$ solution to the equations\break
$~E_j(\vc\w)=0$,~ and that $~\vc\ub\,=\,\vc\w+\vc\delbold~$ is a period
$~p~$ solution, with $~p\,>\,n$.\break Assume that
 $~\delbold_j=o(\w_j)$~ for $~j\not\equiv 0~({\rm mod}~n)$,~
and that $~\ub_j=\delbold_j=o(1)~$ for  $~j\equiv0~({\rm mod}~n)$.~
Then for $~j\equiv0~({\rm mod}~n)$, the leading monomial for
 $~\ub_j=\delbold_j~$ has the form
\begin{equation}\label{e-np1}
 \m_{n\,i}~=~c_i\,\blam\,,
\end{equation}
where $~ 0~\mapsto ~c_1~\mapsto ~c_2~\mapsto~\cdots
 ~\mapsto ~c_{p/n-1}~\mapsto~0$~
is the period $~p/n~$ critical orbit for an associated quadratic map,~
and where $~\blam~$ is the fixed monomial
\begin{equation}\label{e-lam}
 \blam~=~-\bxi^{2n}/(\m_1^\bstar\,\m_2^\bstar\,\cdots\,\m_{n-1}^\bstar)\,.
\end{equation}\end{theo}


Here $~\m_j^\bstar=(3\sigma_j-2)\m_j~$ as in Equation (\ref{e-m*}).
Thus, given the leading monomials for the $\w_j$, and given the 
complex constant $c_1$, 
we can compute the leading monomials for the $~\ub_j$.~ According to
Theorem \ref{t-lmon}, the escape region $~\cE_h~$ is then  uniquely
determined. Note that
\begin{equation}
\m_j~=~\m_{n\,i+j}~\sim~\w_{n\,i+j}\qquad{\rm for}\qquad~ 0<j<n~\quad{\rm and}
\quad 0<n\,i<p\,, 
\end{equation}
but that
\begin{equation}\label{e-m=del}
 \m_{n\,i}~=~\bdel_{n\,i}\qquad{\rm with}\qquad \w_{n\,i}=0\,.
\end{equation}
\smallskip

The following identity is an immediate consequence of Theorem \ref{t-np},
and will be useful
in \S\ref{s-can-esc}. Since $~\m_n^\star=-2\m_n$,~
it follows that
\begin{equation}\label{e-prod=c1} 
 \frac{\m_1^\star}{\bxi^2}\cdots\frac{\m_n^\bstar}{\bxi^2}~=~ 2\,c_1\,.
\end{equation}\medskip

The proof of Theorem \ref{t-np} will depend on the following lemma.

\begin{lem}\label{l-kiwi}  
For an escape region in $\cS_p$ with grid period $~n$~ and for\break
 $0<k<p$,~ we have
 $$~\sum_{j=1}^k \big(2-\ord(\m_j)\big)~\ge ~0\,,$$
where this sum is is zero if and only if
  $~k~$ is congruent to zero modulo $~n\,$.\break
It follows that the holomorphic function $~\prod_{j=1}^k 
(\xi^2/m_j)$~ of the local parameter  $~\xi^{1/\mu}~$
tends to a well defined finite limit as the local parameter
 $~\xi^{1/\mu}~$ tends to zero,~ and that
this limit is non-zero if and only if $~k\equiv 0~({\rm mod}~n)$.
\end{lem}
\medskip

{\bf Proof.} 
With notation as in Equation (\ref{e-MOD}), let
$$ {\mathfrak S}_j~=~\sum_{\ell=0}^\infty\,{\rm MOD}\big(A_\ell(a_j)\big)\,.$$
(Intuitively, the number $~{\mathfrak S}_j~$ is ``large'' if and only if the
intersection of the puzzle pieces containing $~a_j~$ is ``small''.)
These numbers have two basic properties:\medskip

\begin{enumerate}
\item[{\bf(a)}] ~~$ {\mathfrak S}_{j+1}\,-\,{\mathfrak S}_j~=~\ord(\m_j)\;
-\;2\,.$
\medskip

\item[{\bf(b)}]~~ If the period of the grid is $~n$,~ then
$~ {\mathfrak S}_1~>~ {\mathfrak S}_j~$ for $~1<j\le n$.
\end{enumerate}
\medskip

\noindent To prove {\bf(a)}, note that
 $~{\rm MOD}\big(A_{\ell-1}(a_{j+1})\big)~$ is equal to either
$~2\;{\rm MOD}\big(A_\ell(a_j)\big)~$ or $~{\rm MOD}\big(A_\ell(a_j)\big)$,~
according as $~{\rm MOD}\big(A_\ell(a_j)\big)~$ does or does not surround
the critical point. Using Theorem \ref{t-kiwi}, it follows that the difference
 $~{\mathfrak S}_{j+1}\,-\,{\mathfrak S}_j$~ is equal to
$$\sum_{\ell=1}^{L_0(a_j)} {\rm MOD}\big(A_\ell\big)~=~\ord(\m_j)\,,$$
minus the term $~{\rm MOD}\big(A_0(a_j)\big)=2$.\medskip

To prove {\bf(b)}, note the formula
 $$ {\rm MOD}\big(A_\ell(a_j)\big)~=~2\;{\rm MOD}\big(A_{\ell+k}
(a_{j-k})\big)\,,$$
where $~k\le j~$ is the smallest positive integer such that 
$~A_{\ell+k}(a_{j-k})~$ surrounds the critical point. It follows that
$$ {\mathfrak S}_1~=~\sum_{\ell=1}^\infty
\,2\;{\rm MOD}\big(A_\ell(a_0)\big)\,,$$
while each $~{\mathfrak S}_j~$ with $~1<j<n~$ can be computed as a similar
 sum, but with some of the terms missing. The inequality {\bf(b)} follows.
\smallskip
Using both {\bf(a)} and {\bf(b)}, we have the inequality
$$ \sum_{j=1}^k\;\big(\ord(\m_j)-2\big)~=~{\mathfrak S}_{k+1}-{\mathfrak S}_1~<~0\qquad
{\rm for}\qquad 1\le k<n\,. $$
On the other hand, if $~n<p$,~ then since $~{\mathfrak S}_{n+1}={\mathfrak S}_1$,
the corresponding sum is zero. This proves Lemma \ref{l-kiwi} in the cases
when $~k\le n$.~ The remaining cases follow immediately since $~{\mathfrak S}_{j+n}
\,=\,{\mathfrak S}_j\,$.
\qed\medskip

\begin{rem}\label{r-mu}
One interesting consequence of this lemma, together with Theorem ~\ref{t-kiwi},
 is that multiplicity of an escape region depends only
on its marked grid, and is independent of the period $p$. In particular,
if a region in $~\cS_p~$ is a ``satellite'' of a region in $~\cS_n$~
with $~n<p$,~
then the two have the same multiplicity. (This statement would also
follow easily from the Branner-Hubbard Theory.)
\end{rem}\smallskip

{\bf Proof of Theorem \ref{t-np}.} Since
 $~E_j(\vc\ub)=0~$ and $~E_j(\vc{\w})=0$,~
we have  
\begin{eqnarray} \label{e-odd}
\idt^2(\delbold_{j+1}-\delbold_1)
&~\sim~&~\m_j^{\textstyle\star}\delbold_{j}\qquad \quad{\rm if}\quad j\not\equiv 0~{\rm mod}~n\,,\\
\label{e-even}
\idt^2(\delbold_{j+1}-\delbold_1)
&~\sim~&~-\delbold_{j}^{\,2}~~\,\qquad\;{\rm if}\quad j=n\,i\,,
\end{eqnarray}
where the first equation follows from Lemma~\ref{l-err-diff} and the second
follows from  Equation~(\ref{e-err}). 
We will prove by induction on $~j~$ that
\begin{equation}\label{e-ind}
\delbold_1~\sim~\frac{\bxi^2}{\m_1^\bstar}\cdots
\frac{\bxi^2}{\m_{j-1}^\bstar}\;\delbold_j
\end{equation}
for $~1\le j\le n\,$. In fact, given this statement for some $~j<n$,~
 it follows from  Lemma \ref{l-kiwi} that
$$ \frac{\bxi^2}{\m_j^\bstar}\;\delbold_1~=~o(\delbold_j)\,.$$
Equation (\ref{e-odd}) then implies that $~\bxi^2\,\delbold_{j+1}
~\sim~\m_j^\bstar\,\delbold_j$.~ Substituting this equation into (\ref{e-ind}),
the inductive statement that
$$\delbold_1~\sim~\frac{\bxi^2}{\m_1^\bstar}\cdots
\frac{\bxi^2}{\m_{j}^\bstar}\;\delbold_{j+1} $$
follows. This completes the proof of Equation (\ref{e-ind}).
\smallskip

It follows from this equation that
$$ \ord(\delbold_1)
~=~2(n-1)\,+\,\ord(\delbold_n)\,-\,\ord(\m_1\cdots\m_{n-1})\,.$$
On the other hand, it follows from Lemma \ref{l-kiwi} that
$$ \ord(\delbold_n)~=~\ord(\m_n)~=~2n\,-\,\ord(\m_1\cdots\m_{n-1})\,.$$
Combining these two equations, it follows that
$$\ord(\bxi^2\,\bdel_1)~=~\ord(\bdel_n^{\,2}).$$
On the other hand, according to (\ref{e-even}) we have
$$ \bxi^2(\delbold_{n+1}-\delbold_1)~\sim~ -\delbold_n^{\,2} \,$$
Using the ultrametric inequality (\ref{e-umi}), this implies that
$~\|\bdel_{n+1}\|\,\le\,\|\bdel_1\|$,~ or in other words
$~\ord(\bdel_{n+1})\,\ge\,\ord(\bdel_1)$.~
 The proof now divides into two cases.
\medskip

{\bf Case 1.}~~ If $\ord(\bdel_{n+1})=\ord(\bdel_1)$,~ then arguing as 
above we see that
$$ \bdel_{n+1}~\sim~\frac{\bxi^2}{\m_1^\bstar}\cdots
\frac{\bxi^2}{\m_{n-1}^\bstar}\;\delbold_{2n}\,. $$
again there is a dichotomy; but if $~\ord(\bdel_{2n+1})=\ord(\bdel_1)$,~
then we can continue inductively. However, this cannot go on forever
since $~\bdel_p=0$. Thus eventually we must reach:
\smallskip

{\bf Case 2.}~~ For some $~i\ge 1$,~ $\ord(\bdel_{n\,i+1})\;>\;
\ord(\bdel_1)$,~
or in other words $~\|\bdel_{n\,i+1}\|~<\|\bdel_1\|$.~ Then a
 straightforward induction, based on the statement that
 $~\bxi^2(\bdel_{n\,i+j+1}-\bdel_1)~\sim~\m_j^\bstar\bdel_{n\,i+j}$,~ shows that
$$ \|\bdel_{n\,i+j}\|~<~\|\bdel_j\|~=~
\left\|\frac{\m_1^\bstar}{\bxi^2}\cdots
\frac{\m_{j-1}^\bstar}{\bxi^2}\bdel_1\right\|\,.$$
Taking $~j=n$,~ since the $~\|\bdel_{n\,i}\|~$ with $~0<n\,i<p~$
all have the same norm by Theorem \ref{t-kiwi}, it follows that $~n\,i+n=p$,~
with $~\bdel_p=0$.\smallskip

Now define the monomial $~\blam\in\C[\bxi^{1/\mu}]$~ by Equation
 (\ref{e-lam}), and define $~c_i$~ by Equation (\ref{e-np1}).
Then it is not hard to check that the $c_i$ are complex numbers.
Now note that
$$ \bxi^2\,\bdel_{n\,i+1}~\sim~\bxi^2\,\frac{\bxi^2}{\m_1^\bstar}\cdots
\frac{\bxi^2}{\m_{n-1}^\bstar}\;\bdel_{n\,i+n}~\sim~\blam^2\,c_{i+1}
\quad{\rm for}\quad n\,i<p-n\,,$$
while
$$\bxi^2\,\bdel_{p-n+1}~=~o(\blam^{\;2})\,.$$
Using the equation
$$~\bxi^2(\bdel_{n\,i+1}-\bdel_1)~\sim~-\bdel_{n\,i}^{\;2}~\sim~
-\blam^2\,c_i^{\,2}\,,$$
it follows that 
$$ c_{i+1}~=~c_i^{\;2}\,+\,c_1\,,$$
with $~c_{p/n}\,=\,0\,$.
This completes the proof of Theorem \ref{t-np}.\qed
\end{ex}\smallskip

\begin{rem}[{\bf An Empirical Algorithm}]\label{r-empir} 
How can we locate examples of maps $~F_{a,\,v}$~ 
belonging to some unknown escape 
region? Choosing an arbitrary value of $~a~$ (with $~|a|~$ not too
small), given the kneading sequence, and given a rough approximation
$~v'\,\approx\,~v$,~
the following algorithm is supposed to converge to the precise value
 of $~v$.~ It seems quite useful in practice.
(In fact, it often eventually
 converges, starting with a completely random  $~v'$.)~
However, we have no proof of convergence, even if the initial $~v'~$ 
is very close to the correct value. Furthermore, there are
 examples (near the boundary of an escape region) where the algorithm converges
 to a map with the wrong kneading sequence.

Given the required kneading sequence $~\vc\sigma$,~ consider the generically
defined holomorphic maps $~\Psi_j\,:\,\C^3~\to~\C~$ defined by
$$ \Psi_j(w_1,\,w_j,\, w_{j+1})~=~\begin{cases}
w_j\sqrt{\frac{\xi^2(w_{j+1}-w_1)}{w_j^{\;2}(w_j-1)}} & {\rm if}\quad \sigma_j=0\,,\\ \\
1\;+\;\frac{\xi^2(w_{j+1}-w_1)}{w_j^{\,2}} & {\rm if}\quad \sigma_j=1\,,
\end{cases}$$
taking that branch of the square root function which is defined on the
right half-plane. Then it is easy to check that
$$ w_j~=~\Psi_j(w_1,\,w_j,\,w_{j+1})\qquad{\rm if~and~only~if}\qquad
E_j(\vc w)=0\,.$$
Choosing some $~\xi=1/(3a)$,~ and
 given some approximate solution to the equations $~E_j(\vc w)=0$,~
proceed as follows. Replace $~w_j~$ by $~\Psi_j(w_1,\,w_j,\,w_{j+1})$,
starting with $~j=p-1$,~ then with $~j=p-2$,~ and so on until $~j=1$.
Then repeat this cycle until the sequence has converged to machine accuracy,
or until you lose patience. If the sequence does converge, then  
the limit will certainly describe some map in  $~\cS_p$,~ and it is easy
 to check whether or not it has the required kneading sequence.\smallskip

If we start only with a pair $~(a,\,v)$~ and a kneading sequence,
 we can first compute the partial orbit
$$ a~=~a_0~\mapsto~a_1~\mapsto~ \cdots~\mapsto a_{p-1}\,,$$
then set $~w_j=(a-a_j)/(3a),$~ and proceed as above.
\end{rem}
\medskip

\setcounter{table}{0}
\section{Canonical Parameters and Escape Regions.}\label{s-can-esc}

Recall from \S\ref{s-can} that the canonical coordinate $~t=t(a,v)~$ on $~\cS_p~$
 is defined locally, up to an additive constant, by the formula
\begin{equation}\label{e-dt}
dt~=~\frac{da}{\partial\Phi_p/\partial v}\qquad{\rm where}\qquad
\Phi_p(a,v)~=~F_{a,v}^{\circ p}(a)-a\,.
\end{equation}
In particular, the differential $~dt~$ is uniquely defined, holomorphic,
and non-zero everywhere on the curve $\cS_p$. We will prove that
the residue
\begin{equation}\label{e-res}
\frac{1}{2\pi i}\oint dt\end{equation}
at the ideal point 
is zero, so that the indefinite integral $~t\;=\;\int dt~$ is
well defined, up to an additive constant, throughout the escape region.
(However, as noted in Remark \ref{r-can-glob}, it usually cannot be defined
as a global function on $~\cS_p$.)\smallskip

Furthermore, whenever the kneading sequence is non-trivial, we will show that
the function $~t~$ has a removable singularity, and hence extends to a
smooth holomorphic function which is defined and finite 
 also at the ideal point.
On the other hand, in the case of a trivial kneading sequence 
 $~000\cdots0$,~ the function $~t~$
 has a pole at the ideal point, and in fact is asymptotic to some constant
multiple of $a$. 
We will first need the following.
\medskip

\begin{definition}\label{d-psi}
For each $~j\ge 0$,~ consider the polynomial expression
\begin{eqnarray*}
\psi_{j+1}(X_1,\,\ldots,\,X_j)~=~\hspace{8cm}\\
 X_1X_2X_3X_4\cdots X_j~+~X_2X_3X_4\cdots X_j~+~X_3X_4\cdots X_j
~+~\ldots~ + X_j~+~1
\end{eqnarray*}
in $j$ variables. These can be defined inductively by the formula
 $$~~\psi_{j+1}(X_1,\,\ldots,\,X_j)~=~
 \psi_j(X_1,\,\ldots,\,X_{j-1})\,X_j\,+\,1\,,$$ 
starting with $~\psi_1(~)~=~1$.
\end{definition}

 We will prove the following
statement. It will be convenient to make use of both of the variables
$~a~$ and $~\xi=1/(3a)$. As in Remark \ref{r-anch}, a choice of anchoring
for $~\cE_h~$ determines a choice of local parameter
 $~\xi^{1/\mu}\,=\,1/(3a)^{1/\mu}~$ near the ideal point.

\begin{theo}\label{t-dt} 
Let $~n\le p$~ be the grid period for the escape region\break
$~\cE_h\subset \cS_p\,$,~ and let
 ~~$0\mapsto c_1\mapsto c_2\mapsto\cdots\mapsto c_{r-1}\mapsto 0$~~
 be the critical orbit for the
associated quadratic map, where $~r=p/n\ge 1$.~
Then the derivative  $~da/dt~$ is given by the asymptotic formula
$$da/dt~\sim~\psi_r(2c_1,\ldots,2c_{r-1})\,
\prod_{j=1}^{n-1}\,\frac{m_j^{\textstyle\star}}{\xi^2}\qquad{\rm
as}\qquad \xi\to 0\quad({\rm or~as}\quad |a|\to\infty)\,,$$
where $~m_j^{\textstyle\star}=(3\sigma_j-2)m_j$~ as in Equation $(\ref{e-m*})$,
and where $~\psi_r(2c_1,\cdots,2c_{r-1})$ is a non-zero complex number.
\end{theo}

Note that
$$ d\xi/dt~=~(d\xi/da)(da/dt)~=~-3\xi^2 \,da/dt\,.$$
Thus it follows from Theorem~\ref{t-dt} that $~t~$ can be expressed
 as an indefinite integral 
\begin{eqnarray}
t~=~\int dt  &\sim &\int\frac{da}{ \psi_r(2c_1,\ldots, 2c_{r-1})}
\prod_{j=1}^{n-1}\frac{\xi^2}{m_j^\bstar}~~x \nonumber \\
&~~=~& \int\frac{-d\xi~}{~3\,\xi^2\,
 \psi_r(2c_1,\ldots, 2c_{r-1})}\prod_{j=1}^{n-1}\frac{\xi^2}{m_j^\bstar}\;.
\label{e-t-sim}
\end{eqnarray}
If the kneading sequence is non-trivial, or in other words if $~n>1$,~
then Lemma~\ref{l-kiwi} implies that the product
$~\prod(\xi^2/m^\bstar_j)~$
is a bounded holomorphic function of the local parameter
$~\xi^{1/\mu}$~ in a neighborhood of the ideal point,
 and hence that $~t~$
can be defined as a bounded holomorphic function of $\xi^{1/\mu}$.~
On the other hand, if the
kneading sequence is trivial so that $~n=1$,~ then this integral has
a pole at the ideal point, and in fact it follows easily that
\begin{equation}\label{e-t-sim2}
 t~~\sim~~\frac{a}{\psi_p(2c_1,\ldots, 2c_{p-1})}~.
\end{equation}

If we choose the constant of integration in the indefinite integral of
Equation~(\ref{e-t-sim})  appropriately, then we can express $t$ as
 a holomorphic function of the form
\begin{equation}\label{e-puis-can}
 t(\xi^{1/\mu}) ~=~\beta\, \xi^{\nu/\mu}
\;+\;\rm{(higher~order~terms)}\;,
\end{equation}
where $~\beta~$ is a non-zero complex coefficient and $~\nu~$
 is a non-zero integer. 

\begin{definition}\label{d-ind-int}
This integer
 $~\nu\ne 0~$ will be called the \textbf{\textit{winding number}}\break (or the
\textbf{\textit{~ramification index}})~ of the escape region over the
 $~t$-plane. 
In fact,\break
as we make a simple loop in the positive
 direction in the $\xi^{1/\mu}$-parameter plane,  the corresponding point
 in the $~t$-plane will wind $~\nu~$ times around the origin. Thus
 $~\nu>0~$ when the kneading sequence $~\vc\sigma~$
is non-trivial;~ but $~\nu=-1$~ when $~\vc{\sigma}\;=\;(0,\ldots,0)$.~
 In either case, it follows that any choice of
$~t^{1/\nu}~$ can be used as a local uniformizing parameter near the ideal
point.
It follows easily from Equation (\ref{e-t-sim}) that this winding number
is given by the formula
\begin{equation}\label{e-nu=} 
 \nu~=~\big(2n\,-\,3\,-\,\ord(m_1\cdots m_{n-1})\big)\,\mu\,.  
\end{equation}
Using Theorem \ref{t-kiwi}, it follows that $~\nu~$ depends only on
 the marked grid.
\end{definition}\smallskip
\medskip

{\bf Proof of Theorem \ref{t-dt}.} The argument
 will be based on the following
 computations. Start with the periodic orbit
 ~~$ a=a_0\mapsto a_1\mapsto\cdots$,~ where
$$  a_{j+1}=F_{a,v}(a_j)~=~ a_j^3-3a^2a_j+2a^3+v\,.$$
It will be convenient to set
\begin{eqnarray}\label{e-Xj}
X_j&~=~&\frac{\partial a_{j+1}}{\partial a_j}~=~3(a_j^2-a^2)\,,\\
\label{e-paj}
Y_{j+1}&~=~&\frac{\partial a_{j+1}}{\partial v}~=
~\frac{\partial a_{j+1}}{\partial a_j}\,\frac{\partial a_j}{\partial v}
\,+\,1\,.
\end{eqnarray}
In other words
\begin{equation}\label{e-Y}
Y_{j+1}~=~Y_jX_j\,+\,1\qquad{\rm with}\qquad Y_1~=~
 \frac{\partial a_1}{\partial v}~=~1\,,
\end{equation}
or briefly
\begin{equation}\label{e-psi}
  Y_{j+1}~=~ \psi_{j+1}(X_1,\,\ldots,\,X_j)\,,
\end{equation}
using Definition \ref{d-psi}.\smallskip

Denoting $~~\Phi_p(a,v)\;=\;F_{a,v}^{\circ p}-a~~$
 briefly by $~a_p-a$,~ and noting that\break
 $~\partial a/\partial v=0$,~
we see from Equation (\ref{e-dt}) that
$$ \frac{\partial\Phi_p(a,v)}{\partial v}~=~\frac{\partial a_p}{\partial v}
~=~Y_p~=~\psi_p(X_1,\,\ldots,\,X_{p-1})\,, $$
and hence that
$$ dt~=~\frac{da}{\partial\Phi_p/\partial
  v}~=~\frac{da}{Y_p}\qquad{\rm or}
\qquad da/dt~=~Y_p\,.$$
Thus, in order to prove Theorem \ref{t-dt}, we must find an asymptotic
 expression for $~Y_p$.\smallskip





Substituting $~a=1/3\,\xi~$ and $~a_j=(1-3u_j)a$~ in Equation~(\ref{e-Xj}),~
 we easily obtain
\begin{equation}\label{e-X-as}
 X_j~=~\frac{3u_j^{\,2}-2u_j}{\xi^2}~\sim~\frac{m_j^{\textstyle\star}}{\xi^2}~,
\end{equation}
where $~m_j^{\textstyle\star}=(3\sigma_j-2)m_j~$ is the leading term in the
 series expansion for
 $~3u_j^{\,2}-2u_j~$ as described in Equation~(\ref{e-m*}),~
so that $~m_j^{\textstyle\star}\sim 3u_j^{\,2}-2u_j$.
We next show by induction on $~j~$ that
$$ Y_j~\sim~X_1\cdots X_{j-1}~\sim~\frac{m_1^\bstar}{\xi^2}\cdots
\frac{m_{j-1}^\bstar}{\xi^2}\qquad{\rm for}\qquad j\le n\;. $$
In fact this expression for a given value of $~j~$ implies that
$$ X_jY_j~\sim~X_1X_2\cdots X_j\,,$$
where this expression has a pole at the origin by Lemma \ref{l-kiwi}.
Hence it is asymptotic to $~X_jY_j+1=Y_{j+1}$,~ which completes the
induction. In particular, it follows that
 $~Y_n\;\sim\;\prod\limits_{j=1}^{n-1}\displaystyle{\big(\frac{m_j^\bstar}{
\xi^2}\big)}$.~ In the primitive case
 $~p=n$,~ this completes the proof of Theorem \ref{t-dt}.\smallskip

Suppose then that $~p>n$.~ By Equations
(\ref{e-np1}) and (\ref{e-lam}) of Theorem \ref{t-np}, we have
$~m_n^\bstar\;=\; -2\,c_1$,~ and hence
$$ X_{n}Y_{n}~\sim~
(m_1^\bstar/\xi^2)\cdots(m_n^\bstar/\xi^2)~=~ 2c_1\,,$$
and therefore $~Y_{n+1}~\sim~2\,c_1+1~=~\psi_2(2\,c_1)\,.$
Suppose inductively that
$$ Y_{n\,i+1}~\sim~\psi_{i+1}(2\,c_1,\ldots,2\,c_i)\,.$$
Then by induction on $~j~$ we see that 
$$ Y_{n\,i+j}~\sim~\psi_{i+1}(2\,c_1,\cdots,2\,c_i)\,
X_1\cdots  X_{j-1}~\sim~\psi_{i+1}(2\,c_1,\cdots,2\,c_i)\,
\frac{m_1^\bstar}{\xi^2}\cdots
\frac{m_{j-1}^\bstar}{\xi^2} $$
 for $~ j\le n\;.~$ Taking $~j=n~$ and $~i=r-1$~ so that $~ni+j=p$,~
this completes the proof.\qed\medskip

\begin{table} [ht] 
\begin{center}
\begin{tabular}{cclcccccccc}
$p$ & $r$ & $\vc{\sigma}$ & $m_1$ & $m_2$ & $m_3$ & $t\sim$ & $\nu$
& $\mu$ & $\#$ & sym\\ \\
2 &1& $\overline{10}$ & 1 & & & $-\xi/3$ &1&1&1 &$\pm$\\
3 &1& $\overline{110}$ & 1 & 1 & &$-\xi^3/9$ &3&1&1&$\pm$\\  
3 &1& $\overline{100}$ & 1 & $\pm\xi$ & & $\pm \xi^2/12$&2&1&2&$+$ \\
3 &1& $\overline{010}$ & $\pm i\xi$ & 1& & $\mp i\xi^2/12$& 2&1& 2&$-$\\
4&2 & $\overline{1010}$ & 1 & $\xi^4$ & 1 & $\xi/3$ & 1&1 & 1 & $\pm$\\ 
4 &1& $\overline{1110}$ & 1 & 1 & 1 & $-\xi^5/15$&5&1&1&$\pm$\\
4 &1& $\overline{1100}$ & 1 & 1 & $\pm \xi$ & $\pm \xi^4/24$&4&1&2&$+$\\
4 & 1&$\overline{0110}$ & $\pm i\xi$& 1 & 1& $\mp i \xi^4/24$&4&1&2&$-$\\
4 &1& $\overline{1000}\,{\rm s}$ & 1 & $\pm \xi$ & $\mp \xi$ & $\;\;\xi^3/36$&3&1&2
&$+$\\
4 &1& $\overline{1000}\,{\rm t}$ & 1 & $\pm \xi$ & $\pm \xi$ & $-\xi^3/36$&3&1&2&$+$\\
4 &1& $\overline{0100}$ & $\omega^2\xi$ & 1 & $-\omega \xi^{3/2}$ & $-\omega
\xi^{5/2}/30$& 5&2&2&$-$\\
4 &1& $\overline{0010}$ & $\omega \xi^{3/2}$ & $-\omega^2\xi$ & 1 & $-\omega
\xi^{5/2}/30$&5&2&2&$-$\\
\end{tabular}
\end{center}
\vspace{-.3cm}
\caption{\textsl{Escape regions in $~\cS_p~$
with nontrivial kneading  sequence for $p \leq 4$.
Here $~r=p/n~$ is the critical period for the associated quadratic map, and
$\#$ is the number of such regions. In the last column, $+$ stands
 for symmetry under complex conjugation
so that $~\cE_h=\overline\cE_h$,
and $-$ stands for symmetry under $\I$ composed with complex conjugation.
 $($Compare \S$\ref{s-real}.)$ The symbol $~\omega~$ stands for any
 $4$-th root of $-1$. Note that
the exponent of $~\xi~$ in the $~t\!\sim~$ column, 
is equal to $\nu/\mu$. In all of these cases, $\nu$ and $\mu$ are
relatively prime.}
\label{t-4}}
\end{table}

\setcounter{table}{0}
\section{The Euler Characteristic}\label{s-Eul}

Using the meromorphic 2-form $~dt$,~ it is not difficult to calculate the
Euler characteristic of $~\cS_p~$ or of $~\overline\cS_p$. Here is a preliminary
result.

\begin{lem}\label{l-mer} 
The Euler characteristic of the compact curve
$~\overline\cS_p~$ can be expressed as
$$ \chi(\overline\cS_p)~=\sum_{\rm escape~regions} (1-\nu)
\,,$$
to be summed over all escape regions in $~\cS_p$.~ Here $~\nu$~ is
 the winding number of Definition $\ref{d-ind-int}$,~ as computed in
Equation $(\ref{e-nu=})$. It follows that the Euler characteristic of
 the affine curve $~\cS_p~$ is given by
$$ \chi(\cS_p) ~=~
-\sum \nu\,.$$
\end{lem}

{\bf Note:} By definition, the Euler characteristic $~\chi(\cS_p)~$
is equal to the difference of Betti numbers, 
$~{\rm rank}\big(H^0(\cS_p)\big)\,-\,{\rm rank}\big(H^1(\cS_p)\big)$.~
If we use {\v C}ech cohomology, then the Betti numbers of $~\cS_p~$
can be identified with the Betti numbers of the connectedness locus
$~\cC(\cS_p)$.\medskip

{\bf Proof of Lemma \ref{l-mer}.}
Given any meromorphic 1-form on a compact (not necessarily connected)
Riemann surface $~\cS$,~
we can compute the Euler characteristic of $~\cS$~
in terms of the zeros and poles
of the 1-form counted with multiplicity, by the formula
$$\chi(\cS)~~=~~ \#({\rm poles})\,-\,\#({\rm zeros})\,.$$
We can apply this formula to the meromorphic 1-form $~dt~$ on the
compact curve $\overline\cS_p$.~ Since the only zeros and poles of $~dt~$
are at the ideal points, we need only compute the order of the zero or
pole at each ideal point, using $~\xi^{1/\mu}$~ as local uniformizing parameter.
In fact, using Theorem \ref{t-dt}, we see easily that $~dt~$ has a zero of
order $\nu-1\ge 0$ whenever the kneading sequence is non-trivial,
and a pole of order $~1-\nu=2~$ when the kneading sequence is trivial.
This completes the proof for $~\overline\cS_p$.
The corresponding formula for $~\cS_p$~ follows, since each point which is
removed decreases $~\chi~$ by one.\qed\smallskip

{\bf Note:} Based on this result,  combined with a
combinatorial method to compute the number of escape regions
 associated with each marked grid,
Aaron Schiff at the University of Illinois Chicago,
 under the supervision of Laura DeMarco, has
 written a program to compute the Euler characteristic of
 $\overline\cS_p$ for many values of $~p$.
\bigskip

 The formula for $~\chi~$ can be simplified as follows. 
Let $~d_p~$ be the degree
of the affine curve $\cS_p$, defined inductively by the formula
$$ \sum_{n|p} d_n~=~3^{p-1}\,.$$

\begin{theo}\label{t-eul}
The Euler characteristic of the affine curve $~\cS_p$ is given by 
$$ \chi(\cS_p)~=~ (2-p)\,d_p\,.$$
Hence the Euler characteristic of $~\overline{\cS_p}$~ is
$$ \chi(\cS_p)~=~N_p\,+\,(2-p)\,d_p\,,$$
where $N_p$ is the number of escape regions $($ = number of puncture points$)$.
\end{theo}
\smallskip

As examples, for $~p\le 4$~ we have the following table.
\medskip

\begin{table}[h]
\begin{center}
\begin{tabular}{llccc}
$p$ &  $d_p$ & $(2-p)\times d_p$ &$~~N_p$~~& $\chi(\overline\cS_p)$\\
1 & 1 & $1\times 1$& 1&2\\
2 & 2 & $0\times 2$& 2 & 2\\
3 & 8 & $-1\times 8$&8& 0\\
4 &  24 & $-2\times 24$& 20 &-28\\
\end{tabular}
\end{center}
\vspace{-.3cm}
\caption{\label{tab-euler}}
\end{table}

The proof will depend on the following.
Let $~\widehat a\in\C~$ be an arbitrary
constant. Then the line $~\{(a,v)~;~ a=\widehat a\}~$ intersects $~\cS_p~$ in~
$~d_p~$ points, counted with multiplicity. 
The corresponding values of $~v~$ will be denoted by 
$~v_1,\,\ldots,\,v_{d_p}$.

\begin{lem}\label{l-sym-prod}
For each $~0<j<p$,~ the
product\begin{equation}\label{e-sym-prod}
 \prod_{k=1}^{d_p}~ 
\Big(\widehat a-F_{\widehat a,v_k}^{\circ j}(\widehat a)\Big)
\end{equation}
is a non-zero complex constant $($conjecturally equal to $+1\,)$.
\end{lem}

{\bf Proof.} The numbers $~v_k$~ are the roots of a polynomial equation
with coefficients in $~\C[\widehat{a}]$.~ Since Equation~(\ref{e-sym-prod})
is a symmetric polynomial function of the $~v_k\,$,~ it can be expressed as
a polynomial in the elementary symmetric functions, and hence can be
expressed as a polynomial function of $~\widehat{a}$.~ Since this function has
no zeros, it is a non-zero constant.\qed\medskip

{\bf Proof of Theorem \ref{t-eul}.} Using Lemma \ref{l-mer},
it suffices to prove the identity
$$ \sum_{h=1}^{N_p} \nu(\cE_h)~=~(p-2)\,d_p\,.$$
 By Equation~(\ref{e-nu=}) 
  together with Lemma~\ref{l-kiwi}, we have
$$\nu/\mu~=~ -1+\sum_1^{n-1}\Big(2-{\rm ord}(\m_j)\Big)~=~-1+\sum_1^{p-1} 
\Big(2-{\rm ord}(\m_j)\Big)\,.$$
Since $~~2-{\rm ord}(\m_j)\,=\,1-{\rm ord}(\ab-\ab_j)$,~ it follows that
\begin{equation}\label{e-numu}
 \nu/\mu~=~ p-2\,-\,{\rm ord}\Big((\ab-\ab_1)(\ab-\ab_2)\cdots(\ab-\ab_{p-1})\Big)\,.
\end{equation}
Choosing some large $\widehat{a}$ and summing
over the $d$ intersection points, each escape region contributes $\mu$
copies of Equation $(\ref{e-numu})$. Therefore the left side of the equation
 adds up to $\sum_1^{N_p} \nu$.
The first term on the right side adds up to\break
 $(p-2)d_p$.~
Since it follows from
Lemma~\ref{l-sym-prod}
that the second term on the right adds up to zero, this proves
 Theorem~\ref{t-eul}.\qed
\medskip

\setcounter{table}{0}
\section{Topology of $\cS_p$}\label{ss-top}

Conjecturally the curve  $\cS_p$ is irreducible (= connected) for all
periods $~p$.
Whenever this curve is connected, we can use the results of \S\ref{s-Eul}
to compute the genus of the closed Riemann surface
$~\overline\cS_p$,~ and the first Betti number of the open surface
$~\cS_p$~ or of its connectedness locus. However, unfortunately we do not
have a proof of connectivity for any period $~p>4$.
(In this connection, it 
is interesting to note that any connected component of $\cS_p$ for $p>1$
must contain hyperbolic components of Types A, B, C and D, as well as escape
components. In fact the zeros and poles of the meromorphic function
 $~(a,v)\mapsto a\in\C~$ yield centers of Type A and ideal points, while the
 zeros of the functions $~v+a~$ and $~F(-a)+a~$ yield centers of Type B and D.
With a little more work, one can find centers of Type C as well.)\medskip

\begin{figure}[ht!]
\centerline{\makebox[0pt][l]{\includegraphics[height=3in]{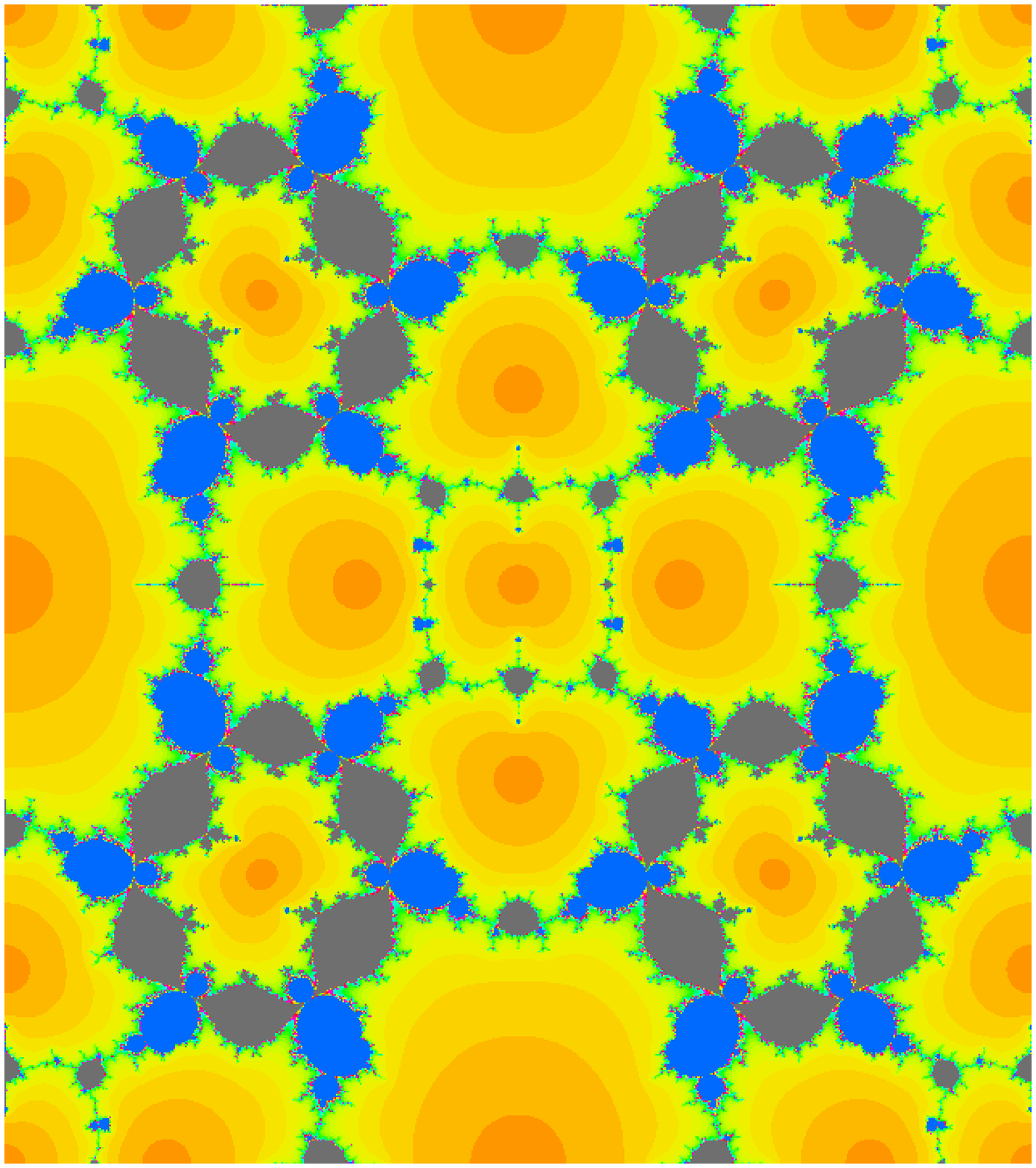}}%
  \includegraphics[height=3in]{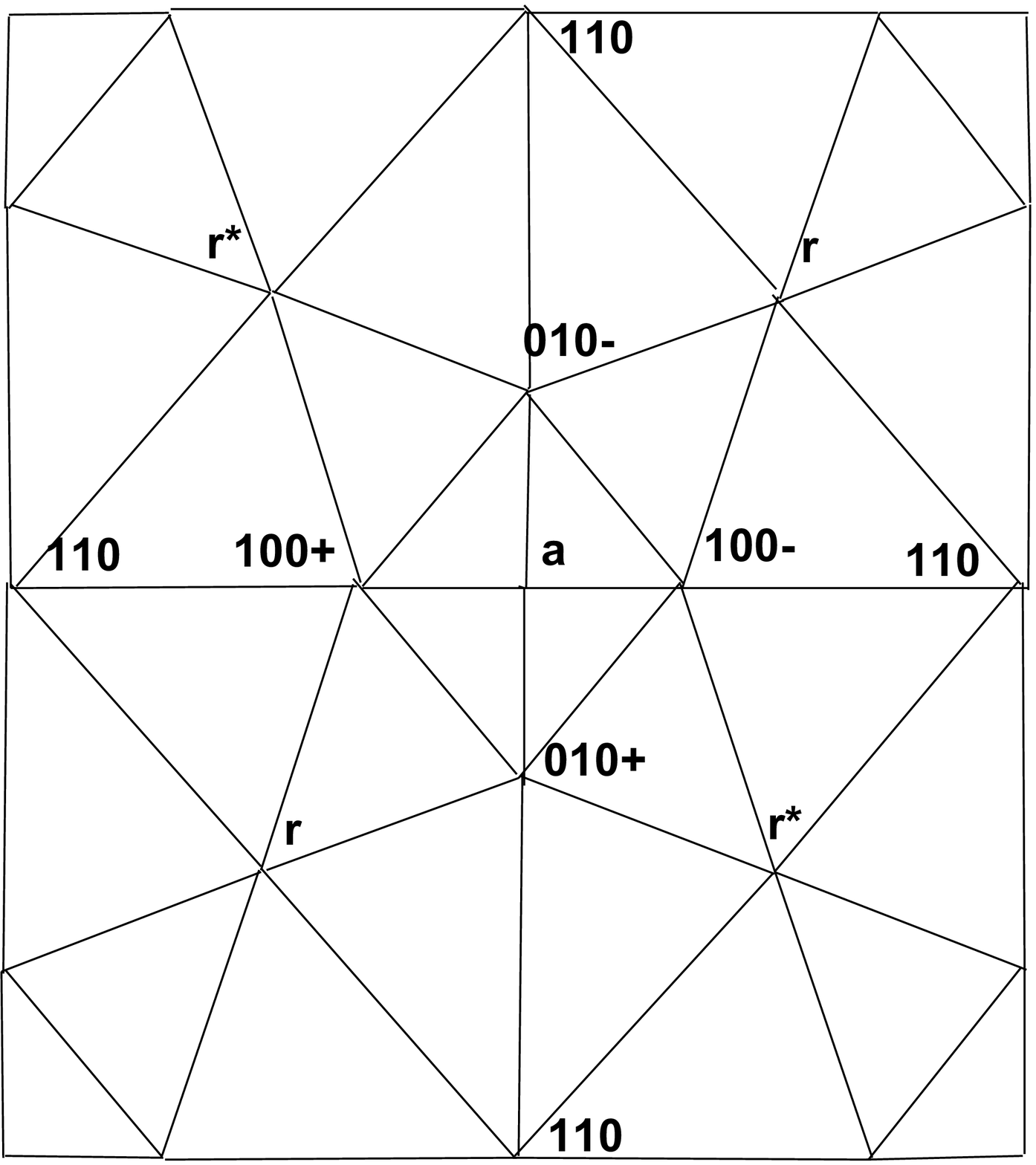}}
\caption{\textsl{
Universal covering space of $\overline\cS_3$. Here $\bf r$ stands for the
ideal point in the $(1/3)$-rabbit region, $\bf r^*$ for its complex conjugate,
and $\bf a$  for the ideal point in the airplane region. Note that the $180^\circ$
 rotation $\I$ interchanges\protect\footnotemark[9]
 $\overline{100}+$ and $\overline{100}-$, and also interchanges 
$\overline{010}+$  and $\overline{010}-$, but fixes the remaining four ideal
 points.}\label{f-3tor}}
\end{figure}

For $p\le 4$, we are
able to prove connectivity, and to provide a more direct computation of
 the genus, by constructing an
explicit cell subdivision of $\overline\cS_p$ with the ideal points as
vertices.\footnote{Compare 
  \cite[\S5D]{M4} which describes a dual cell structure,
 with one 2-cell corresponding to each ideal point.}
As an example, Figure~\ref{f-3tor} shows part of the universal covering
space of the torus $~\overline\cS_3$.~ The ideal point in each escape region
has been marked with a black dot. 
Joining neighboring ideal points by more or less straight lines, we can easily
 cut the entire plane up into triangular 2-cells. The corresponding
triangulation of $~\overline\cS_3~$ itself has 16 triangles, 24 edges,
and 8 vertices.

\footnotetext[9]{We can distinguish the two regions
$\overline{100}\pm$ by the convention that the difference $a_2-a_0$ between 
marked orbit points converges to $\pm 1$ as $|a|\to\infty$.
 Similarly $a_1-a_0\sim\pm i$ in the $\overline{010}\pm$ region.
(Compare Remark 4.3.)}

This cell structure can be described as follows.
 Let $~\cS~$ stand either for the curve
 $~\overline\cS_p~$ or for the quotient curve $~\overline\cS_p/\I$,~ where
$~\I:F_{a,v}\mapsto F_{-a,-v}~$ is the canonical involution.
Let  $~V~$ be the finite set consisting of all ideal points in $~\cS$,
 and for each $~\vb\in V~$ let $~\cE_\vb\subset\cS$~
be the associated escape region.
Given two points $~\vb,\,\vb'\in V$,~
consider all paths $P\subset\cS$ from $\vb$ to $\vb'$ which satisfy the
 following two conditions:

\begin{itemize}
\item[$\bullet$] $~P~$ is disjoint from the closure $\overline\cE_\w$ 
for all $\w$ other than $\vb$ and $\vb'$.

\item[$\bullet$] The intersections  $P\cap\cE_\vb$ and $P\cap\cE_{\vb'}$
are both connected.
\end{itemize}
\noindent
Choose one representative from each homotopy class of such paths $P$
as an edge from $\vb$ to $\vb'$. Conjecturally, these vertices and edges
form the 1-skeleton of the required cell subdivision. In any case, this
certainly works for periods $p\le 4$.                            

\begin{figure}[ht]
\centerline{\makebox[0pt][l]{\scalebox{0.75}{\includegraphics[width=\textwidth]{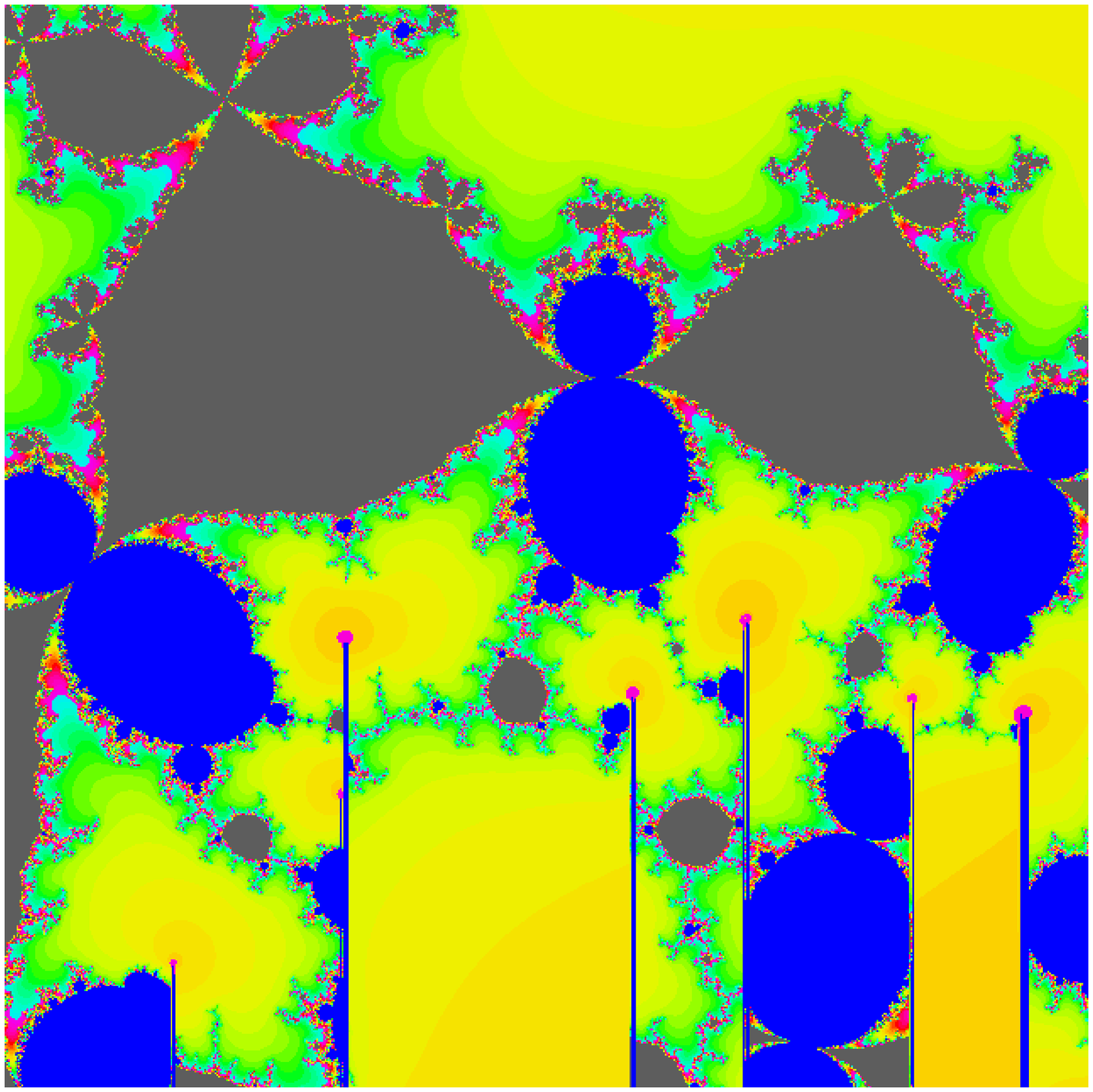}}}%
  {\scalebox{0.75}{\includegraphics[width=\textwidth]{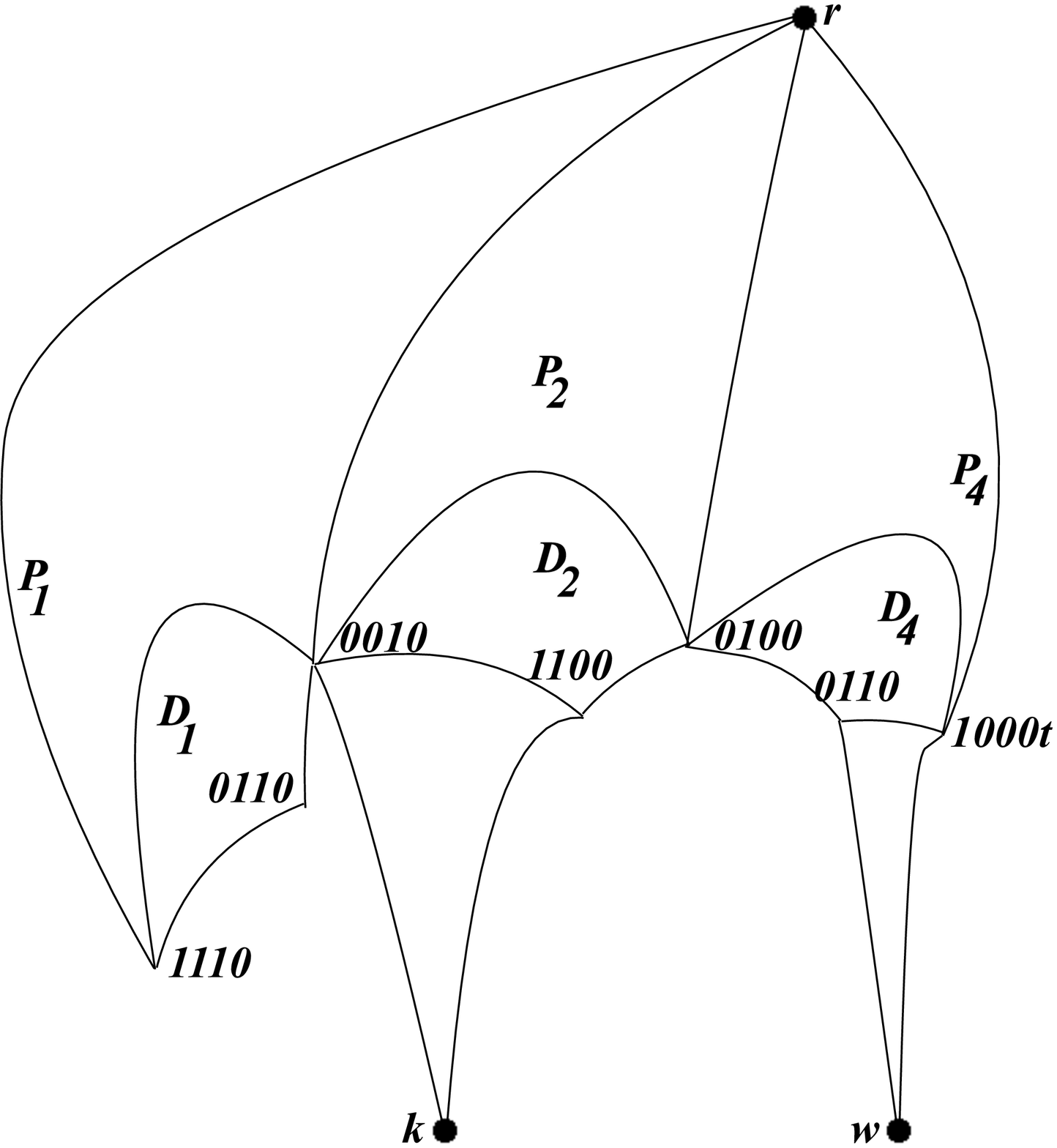}}}}
\caption{\textsl{Part of $\overline\cS_4$, 
 illustrating the cell structure. Here $k,\,r,\,w$ stand for the
kokopelli, $(1/4)$-rabbit, and worm centers. The three top triangles
are centered at parabolic points $P_i$, the next three at Type D centers $D_i$,
 and the bottom two at Type A centers.\protect\footnotemark[10]}\label{f-D2}}
\end{figure}
\footnotetext[10]{By definition, a hyperbolic component $H$ in the
  connectedness locus has 
 \textbf{\textsl{Type A, B, C,}} or \textbf{\textsl{D}}
according as the free critical point $-a$ for a map $F\in H$
either {\bf(A)} belongs to the same Fatou component as the marked
 critical point, {\bf(B)} belongs to a different component but in
 the same cycle of Fatou components, {\bf(C)} has orbit which lands in
 this cycle of Fatou components only after iteration, or {\bf(D)}
 has orbit in a disjoint cycle of attracting Fatou components. Note
 that every component of Type {\bf D} is contained in a complete copy
 of the Mandelbrot set. (See \cite{IK}, \cite{IKR}.)}
\setcounter{footnote}{10}

As an example, the cell structure for the curve $~\overline\cS_4/\I~$ has:
\smallskip

\centerline{14 vertices (= ideal points),\quad  68 edges, and \quad
44 faces.}\smallskip

\noindent See Figure~\ref{f-s4plot} for a diagram illustrating part
of this cell complex, including 22 of the 44 faces and 12 of the
14 vertices. (Most of these vertices are shown in multiple copies.) The only
vertices which are not included are the complex conjugates of $~r~$ and $~k$,
which will be denoted by $~r^*~$ and $~k^*$.~ (All of the other vertices
are self-conjugate.)
One can obtain the full cell complex from this partial diagram in three
steps as follows:

\begin{figure}[ht]
\centerline{\psfig{figure=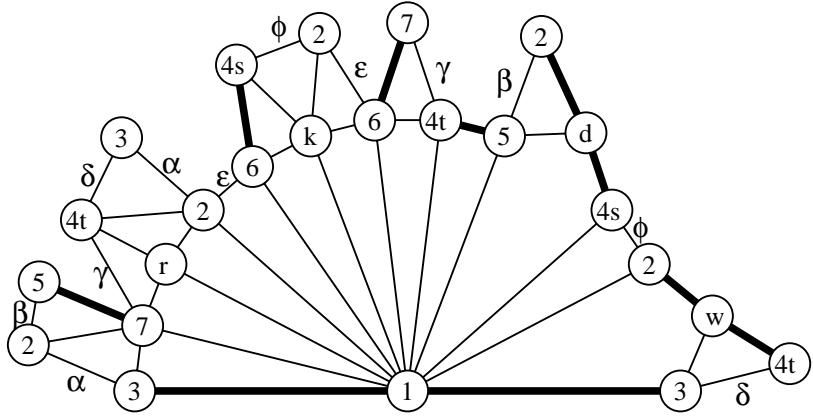,height=2.5in}}
\caption{\textsl{Schematic picture of the cell structure for half of
$~\overline\cS_4/\I$.~ Here each vertex  is represented by a circle containing
a brief label for the corresponding escape region. In particular, the kneading
sequences $~\overline{0010}\,,~\overline{0100}\,,~ \ldots\,,~\overline{1110}~$
 in lexicographical order are labeled
 simply by the integers $1,\,2\,,\ldots,\,7$.\label{f-s4plot}}}
\end{figure} 

\begin{itemize}
\item[$\bullet$] Take a second ``complex conjugate'' copy of this diagram, with
the two interior vertices of the complex conjugate diagram labeled
as $~r^*$~ and $~k^*$~ respectively.

\item[$\bullet$] Glue these two diagrams together along the ten emphasized 
edges. (Each emphasized edge represents a one real parameter family of mappings
$~F_{a,v}~$, where $~a~$ and $~v~$ are either both real or both pure
imaginary. Compare \S\ref{s-real}.) The
 resulting surface will have twelve temporary boundary
circles, labeled in pairs by Greek letters.

\item[$\bullet$] Now glue each of these boundary circles onto the other
 boundary circle
which has the same label, matching the two so that the result will be
an oriented surface. (In most cases, this means gluing an edge from the
original diagram onto an edge from the complex
conjugate diagram. However, there is
one exception: Each edge with label $\beta$ must be glued to the other $\beta$
edge from the same diagram.)
\end{itemize}
\smallskip

\noindent The result of this construction will be a closed oriented
surface.  In particular, the given edge identifications will force the
required  vertex identifications.\smallskip

Evidently complex conjugation operates
 as an orientation reversing involution of this surface.
The fixed points of this involution (made up out of the ten emphasized edges
in Figure~\ref{f-s4plot}) consist of two circles, which
represent all possible ``real'' cubic maps in $~\cS_4/\I$. (Compare
\S\ref{s-real}.)  One of these circles has
just two vertices  $\,1\,$ and $\,3\,$,~ while the other has eight vertices,
listed in cyclic order as
$~~2\,,~{\rm w}\,,~4_{\rm t}\,,~5\,,~7\,,~6\,,~4_{\rm s}\,,~{\rm d}\,$.

Note that this combinatorial structure is closely related to the dynamics:\break
At the center of each triangular face, three escape regions come
together, either
at a parabolic point  or around a hyperbolic component of Type A or D.
Similarly, at the center of each quadrilateral face,
 four escape regions come together around a hyperbolic component of Type B.
\footnote{For higher periods, there may be more possibilities. Perhaps
three or more escape regions may come together at a critically
finite point, or around a hyperbolic component of Type C.}
\smallskip

Given this cell structure, we can easily compute the Euler characteristic
$$ \chi(\overline\cS_4/\I) ~=~ 
 14-68+44 ~=~ -10\,.$$
Setting $~\chi=2-2g$,~ this corresponds to a genus of
 $~g(\overline\cS_4/\I)\,=\,6$.
The curve $~\overline\cS_4~$ itself is a two-fold branched
covering of $~\overline\cS_4/\I$, ~
ramified over the eight of the vertices (namely, those
labeled ~d, k, r, w, 5, and 7 in Figure \ref{f-s4plot}, as well as the
conjugate vertices $k^*$ and $r^*$).
Thus the induced cell subdivision has Euler characteristic
$$ \chi\big(\,\overline\cS_4\big)~=~2-2g(\overline\cS_4)~=~
 (2\cdot 14-8)- 2\cdot 68+2\cdot 44~=~ -28\,,$$
corresponding to genus 15. (Compare Table~\ref{tab-euler}.)


\def\c{{\mathfrak C}}
\setcounter{table}{0}
\section{Real Cubic Maps}\label{s-real}

It will be convenient to say that the cubic map $F_{a,\,v}$ is
\textbf{\textit{~real~}} if $a$ and $v$ are real, so that $~F_{a,v}(\R)=\R$,~
and \textbf{\textit{~pure imaginary~}} if $~a~$ and $~v~$ are pure
imaginary,~ so that $F_{a,v}(i\R)=i\R$.\smallskip

(In fact,  every ``pure imaginary'' map $F_{ia,\,iv}$ is conjugate to
a map
$$ x~\mapsto~ F_{ia,\,iv}(ix)/i~=~ -x^3\;+\;3a^2\,x\;+\;(-2a^3+v) $$
with real coefficients. However,
 the coefficient of the cubic term in this conjugate map
is $-1$ rather than $+1$. This leads
to a drastic difference in real dynamical behavior.)
\smallskip

Note that the space of monic centered complex cubic polynomials
has two commuting anti-holomorphic involutions: 
namely complex conjugation, which will be denoted by
$$\c^+\,:\,F_{a,v}\mapsto F_{\,\overline a,\,\overline v}\,,$$
and the composition of $\c^+$ with the standard involution $\I$,
which will be denoted by
$$ \c^-\,:\,F_{a,\,v}~\mapsto~F_{-\overline a,\,-\overline v}\,.$$
The fixed points of $\c^+$ are precisely
the ``real'' polynomials, and the fixed points of $\c^-$
are the ``pure imaginary'' polynomials. (There is only one common fixed point,
namely the polynomial $F_{0,0}(z)=z^3$.)
\smallskip

Now let us specialize to the period $p$ curve $\cS_p$, or to its smooth
compactification $\overline\cS_p$. 
 Each of the involutions $\c^\pm$
operates smoothly on the complex curves $\cS_p\subset\overline\cS_p$,
 reversing orientation.
The fixed point set is a real analytic curve which will
be denoted by $\cS_p^{\R\pm}\subset\overline\cS_p^{\R\pm}$.
It is easy to check that the curves $\cS_p^{\R+}$ and $\cS_p^{\R-}$
are disjoint, provided that $p>1$. However the compactifications
$\overline\cS_p^{\R+}$ and $\overline\cS_p^{\R-}$ intersect transversally
at a number of ideal points.\smallskip

 We will prove the following.

\begin{theo}\label{t-rc} 
Each connected component of $\cS_p^{\R\pm}$ is a
path joining two distinct ideal points $($or joining the unique ideal
point to itself in the special case $p=1)$. Furthermore,  each component
intersects one and only one hyperbolic component of Type {\rm A} or
{\rm B}, and contains the center point of this component.
 Thus there is a one-to-one correspondence between connected
 components of $\cS_p^{\R\pm}$ and\break $\c^\pm$-invariant hyperbolic
 components of Type {\rm A} or {\rm  B}.
\end{theo}

In fact,  any map which is the center of a hyperbolic component of Type
A is strictly monotone. In the $\c^+$-invariant case, it is monotone
increasing, and hence cannot have any periodic point of period $p>1$.
In the $\c^-$-invariant case it is monotone decreasing, and hence cannot
have any point of period $p>2$. Thus, if we exclude the very special cases
$\cS_1^{\R\pm}$ and $\cS_2^{\R-}$, then we must have a hyperbolic component of
 Type B.


\begin{figure}[ht]
\centering
\subfigure[\textsl{A component of $\cS_4^{\R+}$} \label{f-B31a}]{\makebox[0pt][l]{\includegraphics[height=2in]{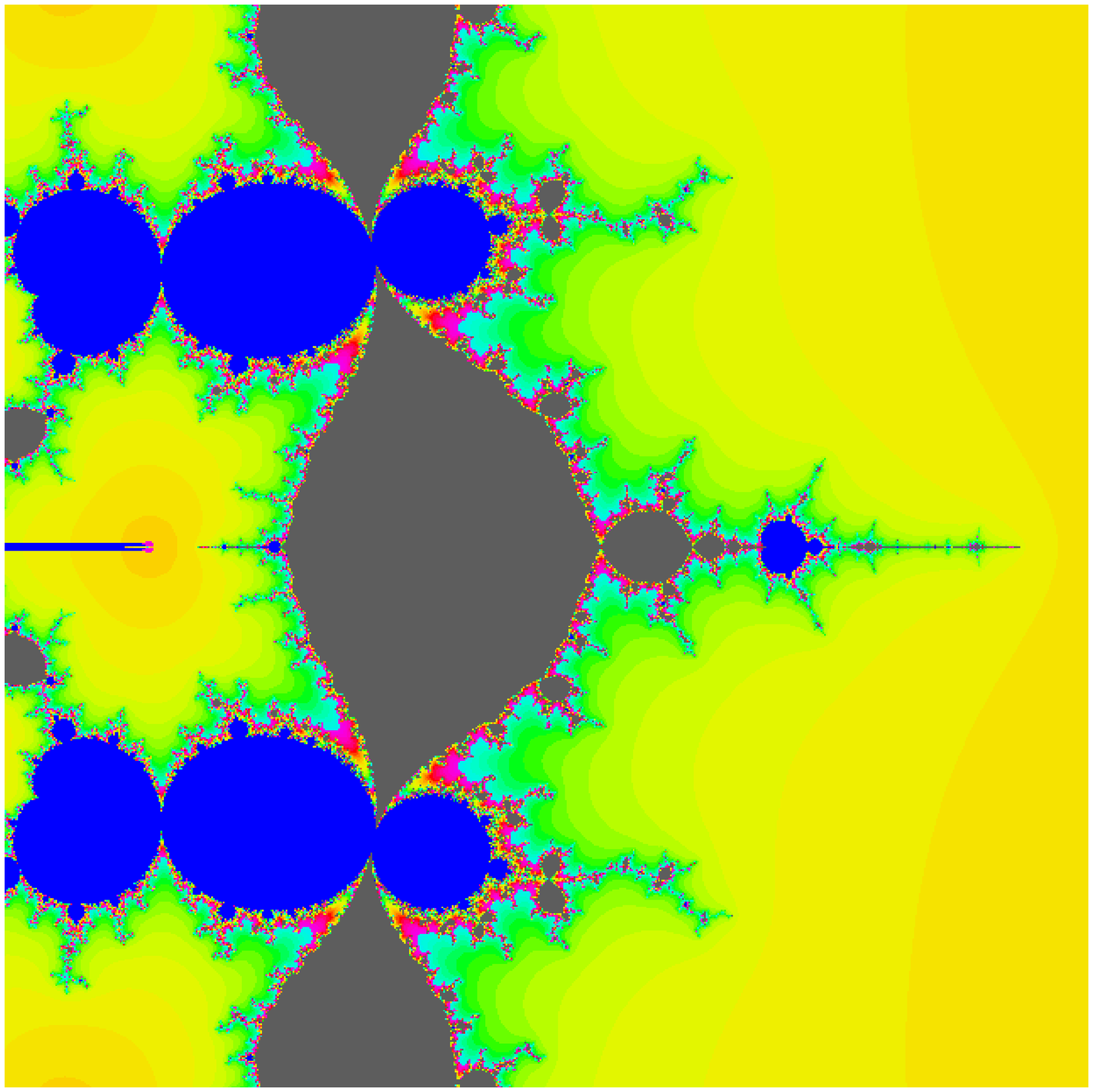}}%
  \includegraphics[height=2in]{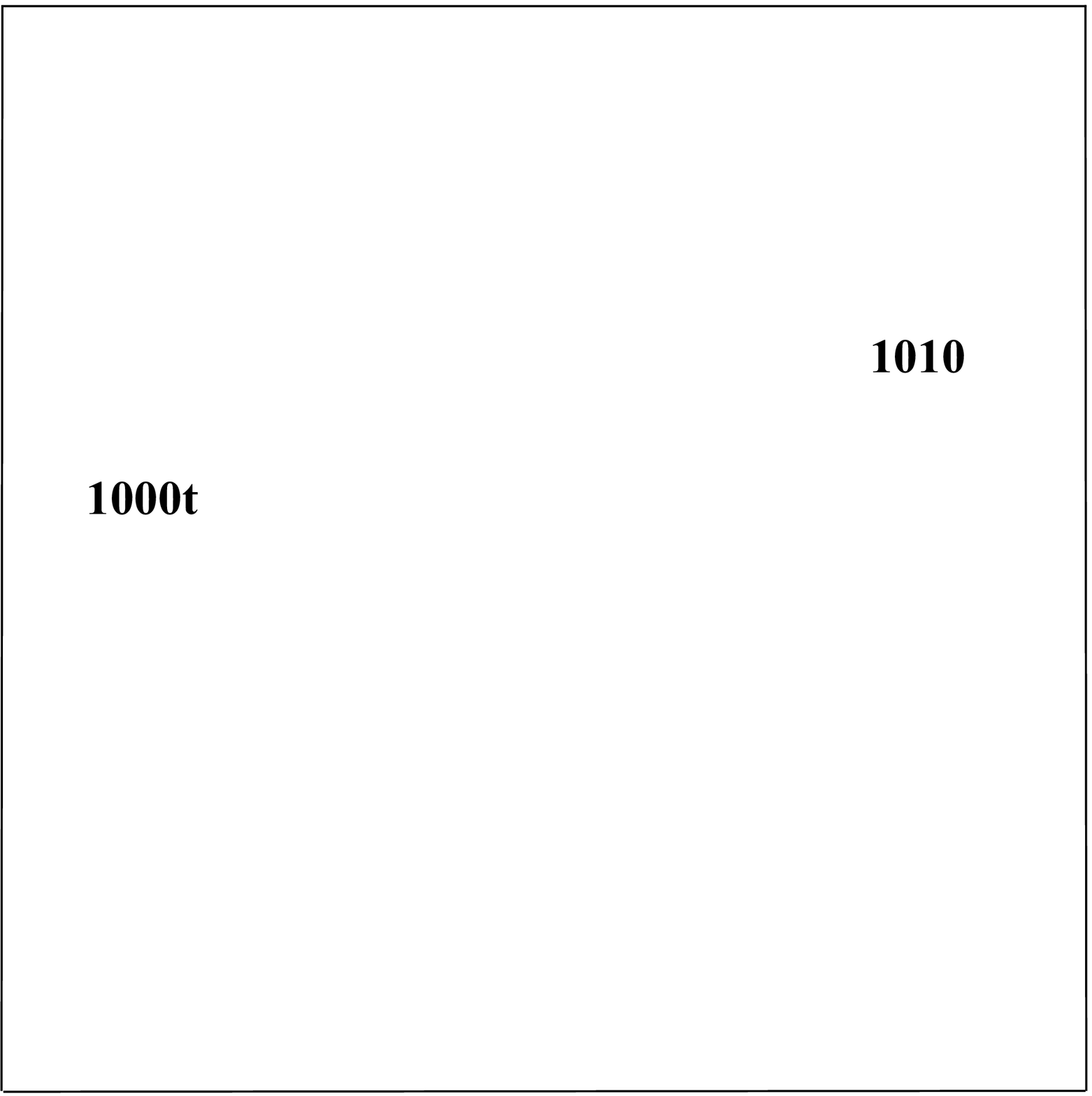}}\qquad
\subfigure[\textsl{A component of $\cS_4^{\R-}$} \label{f-B31b}]{\makebox[0pt][l]{\includegraphics[height=2in]{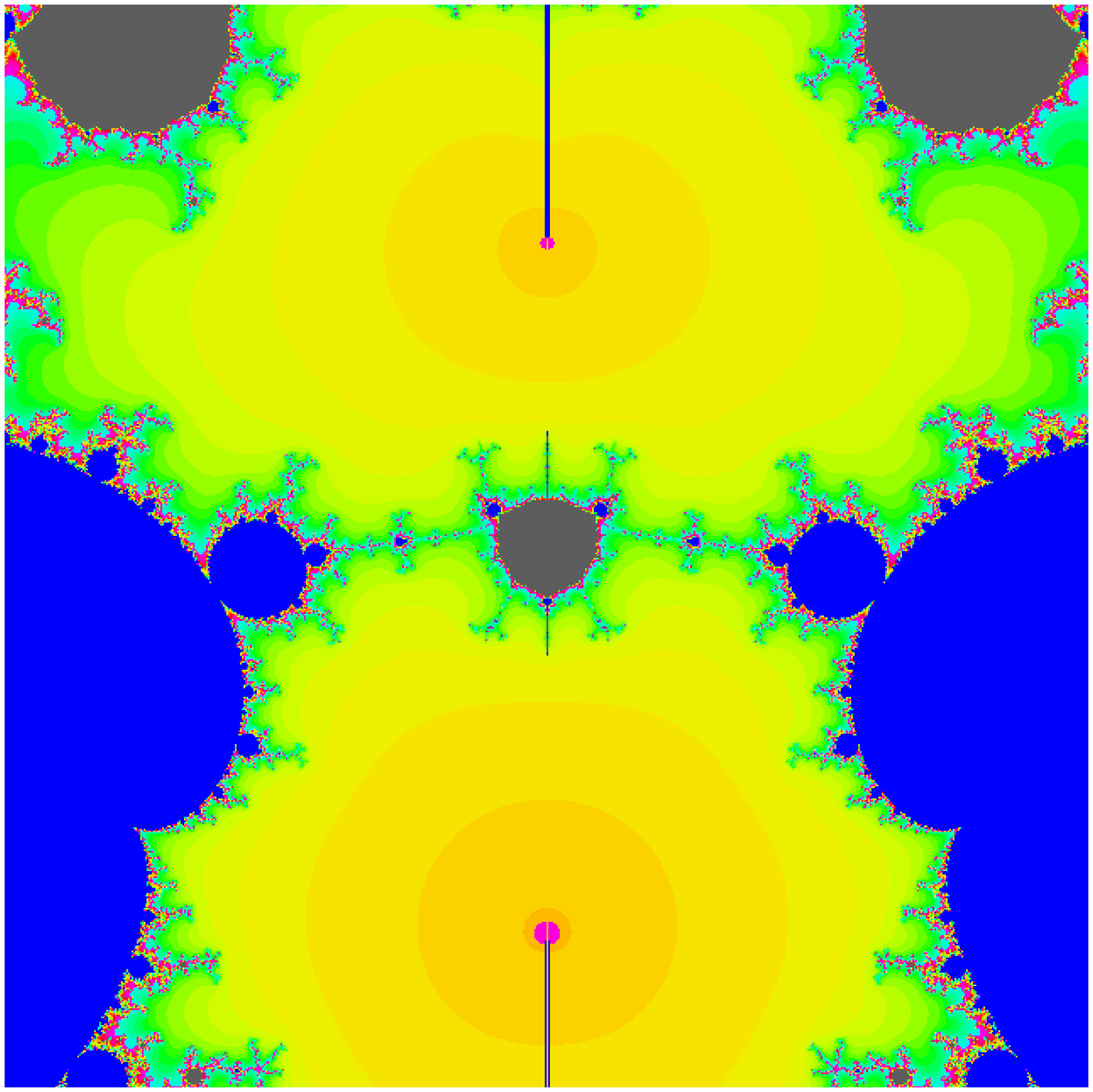}}%
  \includegraphics[height=2in]{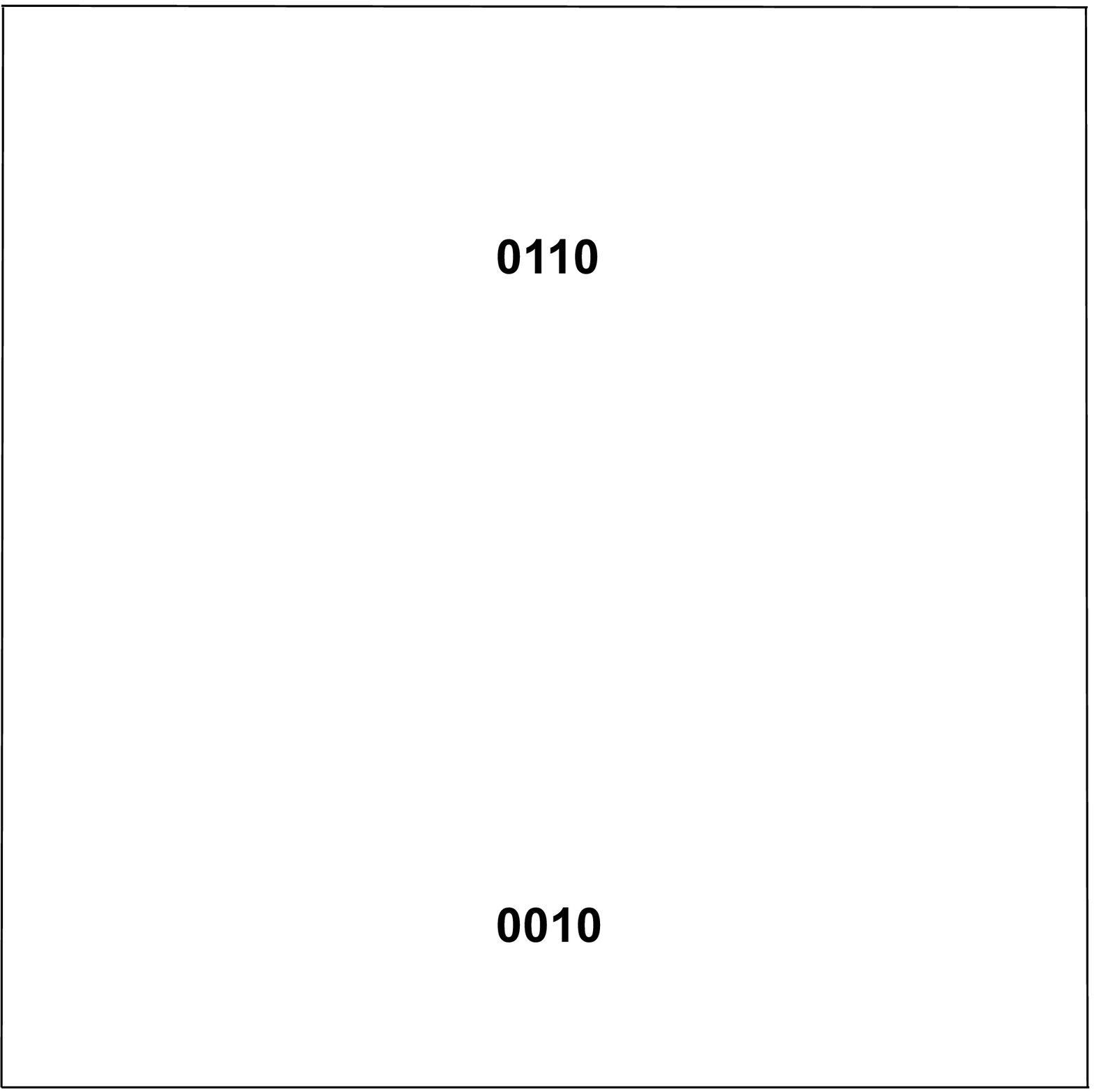}}
\caption{\textsl{Examples of
 $\c^\pm$ invariant components of Type B.  In Case (a), 
the associated component of $\cS_4^{\R+}$ is a horizontal line leading
 from the $\overline{1000}$t ideal point on the left
to the $\overline{1010}$ ideal point which is off-screen on the right. (There
 is a slit to the left, since the $\overline{1000}$t escape region is locally 
a 3-fold branched covering of the $t$-parameter plane.)  In Case (b),
 the associated component of $\cS_4^{\R-}$ is a vertical line leading
from the $\overline{0010}$ point below to the $\overline{0110}$ point above.
 $($These pictures
are centered at $\,(a,v)~=(-.79263,\,1.19929)$ and  $(a,v)~=~(0.743709\,i,\,
0.029712\,i)$
 respectively.\,$)$} \label{f-B31}}
\end{figure} 

\newcommand\boldf{{\mathbf f}}

Each component of Type B can be conveniently labelled by its
Hubbard tree, together with a specification of the marked critical
point $+a$ and the free critical point $-a$. This tree
 can be described as a piecewise linear map $~\boldf~$
from the interval $[0,\,p-1]$ to itself which takes integers to
integers and is linear on each subinterval $[j,\,j+1]$.
 It will be
convenient to extend $~\boldf~$ to a piecewise linear map from $\R$ to
itself which has constant  slope $>1$  outside of the interval
$[0,\,p-1]$ in the $\c^+$ invariant case, or constant slope $<-1$
in the $\c^-$ invariant case. A map $~\boldf~$ of this type
corresponds to a pair of components of Type B if and only if:

\begin{itemize}
\item[$\bullet$] $\boldf$ permutes the integers $\{0,\,1,\,\ldots,\,p-1\}$
cyclically, and

\item[$\bullet$] $\boldf$ is \textit{\textbf{~bimodal\,}}, i.e.,
 with one local minimum and one local maximum.
\end{itemize}

\def\a{{\bf a}} 

\noindent
To specify a unique component of Type B, we must also specify which
of the two local extrema is to be the marked critical point. Let $\a_0$
be the marked critical point of $\boldf$ (or in other 
words the integer corresponding to $a_0$), and let $\widehat\a_0$
be the free critical point for $\boldf$. We will also use the notations
 $\a_j=\boldf^{\circ j}(\a_0)$ and $\widehat\a_j
=\boldf^{\circ j}(\widehat\a_0)$.

\begin{figure}[ht]
\centering
\subfigure[with~~ $\widehat \a_0=\a_2$\label{f-rgf0}]{\psfig{figure=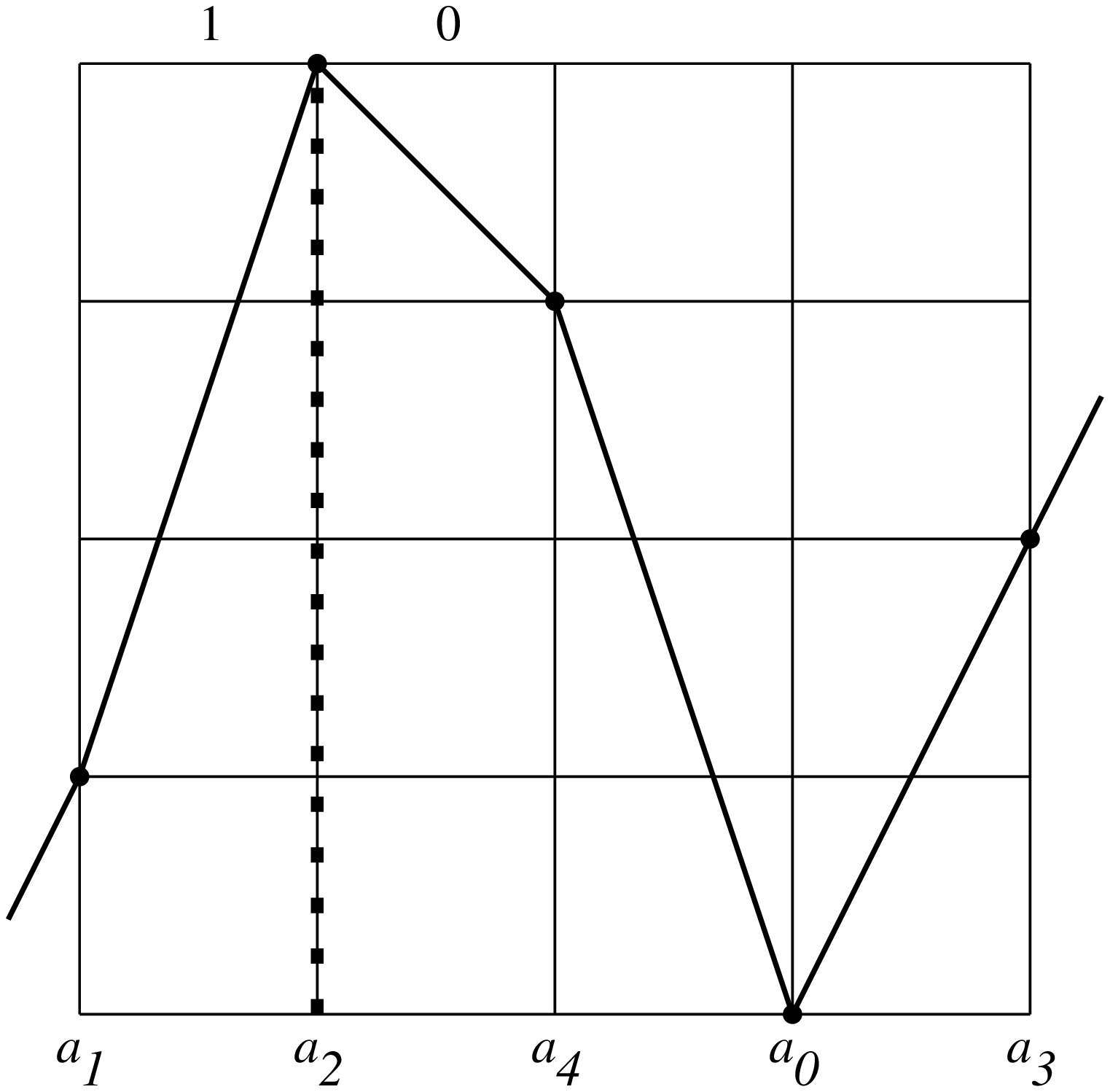,height=2.in}}\qquad
\subfigure[with~~ $\a_2<\widehat \a_0$\label{f-rgf1}]{\psfig{figure=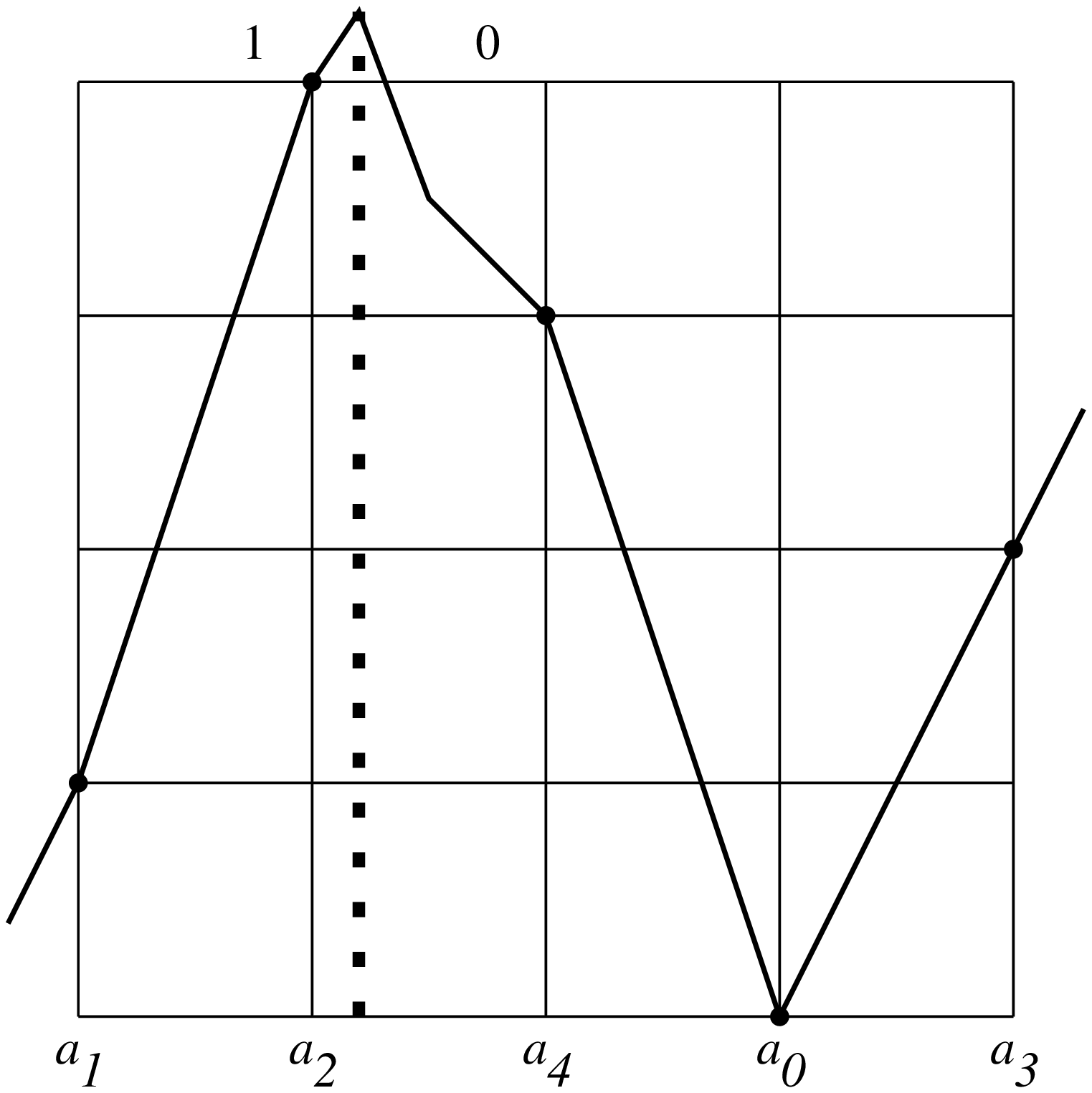,height=2.in}}
\caption{\label{gra}}
\end{figure}

As an example, Figure~\ref{f-rgf0}
 shows the graph for a typical Type B center in
$\cS_5^{\R+}$. The marked critical orbit is labelled as
$\a_0\mapsto \a_1\mapsto \cdots\mapsto \a_{p-1}$, and
 the free critical point is labelled 
by a vertical dotted line. 
The information in this graph can be summarized by the inequalities
$$ \a_1<\a_2=\widehat \a_0<\a_4<\a_0<\a_3 \,.$$
(Since this example is $\c^+$ invariant, this is completely equivalent
to the set of inequalities $a_1<a_2=\widehat a_0<a_4<a_0<a_3$. ~
In the $\c^-$ invariant case, with the $a_j$ pure imaginary,
 we would write corresponding inequalities
for the real numbers $a_j/i$.)\smallskip

{\bf Definition.} 
The \textit{\textbf{~kneading sequence~}} associated with
such a bimodal map with periodic marked critical point
 is a sequence of $p$ symbols which
can be defined as follows. The free critical point $\widehat\a_0$ divides the
real line into two halves. {\it Assign 
the \textit{\textbf{address}} {\bf 0} to every number on the same side as
the marked critical point $\a_0$, and {\bf 1} 
 to every point on the opposite side.\/}\footnote{As in Definition
 \ref{d-kn}, our kneading sequence describes the location of orbit points
only in comparison with the free critical
point. This should not be confused with the kneading sequences in \cite{MTh}
which describe location with respect to {\it both\/} critical points.}
(As an example,
 in Figures \ref{f-rgf0} and \ref{f-rgf1}, everything to the right of the
dotted line has address 0, and everything to the left has address 1.)
 Furthermore, define the address of the
point $\widehat\a_0$ itself to be the symbol \ding{72}. The kneading
sequence is then defined to be the sequence of addresses of the orbit
 points $\a_1,\,\a_2,\,\ldots,\,\a_p$ (where $\a_p=\a_0$). Thus
for a generic map, with $\widehat\a_0$ disjoint from the marked critical
orbit, we obtain a sequence of zeros and ones. However, for
 the center point of any component of Type B (or in
degenerate cases of Type A),
the sequence will contain exactly one \ding{72}.

As an example, for Figure \ref{f-rgf0}
  the kneading sequence is {\bf 1\ding{72}000}.
However, if we move the free critical point and the free critical value
a little to the right, as indicated in Figure \ref{f-rgf1},
 leaving the rest of the 
graph unchanged,  then the ordering becomes
$$ \a_1<\a_2<\widehat \a_0<\a_4<\a_0< \a_3\,,$$
and the kneading sequence will change to {\bf 11000},
replacing \ding{72} with {\bf 1}. Similarly, if we move the free critical
 point a little to the left (still moving the free critical value to the
right), we can replace \ding{72} with {\bf 0}. {\it A similar
argument applies to any component of Type B which is
$~\c^\pm$-invariant.}

\begin{figure} [ht]
\centerline{\psfig{figure=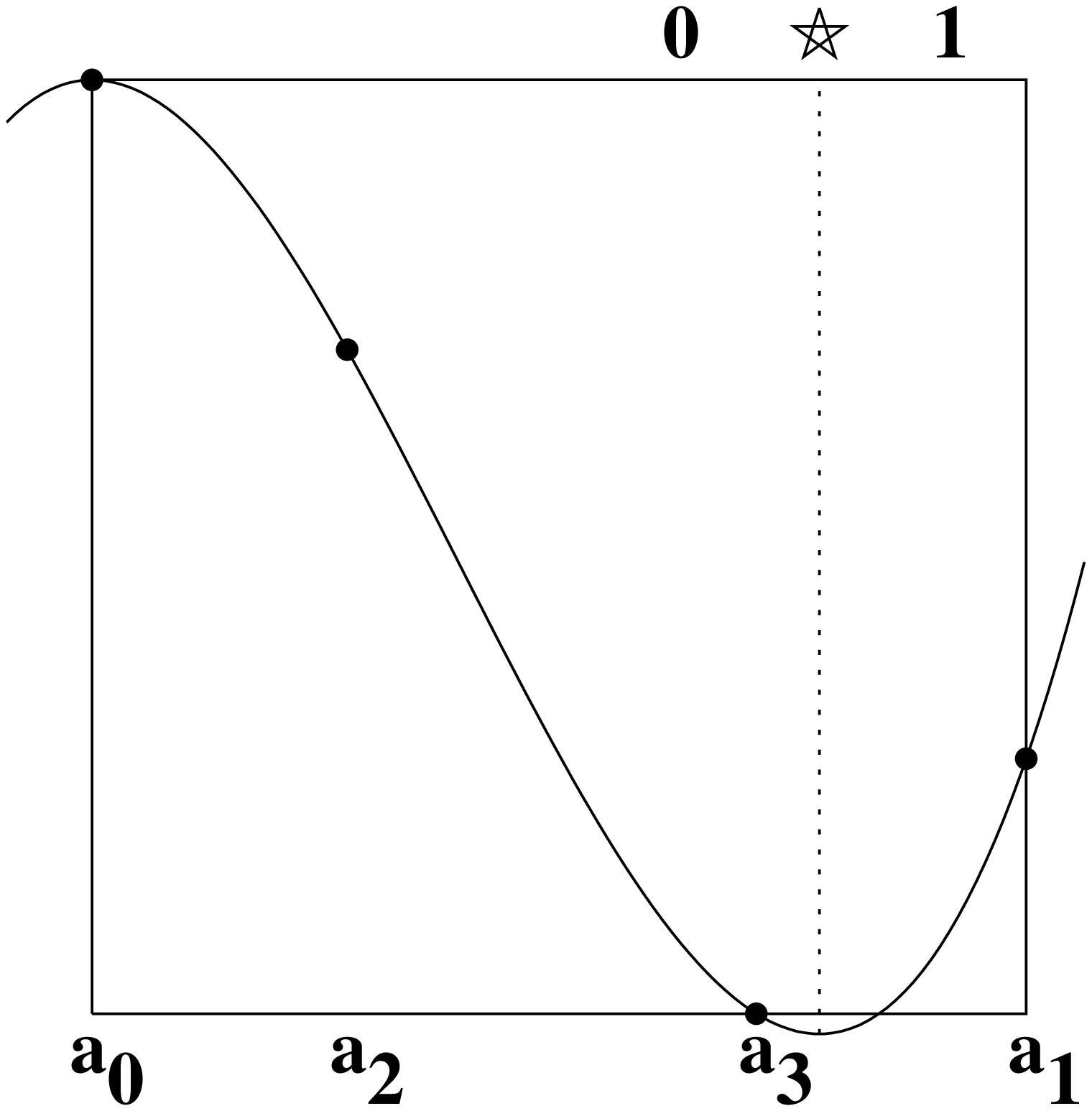,height=1.7in}\quad
\psfig{figure=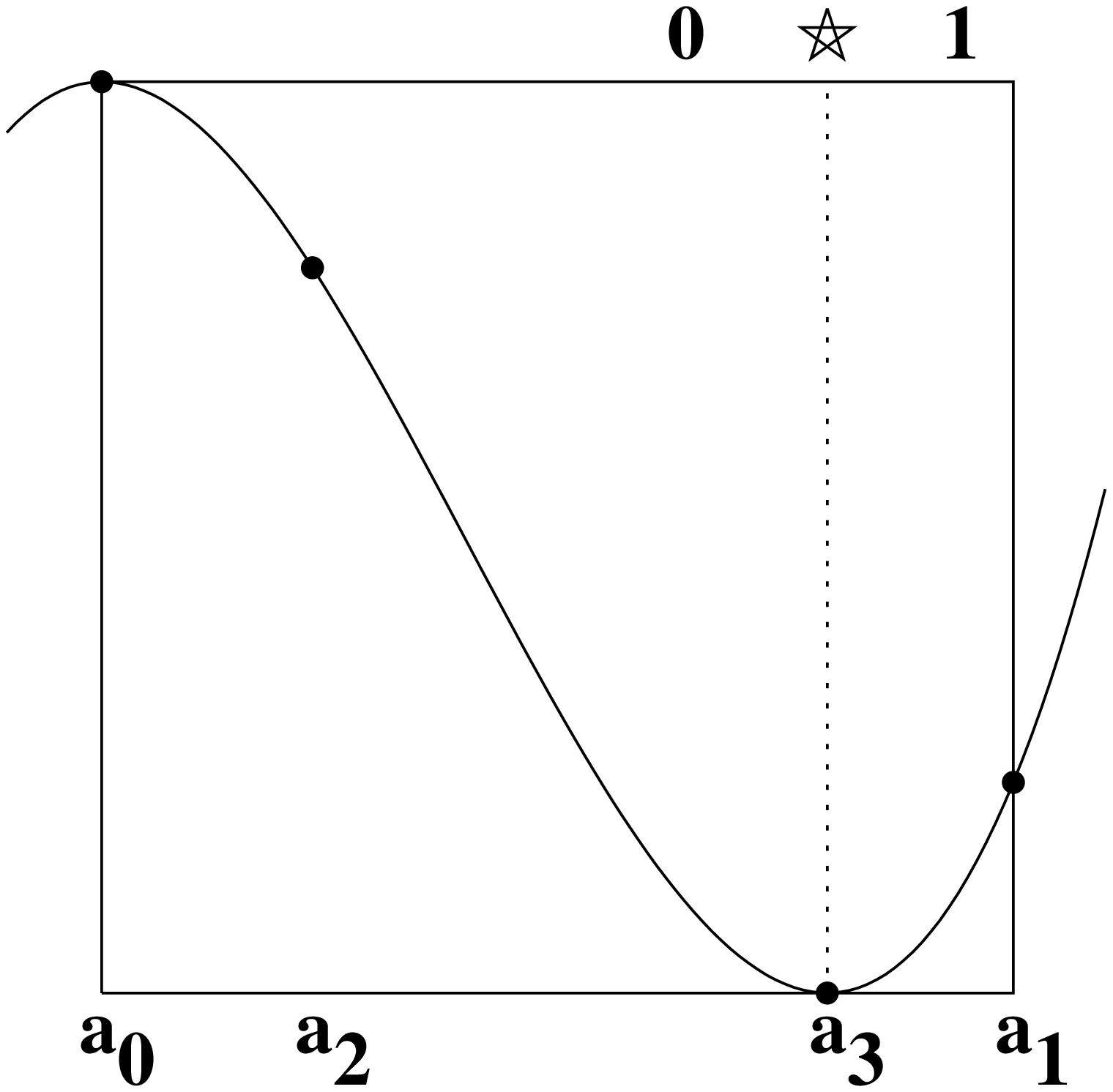,height=1.7in}\quad
\psfig{figure=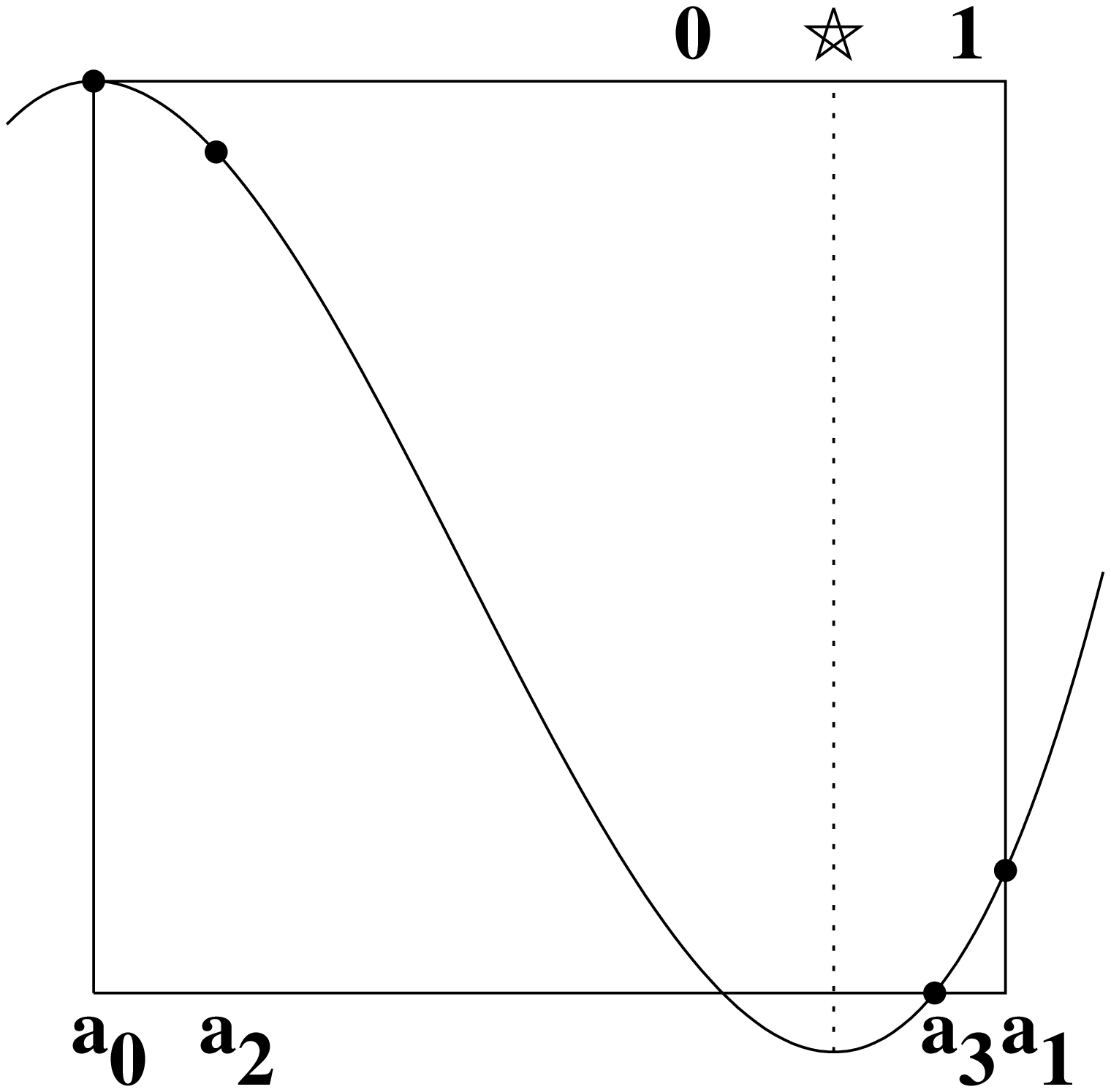,height=1.7in}}
\caption{{\textsl Graphs illustrating the transition from
kneading sequence {\bf 1000} through~ {\bf{10\ding{72}0}}~ to ~{\bf 1010}.}
\label{f-p4grafs}}
\end{figure}

\medskip
 As a period 4 example, the relevant kneading sequences as we move from left
to right near the middle of Figure \ref{f-B31}(a)
 are illustrated in Figure \ref{f-p4grafs}.
\medskip

The proof of Theorem \ref{t-rc} will depend on three lemmas. 
\smallskip

\begin{lem}[{\bf Heckman}]\label{l-heck}
No connected component of $\cS_p^{\R\pm}$ can be
a simple closed curve.\footnote{However, \textbf{every} component of the
compactified locus $\overline\cS_p^{\R+}$ or
 $\overline\cS_p^{\R-}$ is a simple closed curve. Compare Figure
 \ref{f-3tor} where the horizontal simple closed curve
 and the vertical simple closed curve intersect at two ideal points, labeled
~{\bf a}~ and ~{\bf 110}.}
\end{lem}

The proof is quite difficult. See \cite{He}.\qed

\begin{figure}
\centerline{\psfig{figure=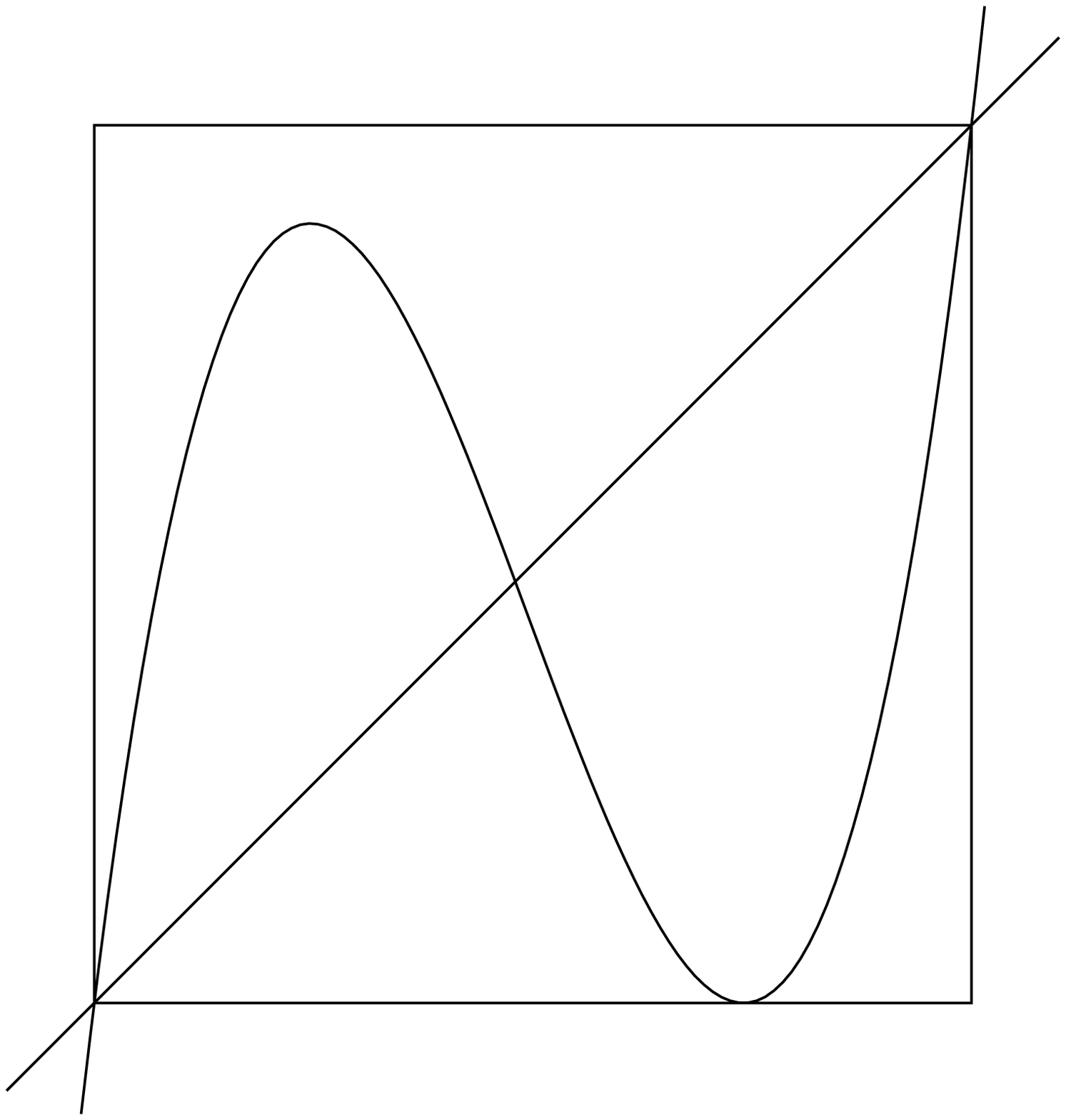,height=1.7in}\qquad\qquad
\psfig{figure=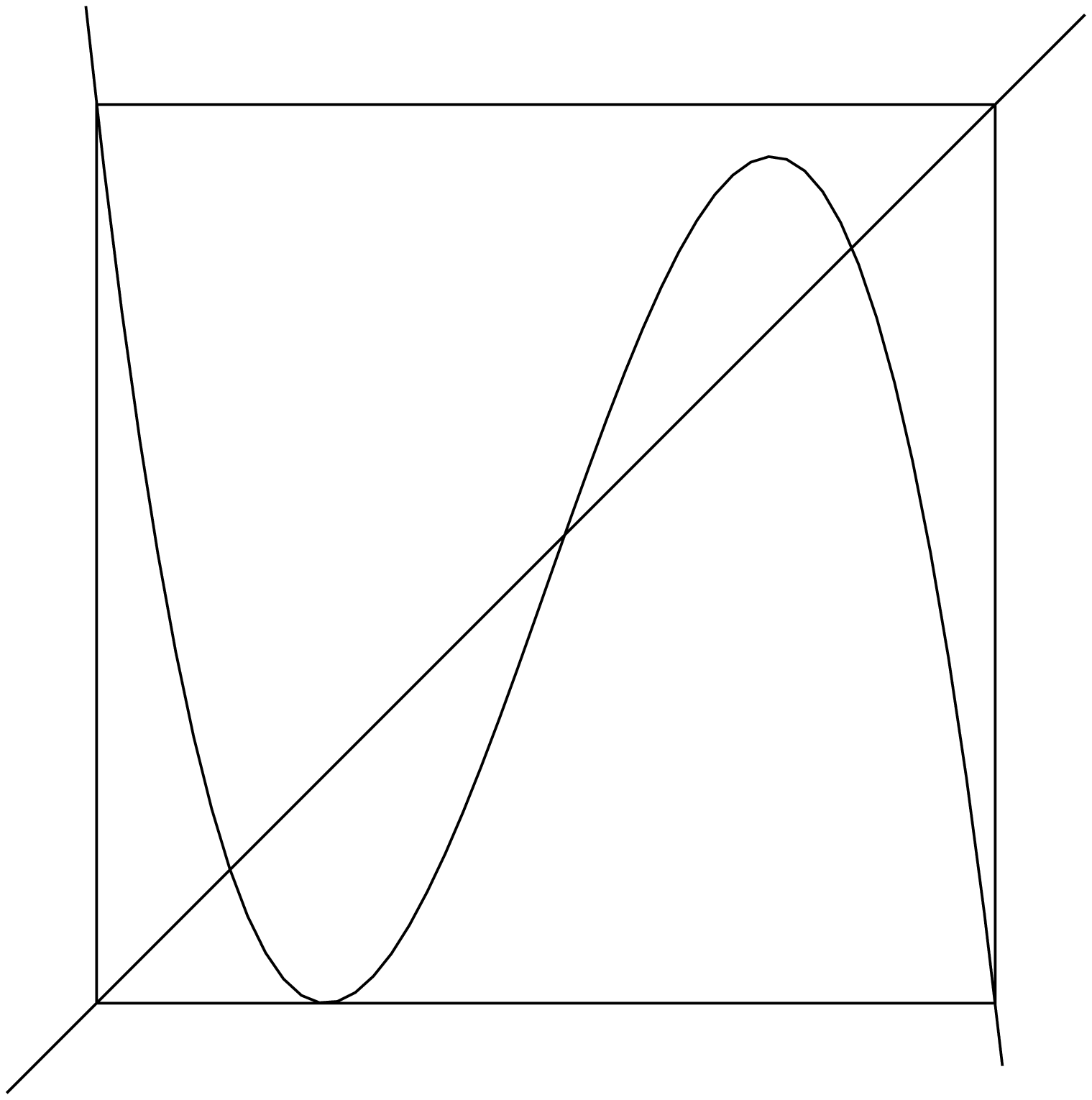,height=1.7in}}
\caption{\textsl Graphs of maps on the boundary of the connectedness
locus within $\cS_p^{\R+}~$ respectively $~\cS_p^{\R-}$.\label{f-bdgrafs}}
\end{figure}

\begin{lem}[{\bf Milnor and Tresser}]\label{l-mtr}
Each connected component of $\cS_p^{\R\pm}$ intersected with
the connectedness locus is homeomorphic to a non-degenerate closed interval
of real numbers. The two endpoints of this interval can be characterized as
those maps in the interval for which the free critical value $~F(-a)~$ is
an extreme point of the real or pure imaginary filled Julia set $($which
 is itself a closed interval---compare Figure \ref{f-bdgrafs}$)$. It follows
that this free critical value is a fixed point in the $\c^+$-invariant case,
and a period two point in the $\c^-$-invariant case.
\end{lem}

{\bf Proof.} See \cite[Lemma 7.2]{MTr}. (The proof makes use of Lemma
 \ref{l-heck}.)\qed

\begin{lem}\label{l-not-cc}
For $~p>1~$ no two maps belonging to the same connected component of
 $~\cS_p^{\R\pm}~$ can be conformally conjugate to each other.
It follows that the connected components $P\subset \cS_p^{\R\pm}$ occur in
pairs, with each $P$  disjoint from $\I(P)$.\end{lem}

{\bf Proof.} Note that two monic
centered cubic polynomial maps are conformally conjugate if and only
if they are either equal to each other, or carried one to the other
by the involution $\I$. But if $\I$ interchanged two points in the
same component $P$, then it would map $P$ onto itself, with a fixed point
in the middle. This is impossible for $p>1$,
 since the only fixed point of $~\I(F_{a,v})=F_{-a,\,-v}$\break
is the point $F_{0,0}\in\cS_1$.\qed\medskip

{\bf  Proof of Theorem~\ref{t-rc}.} If we use the canonical local parameter 
$t$, then
the involution $\c^\pm$ will transform $~t~$ to $~\pm\overline t+({\rm
constant})$.~
Hence the invariant curve $\cS_p^{\R+}$ is represented by a horizontal
line $\Im(t)={\rm constant}$, 
while $\cS_p^{\R-}$ is represented by a vertical line $\Re(t)={\rm
constant}$. (Compare Figure \ref{f-B31}.)
 Clearly such a line can not have any limit point within
the open surface $\cS_p$. Since the connected component cannot be a
simple closed curve, it can only be an infinite path $P$  which leads
from one ideal point to another. 
 We will see that these two ideal
points must always be distinct when $p>1$, 
 since they have different kneading invariants.
\smallskip

We next show that each component $P\subset\cS_p^{\R\pm}$
contains at least one Type A or B center,
that is, at least one map for which the free critical point
$~\widehat a_0$ is equal to some point $~a_j~$  on the marked critical
orbit. Otherwise, if the points
$~a_0,\,\ldots,\,a_{p-1}~$ and $~\widehat a_0~$ were pairwise distinct
for all maps $~F\in P$,~ then their ordering would be the same
for all $~F\in P$.\smallskip

We know that the left hand end
of $P$ is contained in some escape component $\cE_1$, and that the right
 hand end is contained in some escape component $\cE_2$. Each of these
components must map to itself under $\c^\pm$. It follows from Lemma 3.1
that each escape region is conformally diffeomorphic to
$\C\ssm\overline{\mathbb D}$. Furthermore, if the escape region is
$\c^\pm$-invariant, then we can choose this conformal diffeomorphism so
that $\c^\pm$ corresponds to the involution $~z\mapsto\pm\overline z$
of $\C\ssm\overline{\mathbb D}$.~
Thus the fixed point set will correspond to the real
axis in the $\c^+$ case, or the imaginary axis in the $\c^-$ case. It follows
that the intersection $\cE_h\cap P$  is always a path leading from the ideal
point to some boundary point of $\cE_h$. It then follows from Lemma \ref{l-mtr}
that the remaining portion of $P$, its intersection with the connectedness
locus, is a closed interval of real numbers, bounded by two post-critically
 finite maps. Hence each of these boundary maps  is uniquely
determined, up to conjugacy, by its Hubbard tree. If there were no
center points of Type A or B within this path $P$, then the two endpoints
of $P$ intersected with the connectness locus would have to have
 isomorphic Hubbard
 trees. But this is impossible by Lemma \ref{l-not-cc}; hence $P$ must contain
a Type A or B center.\smallskip

Finally, we must prove that there cannot be two such Type A or B centers on
the path $P$. Otherwise, as $F$ varies over the subinterval $P'$
joining two consecutive centers,  the $p+1$ points
 $~a_0,\, a_1,\,\ldots,\, a_{p-1}~$ and $\widehat a_0$ would vary smoothly
and remain distinct. But, among all possible variations, there is only one way
of converging to a Type A or B center. Suppose for example that one of the
endpoints of $P'$ is topologically conjugate to the map of Figure~\ref{f-rgf0},
 and suppose
as in Figure~\ref{f-rgf1} that we have $a_1<a_2<\widehat a_0<a<a_4<a_0<a_3$ for
 the points
$F\in P'$. Since only one of the two endpoints of the interval $[a_2,\,a_4]$
is a local maximum, we can converge to a Type B center only by letting
$\widehat a_0$ converge to $a_2$. This shows that the two endpoints of the
interval $P'$ must have the same Hubbard tree. But this is impossible
by Lemma~\ref{l-not-cc}. The argument in the general case is completely analogous.

Finally, it follows from the discussion of Figures~\ref{f-rgf0} and 
\ref{f-rgf1} that
the two endpoints of the path $P$ must be ideal points with kneading sequences
which differ in exactly one entry. Hence these two endpoints are certainly
distinct.\QED\medskip


\begin{coro} A connected component in $\cS_p^{\R\pm}$ is uniquely
determined by the ordering of the points
 $\a_0,\,\a_1,\,\ldots,\,\a_{p-1}$ in the marked critical orbit.\end{coro}

The proof is immediate. If the graph associated with this ordering has an
interior local maximum or minimum  at a point $\a_j$ with $j\ne 0$,
then we must set $~\widehat \a_0=\a_j$.~ Otherwise, $\widehat \a_0$ must be one
of the two endpoints, and only one of the two is compatible with $\c^\pm$
invariance. In either case, the Hubbard tree and hence the associated
Type B component is uniquely determined.\QED\medskip

\begin{rem} Notice in Figure~\ref{f-B31a}  [and \ref{f-B31b}] that there is 
also a small component
of Type D on the real [or imaginary] axis, immediately to the left
[or below] the big
component. {\it In fact, for periods $p\ge 3$,
 every $\c^\pm$-invariant Type B component is associated with a
$\c^{\pm}$-invariant Type D component which is immediately adjacent to it,
consisting of maps for which both attracting orbits have period $p$.
The common boundary point is a parabolic map, with parabolic orbit of period
$p$ on the boundary of the cycle of attracting Fatou components.}
(Compare \cite[Lemma 7.1]{MTr}.)

The associated Hubbard tree (or graph) can easily be constructed from
the piecewise linear map $~f~$ for the Type B center.
Suppose for example, as in Figure~\ref{f-rgf0}, that $f$ has a local maximum
$\widehat\a_1=f(\widehat\a_0)$ at the free critical point. 
 Consider a one parameter family
of modified piecewise linear functions, as in Figure \ref{f-rgf1}, as follows.
Choose some small $\epsilon>0$, and replace the original free critical point
 $\widehat\a_0=\a_k$ by $~\widehat\a_0=\a_k\pm\epsilon\,$;~
 but set
the free critical value $f_t(\widehat\a_0)$ equal to $f(\a_k)+t=\a_{k+1}+t$
for small $t>0$, leaving all of the marked orbit points $\a_j$ unchanged.
Then as $t$ increases, for each $j\le p$ the
forward image $f_t^{\circ j}(\widehat\a_0)$
 will move linearly to the left or
right. Choose the sign in the equation $\widehat\a_0=
\a_k\pm\epsilon$ so that $f_t^{\circ p}(\widehat a_0)$ will move towards
$\widehat\a_0$. For some first value  $t_0$, we must have the precise
 equality $f_{t_0}^{\circ p}(\widehat \a_0)=\widehat\a_0$.
Then $f_{t_0}$ represents
the required Type D center, with both periodic orbits periodic of period $p$.
\end{rem}
\medskip

\vspace{1.cm}

\noindent
Araceli Bonifant; Mathematics Department, University of Rhode Island, Kingston,
 R.I.,  02881. {\it e-mail:} bonifant@math.uri.edu
\bigskip

\noindent
Jan Kiwi; Facultad de Matem\'aticas, Pontificia Universidad Cat\'olica, Casilla
 306, Correo 22, Santiago de Chile, Chile. {\it e-mail:} jkiwi@puc.cl
\bigskip

\noindent           
John Milnor;  Institute for Mathematical Sciences, Stony Brook University,
Stony Brook, NY. 11794-3660. {\it e-mail:} jack@math.sunysb.edu

\end{document}